\documentclass[11pt, leqno]{article}

\usepackage{amsmath, amssymb}
\usepackage[hang,flushmargin]{footmisc}
\numberwithin{equation}{section}

\begin{document}
\baselineskip 1.7em

\centerline{\bf\Large Tensor Products, Positive Linear Operators,}

\centerline{\bf\Large and Delay-Differential Equations}
\bigskip

\centerline{\small John Mallet-Paret\footnote{Partially supported by NSF DMS-0500674 and by The Center for Nonlinear Analysis at Rutgers University.} and Roger D. Nussbaum\footnote{Partially supported by NSF DMS-0701171 and by The Lefschetz Center for Dynamical Systems at Brown University.}}
\medskip

\begin{abstract}
\noindent We develop the theory of compound functional differential equations,
which are tensor and exterior products of linear functional differential
equations. Of particular interest is the equation
$$
\dot x(t)=-\alpha(t)x(t)-\beta(t)x(t-1)
$$
with a single delay, where the delay coefficient is of one sign, say
$\delta\beta(t)\ge 0$ with $\delta\in\{-1,1\}$. Positivity properties are
studied, with the result that if $(-1)^k=\delta$ then the $k$-fold
exterior product of the above system generates a linear process which
is positive with respect to a certain cone in the phase space.
Additionally, if the coefficients $\alpha(t)$ and $\beta(t)$
are periodic of the same period, and $\beta(t)$ satisfies a uniform
sign condition, then there is an infinite set of
Floquet multipliers which are complete with respect to an associated
lap number. Finally, the concept of $u_0$-positivity of the exterior product
is investigated when $\beta(t)$ satisfies a uniform sign condition.
\\
{\bf Date:} May 1, 2012.
\end{abstract}

\newcommand{\tint}{\sim\hspace*{-13pt}\int}
\newcommand{\eps}{\varepsilon}
\newcommand{\alp}{\alpha}
\newcommand{\bet}{\beta}
\newcommand{\gam}{\gamma}
\newcommand{\Gam}{\Gamma}
\newcommand{\del}{\delta}
\newcommand{\kap}{\kappa}
\newcommand{\lam}{\lambda}
\newcommand{\sig}{\sigma}
\newcommand{\Sig}{\Sigma}
\newcommand{\ome}{\omega}
\newcommand{\Ome}{\Omega}
\newcommand{\vph}{\varphi}
\newcommand{\calh}{{\cal H}}
\newcommand{\calk}{{\cal K}}
\newcommand{\cale}{{\cal E}}
\newcommand{\calg}{{\cal G}}
\newcommand{\calj}{{\cal J}}
\newcommand{\calb}{{\cal B}}
\newcommand{\cali}{{\cal I}}
\newcommand{\call}{{\cal L}}
\newcommand{\cals}{{\cal S}}
\newcommand{\calu}{{\cal U}}
\renewcommand{\AA}{{\bf A}}
\newcommand{\BB}{{\bf B}}
\newcommand{\PP}{{\bf P}}
\newcommand{\QQ}{{\bf Q}}
\newcommand{\NN}{{\bf N}}
\newcommand{\IR}{{\bf R}}
\newcommand{\IC}{{\bf C}}
\newcommand{\UU}{{\bf U}}
\newcommand{\WW}{{\bf W}}
\newcommand{\XX}{{\bf X}}
\newcommand{\XXI}{\mbox{\boldmath $\Xi$}}
\newcommand{\ZZ}{{\bf Z}}
\newcommand{\iop}[2]{{I^{#1}_{#2}}}
\newcommand{\jop}[2]{{J^{#1}_{#2}}}
\newcommand{\dist}{\mathop{\rm dist}}
\newcommand{\spann}{\mathop{\rm span}}
\newcommand{\spec}{\mathop{\rm spec}}
\newcommand{\espec}{\mathop{\rm ess~spec}}
\newcommand{\kr}{\mathop{\rm ker}}
\newcommand{\ran}{\mathop{\rm ran}}
\newcommand{\sgn}{\mathop{\rm sgn}}
\newcommand{\codim}{\mathop{\rm codim}}
\newcommand{\re}{\mathop{\rm Re}}
\newcommand{\sch}{\mathop{\rm sc}}
\newcommand{\tr}{|\hspace*{-1.5pt}|\hspace*{-1.5pt}|}
\newcommand{\wt}{\widetilde}
\newcommand{\wh}{\widehat}
\newcommand{\ovrl}{\overline}
\newcommand{\ds}{\displaystyle}
\newcommand{\bxx}{\begin{it}}
\newcommand{\exx}{\end{it}}
\newcommand{\xx}{\it}
\newcommand{\vv}{\vspace{.2in}}
\newcommand{\bb}{\hspace*{-.08in}}
\newcommand{\fp}{\hspace{1mm}\rule{2mm}{2mm}}
\newcommand{\be}[1]{\begin{equation}}
\newcommand{\ee}[1]{\label{#1}\end{equation}}
\newcommand{\prodform}[3]{\newlength{#1}\settowidth{#1}{$\ds{\prod_{#2}\hspace{2em}}$}\makebox[.5#1][l]{$\ds{\prod_{{#2}}}$}#3}
\newcommand{\prodformtwo}[4]{\prodform{#1}{{\scriptstyle #2}\atop{\scriptstyle #3 \rule{0pt}{8pt}}}{#4}}
\newlength{\mlength}
\settowidth{\mlength}{$M$}
\newcommand{\mhat}{\makebox[\mlength][c]{$\wh{\makebox[0pt][c]{$M\,$}}$}}
\newcommand{\mtilde}{\makebox[\mlength][c]{$\wt{\makebox[0pt][c]{$M\,$}}$}}
\newlength{\wwlength}
\settowidth{\wwlength}{$\WW$}
\newcommand{\wwhat}{\makebox[\wwlength][c]{$\wh{\makebox[0pt][c]{$\WW$}}$}}
\newcommand{\wwtilde}{\makebox[\wwlength][c]{$\wt{\makebox[0pt][c]{$\WW$}}$}}

\section{Introduction}

In this paper we develop the theory of compound functional differential equations.
This is in the spirit of compound ordinary differential equations and dynamical
systems as developed by J.S.~Muldowney~\cite{muld} and Q.~Wang~\cite{wang}. Broadly,
this topic concerns tensor products and exterior products of linear nonautonomous evolutionary systems.
With this approach, criteria for the nonexistence of periodic solutions, and also for the asymptotic
stability of periodic solutions, were obtained in~\cite{muld} for nonlinear ordinary differential equations in $\IR^n$.
(See also Y.~Li and J.S.~Muldowney~\cite{limul}.) Extensions of these results to classes
of nonlinear partial differential equations were given in~\cite{wang}. Overall, their results generalize
classical two-dimensional results of Poincar\'e, Bendixson, and Dulac.

A central feature of our investigation involves
positivity issues connected with compound systems arising from the delay-differential equation
\be{int1}
\dot x(t)=-\alp(t)x(t)-\bet(t)x(t-1).
\ee{int1}
Note that equation~\eqref{int1} has a single delay. Typically we assume a signed feedback, that is,
$\beta(t)$ is of constant sign, either positive or negative, for almost every $t$.

Delay-differential equations have been studied for at least 200 years.
While some of the early work had its origins in certain types of
geometric problems and number theory, much of the impetus for the development
of the theory came from studies of viscoelasticity
and population dynamics (notably by Volterra; see~\cite{vo1}, \cite{vo2}, and~\cite{vo3}),
and from control theory (see, for example, Bellman~\cite{bell}).
More recent work has involved models
from a wide variety of scientific fields, including
nonlinear optics~\cite{gi}, \cite{ik1}, \cite{ik2},
economics~\cite{bel}, \cite{mackey1},
biology~\cite{gl}, \cite{mackey2}, and as well population dynamics~\cite{ho}, \cite{smi}.
Classic references for much of the fundamental theory are
the books of J.K.~Hale and S.M.~Verduyn Lunel~\cite{hale} and of
O.~Diekmann, S.A.~van Gils, S.M.~Verduyn Lunel, and H.-O.~Walther~\cite{di}.
For additional material see~\cite{ba}, \cite{ha2}, \cite{ha3}, \cite{hen1}, \cite{hen2}, and~\cite{wu}.
Equations such as~\eqref{int1} can arise as linearizations around solutions
of nonlinear equations. Two such very classic yet still challenging nonlinear equations are Wright's Equation
$$
\dot x(t)=-\alp x(t-1)(1+x(t)),
$$
and the Mackey-Glass Equation
$$
\dot x(t)=-\alp_1x(t)+\alp_2f(x(t-1)),
$$
where $f(x)=x/(1+x^n)$. See~\cite{wr} and~\cite{mackey2}, as well as the references in~\cite{hale} for further information
on such equations.

Abstractly, a linear (evolutionary) process $U(t,\tau):X\to X$ on a Banach space $X$ is a collection
of bounded linear operators $U(t,\tau)$, for $t\ge\tau$, for which $U(\tau,\tau)=I$ and
$U(t,\sig)U(\sig,\tau)=U(t,\tau)$ whenever $t\ge\sig\ge\tau$, with $U(t,\tau)x$ varying continuously
in $(t,\tau)$ for each fixed $x$. Linear processes occur as solution maps of a wide variety of
nonautonomous linear equations, including of course the finite-dimensional case $\dot x=A(t)x$
of an ordinary differential equation.
In the case of the delay-differential equation~\eqref{int1}
the underlying Banach space is $X=C([-1,0])$, and we assume that $\alpha:\IR\to\IR$ and $\beta:\IR\to\IR$
are locally integrable functions.

Given an abstract linear process $U(t,\tau)$ as above, and given an integer $m\ge 1$,
one obtains the so-called compound processes
$$
\UU(t,\tau)=U(t,\tau)^{\otimes m},\qquad
\WW(t,\tau)=U(t,\tau)^{\wedge m},
$$
by taking the $m$-fold tensor product and $m$-fold exterior (wedge) product, respectively, of the operator $U(t,\tau)$.
These compound processes are themselves linear processes on the tensor and exterior products
$X^{\otimes m}$ and $X^{\wedge m}$ of the space $X$.

In general, if $X_j$ for $1\le j\le m$ are
Banach spaces, then one may consider the tensor product
$$
X_0=X_1\otimes X_2\otimes\cdots\otimes X_m
$$
of these spaces. For infinite dimensional spaces, there are typically many natural but inequivalent
norms for $X_0$ arising from the norms on the $X_j$. For our purposes, the so-called
injective cross norm is the suitable choice of a norm for $X_0$, and it is used throughout this paper.
In a natural way, if $A_j$ are bounded linear operators on $X_j$ for $1\le j\le m$, one
obtains a bounded linear operator $A_0=A_1\otimes A_2\otimes\cdots\otimes A_m$ on $X_0$.
In the case all $X_j=X$ are the same space one writes $X_0=X^{\otimes m}$, and also
$A_0=A^{\otimes m}$ if all operators $A_j=A$ are the same. The exterior
product $X^{\wedge m}\subseteq X^{\otimes m}$ is the subspace of $X^{\otimes m}$ consisting of elements
which satisfy a certain anti-symmetry property in a fashion analogous to the well-known
finite-dimensional case. The subspace $X^{\wedge m}$ is invariant for the operator $A^{\otimes m}$,
and one denotes by $A^{\wedge m}=A^{\otimes m}|_{X^{\wedge m}}$ the restriction of this operator
to this subspace.

A key point connected with tensor and exterior products of operators is the behavior of their
spectra. Suppose that the essential spectral radius $\rho(A)$ of $A$ satisfies $\rho(A)=0$; in particular, this is
the case if either $A$ or some power $A^n$ of $A$ is compact, and this in turn is the case
for the solution operators $U(t,\tau)$ associated to the delay equation~\eqref{int1} for $t>\tau$. Then for any $m\ge 1$ the
spectrum of $A^{\wedge m}$ consists of all products $\lam_1\lam_2\cdots\lam_m$
where $\{\lam_k\}_{k=1}^\infty$ are elements of the spectrum of $A$ (there may be only finitely many),
and where the number of repetitions of a given $\lam_k$ in this product cannot exceed
the multiplicity of $\lam_k$ as an element of the spectrum of $A$.
(For a precise statement of this result, including a formula for the multiplicity
of $\lam_1\lam_2\cdots\lam_m$ as an element of $\spec(A^{\wedge m})$ and a description
of the eigenspace, see Proposition~2.3 below, along with Corollary~2.2.)
As was observed by Muldowney~\cite{muld}, this
fact has ramifications for the stability of periodic orbits of nonlinear systems as it is
applied to the Floquet analysis of the linearized system; see the final remark of Section~3 below.

A surprising aspect of compound systems for equation~\eqref{int1} relates to positivity properties
when the feedback coefficient $\bet(t)$ is of constant sign. A main result of this paper is Theorem~4.1,
the Positivity Theorem. This states that if $(-1)^m\bet(t)\ge 0$ almost everywhere, then for any
$t$ and $\tau$ with $t\ge\tau$, the operator
$\WW(t,\tau)=U(t,\tau)^{\wedge m}$ associated to equation~\eqref{int1} is a positive operator
with respect to an appropriate cone in $X^{\wedge m}=C([-1,0])^{\wedge m}$,
specifically, to the cone $K_m$ given in~\eqref{km}, with~\eqref{tri}.
Note that the integer $m$ is restricted to a given parity
type, namely, $m$ is even in the case of so-called negative feedback ($\beta(t)\ge 0$), while
$m$ is odd in the case of positive feedback ($\beta(t)\le 0$).

(A cone $K$ in a Banach space $X$ is defined to be a closed convex set for which $\sigma u\in K$
whenever $u\in K$ and $\sigma\ge 0$, and for which $u,-u\in K$ only if $u=0$. A linear operator
$A$ on $X$ is called positive if $A$ maps $K$ into $K$. For some general references, see
the original paper of Kre\u\i n and Rutman~\cite{okr},
as well as~\cite{bonsall}, \cite{deim}, \cite{31}, \cite{nw} and~\cite{schwlf}.
More general results can be found in~\cite{kr}.)

If the coefficients in~\eqref{int1} are periodic, say if
$\alp(t+\gam)\equiv\alp(t)$ and $\bet(t+\gam)\equiv\bet(t)$
hold identically for some $\gam>0$, and if also for the second
coefficient there is a uniform lower bound of the form
$(-1)^m\bet(t)\ge(-1)^m\bet_0>0$ for almost every $t$ and some integer $m$, then computable lower
bounds on the norms $|\lam|$ of the Floquet multipliers (characteristic multipliers)
can be obtained; see Corollary~5.4 below. More precisely, the set of (nonzero) Floquet multipliers $\{\lam_k\}_{k=1}^\infty$
is a countably infinite set. If it is ordered so that
$$
|\lam_1|\ge|\lam_2|\ge|\lam_3|\ge\cdots
$$
with repetitions according to algebraic multiplicity, then an explicit lower bound for each $|\lam_k|$
can be given. Further, the strict inequality
$$
|\lam_k|>|\lam_{k+1}|
$$
holds for each $k$ for which $k-m$ is even
(that is, for a particular parity class, odd or even, of $k$), and for each $k$ of either parity there is
associated to $\lambda_k$ a so-called lap number $\calj(|\lambda_k|)$ which measures the rate of oscillation
of the eigensolutions. A precise statement of these results is given in Theorem~5.1, which confirms conjectures
and extends results of G.R.~Sell and one of the authors~\cite{sell} which provided partial information on
Floquet multipliers. Let us also mention earlier results of S.-N.~Chow
and H.-O.~Walther~\cite{chow} in this spirit (see also~\cite{diek}), as well as
related results of H.L.~Smith and one of the authors~\cite{smith} for cyclic systems of ordinary differential equations.

For $m$ of the above parity, the monodromy operator $U(\tau+\gam,\tau)^{\wedge m}$ possesses a positive eigenvector,
here meaning positive with respect to the cone $K_m$. In concrete terms, as given in Proposition~5.2 below,
this means there exist independent solutions $x^j(t)$, for $1\le j\le m$, each one of which is a linear combination
of Floquet solutions corresponding to Floquet multipliers $\lambda_k$ for $1\le k\le m$, such that
\be{dt}
\det\left(\begin{array}{ccc}
x^1(t+\theta_1) &\bb \cdots &\bb x^1(t+\theta_m)\\
\vdots &\bb &\bb \vdots\\
x^m(t+\theta_1) &\bb \cdots &\bb x^m(t+\theta_m)
\end{array}\right)\ge 0,
\qquad
-1\le\theta_1\le\theta_2\le\cdots\le\theta_m\le 0,
\ee{dt}
for every $t\in\IR$.
Indeed, the left-hand side of~\eqref{dt} is simply the leading eigensolution of the positive process
$\WW(t,\tau)=U(t,\tau)^{\wedge m}$.

Interestingly, to our knowledge the positivity of $\WW(t,\tau)$ and the inequality~\eqref{dt} are new
results even in the
case of the constant coefficient equation
$$
\dot x(t)=-\alpha_0x(t)-\beta_0x(t-1).
$$
Here the solutions $x^j(t)$ in~\eqref{dt} are simply linear combinations of $e^{\zeta_kt}$
with $1\le k\le m$, where $\{\zeta_k\}_{k=1}^\infty$ are the roots of the characteristic equation $\zeta=-\alpha_0-\beta_0e^{-\zeta}$
ordered so that $\re\zeta_1\ge\re\zeta_2\ge\re\zeta_3\ge\cdots$.
Even for this constant coefficient case, we know of no elementary proof of~\eqref{dt}.

A key technique used in proving these results is the so-called lap number, which measures
the number of sign changes per delay interval, or equivalently, the rate of oscillation, of a solution $x(t)$.
Lap numbers were used in~\cite{morse} to obtain a Morse decomposition of the attractor for certain nonlinear
feedback systems, and also in~\cite{sell}, \cite{pb}, and~\cite{smith}
to prove Poincar\'e-Bendixson Theorems for high dimensional cyclic systems with a monotone feedback.
The essential idea of a lap number goes back to A.D.~Myschkis~\cite{myschkis}.

An alternate approach to these results, which also provides refinements of them, is via
$u_0$-positivity, a useful idea introduced by
Krasnosel'ski\u\i; see~\cite{26} and~\cite{27}, and also~\cite{eveson}. Generally, if
an operator $A$ is positive with respect to a cone $K$, and if $u_0\in K\setminus\{0\}$,
then $A$ is called $u_0$-positive if there exists a positive integer $k_0$ such that $A^ku\sim u_0$
for every $k\ge k_0$ and every $u\in K\setminus\{0\}$. Here one writes $u\sim v$ if there exist
positive quantities $C_1$ and $C_2$ such that $C_1v\le u\le C_2u$, with the ordering here with respect
to the cone. The point of $u_0$-positivity is that it allows one to conclude additional facts about
the spectrum and the eigenvectors. Specifically,
by a generalization~\cite{31} of the Kre\u\i n-Rutman Theorem it is known that if an operator $A$ is positive with
respect to a cone $K$ which is total (meaning the span of $K$ is dense in $X$),
and if the essential spectral radius of $A$ is less than the spectral
radius $r$ of $A$, then $\lambda=r$ is an eigenvalue of $A$ for which there exists a positive eigenvector,
that is, $Av=rv$ with $v\in K\setminus\{0\}$. However, it need not be the case that $\lambda=r$ is
a simple eigenvalue. But if it further holds that
$A$ is $u_0$-positive for some $u_0\in K\setminus\{0\}$ and $K$ is reproducing (meaning the span of $K$ equals $X$), then
all other spectral points of $A$ satisfy $|\lambda|<r$, the eigenvalue $\lambda=r$ is algebraically simple,
and no other nonzero spectral points have eigenvectors in the cone $K$. Further,
$v\sim u_0$ for the positive eigenvector $v$, and this provides upper and lower bounds for $v$.

The above conclusions are essentially the same as for a finite-dimensional positive operator which is primitive,
meaning some power of $A$ takes the nonzero cone into its interior. Thus $u_0$-positivity can be regarded
as an infinite-dimensional generalization of primitivity. Let us mention that in
our applications to delay-differential equations, the cone $K_m$, while reproducing, in fact has empty interior.

Proving $u_0$-positivity seems quite challenging, and indeed, while we establish it for the cases $m\le 3$ (see Theorem~7.1),
it remains open for $m\ge 4$. A consequence of $u_0$-positivity,
as described in Proposition~7.2, is a refinement of~\eqref{dt} to
$$
C_1u_0(\theta)\le
\det\left(\begin{array}{ccc}
x^1(t+\theta_1) &\bb \cdots &\bb x^1(t+\theta_m)\\
\vdots &\bb &\bb \vdots\\
x^m(t+\theta_1) &\bb \cdots &\bb x^m(t+\theta_m)
\end{array}\right)\le
C_2u_0(\theta)$$
for some positive and $t$-dependent $C_1$ and $C_2$, again for
$-1\le\theta_1\le\theta_2\le\cdots\le\theta_m\le 0$.
The function $u_0(\theta)=u_m(\theta)$, where $\theta=(\theta_1,\theta_2,\ldots,\theta_m)$, depends on $m$ and is given by
$$
u_2(\theta_1,\theta_2)=\theta_2-\theta_1,\qquad
u_3(\theta_1,\theta_2,\theta_3)=(\theta_2-\theta_1)(\theta_3-\theta_2)(\theta_3-\theta_1)(1+\theta_1-\theta_2),
$$
for $m\le 3$, which gives strict positivity when the $\theta_j$ are distinct.
For $m\ge 4$ we conjecture $u_0$-positivity with the function $u_m(\theta)$ in~\eqref{u0}.
The above determinant necessarily vanishes when $\theta_j=\theta_k$ for $j\ne k$, and so the factors $\theta_k-\theta_j$
in $u_m(\theta)$ are to be expected. However, the appearance of $1+\theta_1-\theta_2$
in $u_3(\theta)$ is surprising, indicating a zero of the determinant at
the point $(\theta_1,\theta_2,\theta_3)=(-1,0,0)$ which is of higher order than one might initially expect.
Indeed, the precise nature of such singularities
of these determinants is explored in depth in Section~8.

The outline of this paper is as follows. We review the basic foundations
of tensor products of Banach spaces and operators in Section~2. We specialize
this to general linear evolution processes in Section~3, and then further
to delay-differential equations in Section~4. One of our main results, the Positivity Theorem (Theorem~4.1),
is stated there. This theory is applied, using lap numbers as a basic tool, in Section~5, to obtain
fundamental results on Floquet solutions of periodic systems with a signed feedback; Theorem~5.1 is a main
result there. Section~6 is devoted to the proof
of Proposition~4.2, which provides the key result used to prove both the Positivity Theorem
and Theorem~5.1. Sections~7 and~8 are devoted to results on $u_0$-positivity. As mentioned, detailed
results for the cases $m\le 3$ are given, but some intriguing conjectures are raised for $m\ge 4$.

\section{Tensor Products of Banach Spaces}

In what follows we let $\call(X,Y)$ denote the space of bounded linear operators
between Banach spaces $X$ and $Y$. We also denote $\call(X)=\call(X,X)$.
For any operator $A\in\call(X)$, we let
$\spec(A)$ and $\espec(A)$ denote the spectrum and the essential spectrum of $A$,
and we let
$$
r(A)=\sup\{|\lam|\:|\:\lam\in\spec(A)\},\qquad
\rho(A)=\sup\{|\lam|\:|\:\lam\in\espec(A)\},
$$
denote the spectral radius and essential spectral radius of $A$, respectively.
We remark that while there are several inequivalent definitions of essential spectrum
(see, for example,~\cite{browder}, \cite{kato}, and~\cite{wolf}),
the quantity $\rho(A)$ is nevertheless the same for all of them.

To begin our discussion of tensor products, let $X$ and $Y$ be Banach spaces, and let $X\odot Y$
denote their algebraic tensor product.
Then $X\odot Y$
is the vector space
consisting of equivalence classes of elements of the form
\be{01}
z=\sum_{i=1}^na_i(x_i\otimes y_i)
\ee{01}
with $x_i\in X$ and $y_i\in Y$, and $a_i\in\IC$, under the equivalence
relation generated by all identities of the form
$$
\begin{array}{l}
(x+x')\otimes y=x\otimes y+x'\otimes y,\qquad
x\otimes(y+y')=x\otimes y+x\otimes y',\\
\\
a(x\otimes y)=(ax)\otimes y=x\otimes(ay),
\end{array}
$$
and only those identities.
There are various possible (generally inequivalent) norms for $X\odot Y$, among which are the so-called
{\bf cross norms}, namely norms for which $\|x\otimes y\|=\|x\|\|y\|$ holds
for every $x$ and $y$, and with the corresponding equation
holding with the dual norms.
In particular, the norm defined by
\be{02}
\|z\|=\sup_{(\xi,\eta)\in\calb}\bigg|\sum_{i=1}^na_i\xi(x_i)\eta(y_i)\bigg|,\qquad
\calb=\{(\xi,\eta)\in X^*\times Y^*\:|\:\|\xi\|=\|\eta\|=1\},
\ee{02}
for $z$ as in~\eqref{01},
where $X^*$ and $Y^*$ are the dual spaces to $X$ and $Y$,
is a cross norm, called the {\bf injective cross norm}. One easily checks that
$\|z\|$ is well-defined, that is, it is independent of the
representation~\eqref{01} of $z$, and that the formula~\eqref{02}
does indeed define a norm on $X\odot Y$.
Now define $X\otimes Y$ to be the Banach space which is the completion of
$X\odot Y$ with respect to this norm. Throughout this paper, we shall always
take the injective cross norm when considering tensor products of Banach spaces.

If $E$ and $G$ are two other Banach spaces, and $A\in\call(X,E)$
and $B\in\call(Y,G)$ are bounded linear operators, then one defines
the tensor product $A\otimes B$
of these operators by $(A\otimes B)(x\otimes y)=(Ax)\otimes(By)$,
and extends this by linearity first to $X\odot Y$,
and then continuously to all of $X\otimes Y$.
It is easily checked that
this construction determines a unique bounded linear operator
\be{03}
A\otimes B\in\call(X\otimes Y,E\otimes G),\qquad
\|A\otimes B\|=\|A\|\|B\|,
\ee{03}
with norm as indicated. One also sees that
\be{07}
(A_1\otimes B_1)(A_2\otimes B_2)=(A_1A_2)\otimes(B_1B_2)
\ee{07}
for operators defined on appropriate spaces.

If we have a direct sum decomposition $X=X_1\oplus X_2$ for $X$,
where $X_1,X_2\subseteq X$ are closed subspaces,
then there is a direct sum decomposition
\be{12}
X\otimes Y=(X_1\otimes Y)\oplus(X_2\otimes Y).
\ee{12}
We note that {\it a priori} there are two possible definitions
for $X_j\otimes Y$. Namely, $X_j\otimes Y$ can be defined either
(a)~directly, by considering $X_j$ as a Banach space in its own right
and taking the tensor product with $Y$, or (b)~by taking the closure
in $X\otimes Y$ of the subspace spanned by elements $x'\otimes y$
with $x'\in X_j$ and $y\in Y$. That these two constructions yield
the same result, namely isometric Banach spaces, follows from the identity
$$
\sup_{(\xi',\eta)\in\calb'}\bigg|\sum_{i=1}^n\xi'(x'_i)\eta(y_i)\bigg|
=\sup_{(\xi,\eta)\in\calb}\bigg|\sum_{i=1}^n\xi(x'_i)\eta(y_i)\bigg|,\qquad
\calb'=\{(\xi',\eta)\in X_j^*\times Y^*\:|\:\|\xi'\|=\|\eta\|=1\},
$$
with $x'_i\in X_j$ and $y\in Y$, and $\calb$ as in~\eqref{02}, which is
an immediate consequence of the Hahn-Banach Theorem.
(For the norm in $X_j$ we always take the norm inherited as a subspace of $X$.) 
In a similar fashion, if $Y=Y_1\oplus Y_2$ then
$X\otimes Y=(X\otimes Y_1)\oplus(X\otimes Y_2)$.

The above constructions extend in the obvious way to products and sums of
several Banach spaces. In particular, if $X$, $Y$, and $Z$ are Banach spaces,
then $(X\otimes Y)\otimes Z$ and $X\otimes(Y\otimes Z)$ are
naturally isometrically isomorphic.
If $X_j$ are Banach spaces for $1\le j\le m$ then
one can define $X_1\otimes X_2\otimes\cdots\otimes X_m$ in
a natural fashion, along with products $A_1\otimes A_2\otimes\cdots\otimes A_m$
of operators $A_j\in\call(X_j,E_j)$ where the $E_j$ are Banach spaces,
with the obvious generalization of~\eqref{03}. Similarly,~\eqref{12} generalizes
to the case of multiple summands and multiple factors.
We also note that for spaces $X_j$ of either finite or infinite dimension, we have that
$$
\dim(X_1\otimes X_2\otimes\cdots\otimes X_m)=\prod_{j=1}^m\dim X_j,
$$
with the convention that $0\times\infty=0$.

The following result on the spectra of tensor products is basic.
\vv

\noindent {\bf Theorem~2.1 (T.~Ichinose~\cite[Theorem~4{.}3]{ich1}; see also~\cite{ich2} and M.~Schechter~\cite{schechter}).} \bxx
Let $X_j$ be a Banach space and $A_j\in\call(X_j)$ for $1\le j\le m$.
Then
\be{14}
\spec(A_1\otimes A_2\otimes\cdots\otimes A_m)
=\{\lam_1\lam_2\cdots\lam_m\:|\:\lam_j\in\spec(A_j)\hbox{ for every }1\le j\le m\}
\ee{14}
for the spectrum of the tensor product.
\exx
\vv

The above theorem can be generalized to count multiplicities,
at least of isolated spectral points, as follows.
\vv

\noindent {\bf Corollary~2.2.} \bxx
Let $A_j$ and $X_j$ for $1\le j\le m$ be as in Theorem~2.1.
Denote $A_0=A_1\otimes A_2\otimes\cdots\otimes A_m$ and take
any $\lam_0\in\spec(A_0)$ with $\lam_0\ne 0$ for which $\lam_0$ is an isolated
point of $\spec(A_0)$. Then there are finitely many
distinct $m$-tuples
\be{13}
(\lam_1^k,\lam_2^k,\ldots,\lam_m^k)\in\IC^m,
\ee{13}
for $1\le k\le p$, such that
\be{17a}
\lam_0=\lam_1^k\lam_2^k\cdots\lam_m^k,\qquad
\lam_j^k\in\spec(A_j),
\ee{17a}
for $1\le j\le m$ and $1\le k\le p$.
Moreover, each such $\lam_j^k$ is an isolated point of $\spec(A_j)$.
Let $G_j^k\subseteq X_j$ denote the
spectral subspace of $A_j$ corresponding to $\lam_j^k$, let
$\nu_j^k=\dim G_j^k$, so $1\le\nu_j^k\le\infty$, and let
\be{18a}
\nu_0=\sum_{k=1}^p\nu_1^k\nu_2^k\cdots\nu_m^k.
\ee{18a}
Then the spectral subspace $G_0\subseteq X_0$ of $A_0$ corresponding to $\lam_0$ is given by
\be{18b}
G_0=\bigoplus_{k=1}^pG_1^k\otimes G_2^k\otimes\cdots\otimes G_m^k,
\ee{18b}
where $\dim G_0=\nu_0$, and where each subspace $G_1^k\otimes G_2^k\otimes\cdots\otimes G_m^k$
is invariant under $A_0$.
\exx
\vv

\noindent {\bf Remark.}
We are assuming that every possible $m$-tuple~\eqref{13} satisfying~\eqref{17a}
has been enumerated and thus occurs for some $k$.
Also, the $m$-tuples~\eqref{13} are geometrically distinct points in $\IC^m$,
with no repetitions for multiplicity as elements of a spectrum, that is,
$(\lam_1^k,\lam_2^k,\ldots,\lam_m^k)=(\lam_1^{k'},\lam_2^{k'},\ldots,\lam_m^{k'})$
as points in $\IC^m$ if and only if $k=k'$. But note
it can still happen that for some $j$, there may be repetitions among the quantities
$\lam_j^1,\lam_j^2,\ldots,\lam_j^p$, say $\lam_j^k=\lam_j^{k'}$ and thus $G_j^k=G_j^{k'}$,
even if $k\ne k'$.
\vv

\noindent {\bf Remark.}
A sufficient condition for $\lam_0$ to be an isolated point of $\spec(A_0)$,
as in the statement of Corollary~2.2, is easily given. Namely, assume that $\lam_0\in\spec(A_0)$ satisfies
$$
|\lam_0|>\max_{1\le j\le m}\{\rho_jr_j^{-1}\}r_1r_2\cdots r_m,
$$
where $r_j=r(A_j)$ and $\rho_j=\rho(A_j)$ are the spectral radii and essential spectral radii,
respectively, of these operators, and where we assume that $r_j>0$ for each $j$.
To prove that such $\lam_0$ is an isolated point of $\spec(A_0)$, it is enough to prove that
for every representation
$\lam_0=\lam_1\lam_2\cdots\lam_m$ where $\lam_j\in\spec(A_j)$, that each $\lam_j$
is an isolated point of $\spec(A_j)$. To this end, it is enough
to prove that $|\lam_j|>\rho_j$ for each $j$. Thus assume that $|\lam_{j_0}|\le\rho_{j_0}$
for some $j_0$. Then as $|\lam_j|\le r_j$ for each $j$, it follows that
$$
|\lam_0|=|\lam_1\lam_2\cdots\lam_m|\le(\rho_{j_0}r_{j_0}^{-1})r_1r_2\cdots r_m,
$$
which is a contradiction. Thus $\lam_0$ is isolated. In fact, one easily sees that $\lambda_0$ has finite
algebraic multiplicity, that is, $\nu_0=\dim G_0<\infty$.
\vv

\noindent {\bf Remark.}
If $\rho(A_j)=0$ for $1\le j\le m$, then
$\rho(A_0)=\rho(A_1\otimes A_2\otimes\cdots\otimes A_m)=0$, as one sees
by the above remark. This is indeed the case in
our analysis of delay-differential equations below.
\vv

\noindent {\bf Proof of Corollary~2.2.}
The fact that $\lam_0$ is a nonzero isolated point of $\spec(A_0)$, along with~\eqref{14}
from Theorem~2.1, implies that if $\lam_0=\lam_1\lam_2\cdots\lam_m$
with $\lam_j\in\spec(A_j)$, then each $\lam_j$ is an isolated point of $\spec(A_j)$.
This in turn implies that there is a finite number $p$ of such representations
of $\lam_0$ as a product. Let us enumerate all such representations, as in~\eqref{17a}
in the statement of the corollary.

For every $j$ satisfying $1\le j\le m$, let $q_j$ be the number of distinct quantities
$\lam_j^k$ for $1\le k\le p$. Here we mean numerically distinct quantities,
that is, without repetitions for multiplicity as an element of $\spec(A_j)$.
Let $\wt\lam_j^i$ for $1\le i\le q_j$ be a renumbering
of these quantities where each occurs only once, and so we have
$$
\{\lam_j^1,\lam_j^2,\ldots,\lam_j^p\}=\{\wt\lam_j^1,\wt\lam_j^2,\ldots,\wt\lam_j^{q_j}\},
$$
with equality as unordered sets.
Let $\wt G_j^i\subseteq X_j$ denote the spectral subspace
of $A_j$ corresponding to $\wt\lam_j^i$.
Then for each $j$ we have a direct sum decomposition
$$
X_j=\wt G_j^0\oplus(\wt G_j^1\oplus \wt G_j^2\oplus\cdots\oplus \wt G_j^{q_j}),
$$
where $\wt G_j^0$ is the spectral subspace of $A_j$ corresponding to
$\spec(A_j)\setminus\{\wt\lam_j^1,\wt\lam_j^2,\ldots,\wt\lam_j^{q_j}\}$.

Now consider all $m$-tuples $\iota=(i_1,i_2,\ldots,i_m)\in\cali$ where
$$
\cali=\{(i_1,i_2,\ldots,i_m)\in\ZZ^m\:|\: 0\le i_j\le q_j\hbox{ for every }j\hbox{ satisfying }1\le j\le m\}
$$
and for each such $\iota\in\cali$ let
\be{16}
\Gam^\iota=\wt G_1^{i_1}\otimes \wt G_2^{i_2}\otimes\cdots\otimes \wt G_m^{i_m}\subseteq
X_1\otimes X_2\otimes\cdots\otimes X_m.
\ee{16}
Then
$$
X_1\otimes X_2\otimes\cdots\otimes X_m=\bigoplus_{\iota\in\cali}\Gam^\iota.
$$
By construction, each subspace $\wt G_j^i\subseteq X_j$ is invariant for the operator $A_j$,
and thus each subspace $\Gam^\iota$ is invariant for $A_0$. Thus the multiplicity of $\lam_0$
as a point in the spectrum of $\spec(A_0)$, namely the dimension of the corresponding
spectral subspace, equals the sum of the multiplicities of $\lam_0$
as a point in the spectrum of $A_0|_{\Gam^\iota}$ for the various $\Gam^\iota$,
where $A_0|_{\Gam^\iota}$ is the restriction of $A_0$ to $\Gam^\iota$.

Note that not every $A_0|_{\Gam^\iota}$ need have $\lam_0$ in its spectrum. In fact,
there are precisely $p$ of the $m$-tuples $\iota\in\cali$ for which
\be{15}
\lam_0\in\spec(A_0|_{\Gam^\iota})
\ee{15}
holds, with these corresponding to the $p$ different $m$-tuples
in~\eqref{13}, \eqref{17a}. Moreover, it is the case that $i_j\ne 0$ for each $i_j$
occurring in such an $m$-tuple $\iota$, that is, the associated
subspace $\wt G_j^{i_j}$ is a spectral subspace of $A_j$ corresponding to $\wt\lam_j^{i_j}$, and
not the complementary space $\wt G_j^0$. The spectral subspace $G_0\subseteq X_0$
of $A_0$ corresponding to $\lam_0$ is thus the direct
sum of those $\Gam^\iota\subseteq X_0$ satisfying~\eqref{15}.
These facts are direct consequences of the
definition of the quantities $\lam_j^k$, along with the
re-labeling of the $\lam_j^k$ as $\wt\lam_j^i$ and the construction of the set $\cali$.
Let us denote by
$$
\iota^k=(i_1^k,i_2^k,\ldots,i_m^k),\qquad 1\le k\le p,
$$
those $\iota\in\cali$ for which~\eqref{15} holds. We may assume these $m$-tuples
are labeled to correspond with the $m$-tuples in~\eqref{13}, \eqref{17a}, namely that
\be{17}
(\wt\lam_1^{i_1^k},\wt\lam_2^{i_2^k},\ldots,\wt\lam_m^{i_m^k})
=(\lam_1^k,\lam_2^k,\ldots,\lam_m^k),\qquad 1\le k\le p.
\ee{17}
Thus the spectral subspace of $A_0$ for $\lam_0$ is the direct sum
\be{19}
G_0=\Gam^{\iota^1}\oplus\Gam^{\iota^2}\oplus\cdots\oplus\Gam^{\iota^p}
\ee{19}
in this notation. Then from~\eqref{16} and using~\eqref{17},
\be{20}
\Gam^{\iota^k}=\wt G_1^{i_1^k}\otimes \wt G_2^{i_2^k}\otimes\cdots\otimes \wt G_m^{i_m^k}
=G_1^k\otimes G_2^k\otimes\cdots\otimes G_m^k.
\ee{20}
Combining~\eqref{19} and~\eqref{20} gives the desired formula in~\eqref{18b} for $G_0$.
Furthermore, as we have defined $\nu_j^k=\dim G_j^k$, we obtain
the formula in~\eqref{18a} for $\nu_0=\dim G_0$.~\fp
\vv

Now take any Banach space $X$ and consider the $m$-fold tensor product, denoted
$$
X^{\otimes m}=X\otimes X\otimes\cdots\otimes X,
$$
with $m$ identical factors on the right-hand side.
Let $\cals_m$ denote the symmetric group on $m$ elements, namely
the set of all maps $\sig:\{1,2,\ldots,m\}\to\{1,2,\ldots,m\}$ which are
one-to-one and therefore onto. Taking any $\sig\in\cals_m$, we define
a linear operator $S_\sig\in\call(X^{\otimes m})$
as follows. Let
\be{cals}
S_\sig(x_1\otimes x_2\otimes\cdots\otimes x_m)
=x_{\sig(1)}\otimes x_{\sig(2)}\otimes\cdots\otimes x_{\sig(m)},
\ee{cals}
then extend $S_\sig$ to all of the algebraic tensor product
$X^{\odot m}=X\odot X\odot\cdots\odot X$ by linearity, and finally extend
$S_\sig$ to all of $X^{\otimes m}$ by continuity. One checks that $S_\sig$
is well-defined, and is an isometry, $\|S_\sig z\|=\|z\|$ for every $z\in X^{\otimes m}$.
Clearly, $S_{\sig_1}S_{\sig_2}=S_{\sig_1\sig_2}$ and $S_\sig^{-1}=S_{\sig^{-1}}$.
We now define the $m$-fold exterior product $X^{\wedge m}$ to be
$$
X^{\wedge m}=\{z\in X^{\otimes m}\:|\:S_\sig z=\sgn(\sig)z
\hbox{ for every }\sig\in\cals_m\},
$$
which is a closed subspace of $X^{\otimes m}$. Here $\sgn(\sig)=\pm 1$ is
the sign of the permutation $\sig$. Equivalently, we may define
$P\in\call(X^{\otimes m})$ by
\be{pdef}
P=\frac{1}{m!}\sum_{\sig\in\cals_m}\sgn(\sig)S_\sig,
\ee{pdef}
which is easily seen to be a projection, $P^2=P$. Then $X^{\wedge m}=PX^{\otimes m}$ is
the range of $P$, and we generally denote
$$
x_1\wedge x_2\wedge\cdots\wedge x_m=P(x_1\otimes x_2\otimes\cdots\otimes x_m).
$$
We note here, for future use, that
\be{psig}
PS_\sig=\sgn(\sig)P
\ee{psig}
for every $\sig\in\cals_m$.
Let us remark also that
\be{wdim}
\dim X^{\wedge m}=\binom{\dim X}{m},
\ee{wdim}
where $\binom{a}{b}$ denotes the binomial coefficient for $1\le a\le\infty$
and $1\le b<\infty$, with $\binom{a}{b}=0$ if $b>a$ and with $\binom{\infty}{b}=\infty$.
One easily checks~\eqref{wdim}, at least if $\dim X=n<\infty$, by noting that if
$e_1,e_2,\ldots,e_n\in X$ is a basis for $X$, then the set of elements
$e_{j_1}\wedge e_{j_2}\wedge\cdots\wedge e_{j_m}$ for $1\le j_1<j_2<\cdots<j_m\le n$
is a basis for $X^{\wedge m}$.

Now denoting
$$
A^{\otimes m}=A\otimes A\otimes\cdots\otimes A\in\call(X^{\otimes m})
$$
for the $m$-fold product of any operator $A\in\call(X)$ on $X$,
we observe that $S_\sig A^{\otimes m}=A^{\otimes m}S_\sig$ for every $\sig\in\cals_m$,
and thus $PA^{\otimes m}=A^{\otimes m}P$.
It follows that
$X^{\wedge m}$ is an invariant subspace of $X^{\otimes m}$ for $A^{\otimes m}$.
With this, it makes sense to study the spectrum of $A^{\otimes m}$ restricted to
$X^{\wedge m}$. Let us denote
$$
A^{\wedge m}=A^{\otimes m}|_{X^{\wedge m}}\in\call(X^{\wedge m})
$$
for this operator so restricted.
\vv

\noindent {\bf Proposition~2.3.} \bxx
Let $X$ be a Banach space and $A\in\call(X)$. Then for every $m\ge 1$
$$
\spec(A^{\wedge m})\subseteq\spec(A^{\otimes m})
=\{\lam_1\lam_2\cdots\lam_m\:|\:\lam_j\in\spec(A)\hbox{ for every }j\mbox{ satisfying }1\le j\le m\}
$$
for the operators $A^{\wedge m}\in\call(X^{\wedge m})$ and
$A^{\otimes m}\in\call(X^{\otimes m})$.
Suppose further that $\lam_0\in\spec(A^{\otimes m})$ is a nonzero isolated point of $\spec(A^{\otimes m})$
with spectral subspace
$G_0\subseteq X^{\otimes m}$. Then
\be{21}
PG_0=G_0\cap X^{\wedge m}
\ee{21}
for the image of this space under $P$.
Moreover, $PG_0\ne\{0\}$ if and only if $\lam_0\in\spec(A^{\wedge m})$,
in which case $\lam_0$ is an isolated point of $\spec(A^{\wedge m})$
with $PG_0$ as its spectral subspace.
\exx
\vv

\noindent {\bf Proof.}
The fact that $P$ commutes with $A^{\otimes m}$ implies that
$G_0$ is invariant under $P$, which in turn implies the equality in~\eqref{21}. The
remaining claims are elementary.~\fp
\vv

Corollary~2.2 may be used to obtain detailed information about the spectrum and spectral subspaces
of $A^{\wedge m}$. In this case each subspace $G_j^k\subseteq X$
in~\eqref{18b} is a spectral subspace of $A$, and it may happen
for a given $k$ that there are repetitions among these spaces, namely that
$G_{j}^k=G_{j'}^k$ and so $\lam_{j}^k=\lam_{j'}^k$, for some $j\ne j'$.
It is also the case that for every subspace $G_1^k\otimes G_2^k\otimes\cdots\otimes G_m^k$
occurring as a summand in~\eqref{18b}, and for every permutation $\sig\in\cals_m$, the space
obtained by permuting the factors $G_j^k$ using $\sig$ must also appear as a summand
in~\eqref{18b}. That is, there exists $k'$ such that
$$
G_1^{k'}\otimes G_2^{k'}\otimes\cdots\otimes G_m^{k'}=
S_\sig(G_1^k\otimes G_2^k\otimes\cdots\otimes G_m^k)=
G_{\sig(1)}^k\otimes G_{\sig(2)}^k\otimes\cdots\otimes G_{\sig(m)}^k.
$$
Of course, it may be the case that $k'=k$ even if $\sig$
is not the identity permutation, due to repetitions among the $G_j^k$.

The following result determines the multiplicity of a point $\lam_0$
in the spectrum of $A^{\wedge m}$, namely the quantity $\dim(PG_0)$
as in the statement of Proposition~2.3. Note that $\dim(PG_0)=0$
is possible, that is, it is possible that $\lam_0\in\spec(A^{\otimes m})$
but $\lam_0\not\in\spec(A^{\wedge m})$.
\vv

\noindent {\bf Proposition~2.4.} \bxx
Let $X$ be a Banach space and $A\in\call(X)$.
Fix $m\ge 1$ and let $\lam_0\in\spec(A^{\otimes m})$ be a nonzero
isolated point of $\spec(A^{\otimes m})$. For $1\le k\le p$ denote
$$
H^k=G_1^k\otimes G_2^k\otimes\cdots\otimes G_m^k,
$$
where we use the notation in the statement of Corollary~2.2.
Define an equivalence relation $\sim$ on the set $\{1,2,\ldots,p\}$
by letting $k\sim k'$ if and only if there exists $\sig\in\cals_m$ such that
$$
H^{k'}=S_\sig H^k,\qquad\hbox{that is,}\qquad G_j^{k'}=G_{\sig(j)}^k\hbox{ for }1\le j\le m.
$$
(Equivalently, $k\sim k'$ if and only if the two $m$-tuples in~\eqref{13} corresponding
to $k$ and $k'$ are obtained from one another by permuting the entries.)
Let $\cale^1,\cale^2,\ldots,\cale^r\subseteq\{1,2,\ldots,p\}$
denote the corresponding equivalence classes of~$\sim$ and let
\be{omq}
\Ome^q=\bigoplus_{k\in\cale^q}H^k
\ee{omq}
for each equivalence class, that is, for $1\le q\le r$. Then
\be{pg0}
PG_0=\bigoplus_{q=1}^rP\Ome^q,\qquad
\dim(PG_0)=\sum_{q=1}^r\dim(P\Ome^q),
\ee{pg0}
where $P$ is as in~\eqref{pdef} and $G_0\subseteq X^{\otimes m}$ is the spectral subspace of $\lam_0$
for $A^{\otimes m}$, as in the statement of Corollary~2.2.

Now fix any $q$ in the range $1\le q\le r$
and select an index $k_*\in\cale^q$ such that $H^{k_*}$ has the form
\be{hkap}
H^{k_*}=C_1^{\otimes\kappa_1}\otimes C_2^{\otimes\kappa_2}\otimes\cdots\otimes C_d^{\otimes\kappa_d},
\ee{hkap}
where for each $i$ we have that $C_i=G_j^{k_*}$ for some $j$, and where $C_{i}\ne C_{i'}$ and thus $C_{i}\cap C_{i'}=\{0\}$
if $i\ne i'$. The integers $\kappa_i\ge 1$ are thus precisely the number
times that $C_i$ occurs as a factor in this product.
(We remark that for any $q$ such $k_*$ exists, and that $k_*$ and $d$, and each $\kappa_i$ and $C_i$, of course depend on $q$.) Then
\be{dform2}
\dim(P\Ome^q)=\prod_{i=1}^d\binom{\dim C_i}{\kappa_i},
\ee{dform2}
with the convention in the above product that $0\times\infty=0$.
\exx
\vv

\noindent {\bf Remark.}
If, in the setting of Proposition~2.4, every nonzero point of $\spec(A)$
is an isolated point of $\spec(A)$, then the same is true for $\spec(A^{\otimes m})$.
In this case the nonzero points in the spectrum of $\spec(A^{\wedge m})$ are
precisely those points $\lam_0$ of the form
\be{lpr}
\lam_0=\lam_1\lam_2\cdots\lam_m,
\ee{lpr}
where each $\lam_j\in\spec(A)$ with possible
repetitions, but where the number of repetitions of each $\lam\in\spec(A)$ in the
product~\eqref{lpr} is less than or equal to the multiplicity of $\lam$ (the dimension of
the spectral subspace) as an element of $\spec(A)$. This means that
in the formula~\eqref{dform2}, one requires that $\kappa_i\le\dim C_i$ for each $i$.
\vv

\noindent {\bf Remark.}
Suppose, in the setting of Proposition~2.4, that every nonzero point of $\spec(A)$ is
an isolated point of simple multiplicity, that is,
an element of the point spectrum of algebraic multiplicity one. Let $\lam_j$ for $j\ge 1$
denote the distinct nonzero elements of $\spec(A)$. Then every nonzero $\lam_0\in\spec(A^{\wedge m})$
has the form
\be{lj}
\lam_0=\lam_{j_1}\lam_{j_2}\cdots\lam_{j_m}
\ee{lj}
for distinct integers $j_i$ satisfying $1\le j_1<j_2<\cdots<j_m$. Moreover, the multiplicity
of $\lam_0$ as an element of $\spec(A^{\wedge m})$ is precisely the number of possible
ways of expressing $\lam_0$ as such a product in this fashion.
One sees this easily from Proposition~2.4,
in particular, upon noting that in order for the quantity in~\eqref{dform2} to be positive
one must have each $\kappa_i=1$, as $\dim C_i=1$ for each $i$. Thus the space $H^{k_*}$
is a tensor product of $m$ spectral subspaces $C_i$ corresponding to distinct
points of $\spec(A)$ whose product is $\lam_0$.
\vv

\noindent {\bf Remark.}
Suppose, again in the setting of Proposition~2.4, that every nonzero point of $\spec(A)$ is
isolated. Suppose further there exists $r>0$ such that
there are exactly $m$ points $\lam\in\spec(A)$ satisfying $|\lam|>r$, and where here
we count multiplicity. That is, the spectral subspace corresponding to
all elements of $\spec(A)$ with $|\lam|>r$ has dimension exactly $m$. Denote these
elements of $\spec(A)$ by $\lam_j$, for $j=1,2,\ldots,m$, listed with repetition in
the case of multiplicity. Then $\lam_0=\lam_1\lam_2\cdots\lam_m$
is an isolated point of $\spec(A^{\wedge m})$ of simple multiplicity, namely its spectral
subspace has dimension $+1$. Further, there exists $\varepsilon>0$ such that
every other $\lam\in\spec(A^{\wedge m})$ satisfies
$|\lam|<|\lam_0|-\varepsilon$. Again, these facts follow easily from Proposition~2.4, where the
spaces $C_i$ are the spectral subspaces of the various $\lam_i$, with dimension equal
to the multiplicity of $\lam_i$, and where $\kappa_i=\dim C_i$.
\vv

\noindent {\bf Proof of Proposition~2.4.}
It is clear from ~\eqref{18b} and from~\eqref{omq}, and the fact that $\sim$ is an equivalence relation, that
\be{gom}
G_0=\bigoplus_{q=1}^r\Ome^q.
\ee{gom}
Further, it is clear using the definition of $\sim$
that $S_\sig\Ome^q=\Ome^q$ for every $\sig\in\cals_m$ and $1\le q\le r$, and so
\be{pomq}
P\Ome^q\subseteq\Ome^q
\ee{pomq}
holds. Thus~\eqref{pg0} follows from~\eqref{gom} and~\eqref{pomq}.

Now let $q$ be fixed, along with $k_*$, and $\kappa_i$ and $C_i$, as in the statement of the proposition.
For any $k\in\cale^q$ there exists $\pi\in\cals_m$ such that
\be{pi}
S_\pi H^k=H^{k_*},
\ee{pi}
and thus from~\eqref{psig}
we  see that $PH^k=PS_\pi H^k=PH^{k_*}$. It follows directly from this,
and from the definition~\eqref{omq} of $\Ome^q$, that
\be{phom}
P\Ome^q=PH^{k_*}.
\ee{phom}
Let us further denote $\Pi\in\call(\Ome^q,H^{k_*})$ to be the canonical projection
of $\Ome^q$ onto $H^{k_*}$ associated to the decomposition~\eqref{omq}. Also define the
isotropy group
$$
\Psi=\{\sig\in\cals_m\:|\:S_\sig H^{k_*}=H^{k_*}\}
$$
associated to the subspace $H^{k_*}$. We claim that
\be{pqp}
P\Pi P=\frac{|\Psi|}{m!}P\qquad\hbox{on }\Ome^q,
\ee{pqp}
where $|\Psi|$ denotes the cardinality of $\Psi$.
To prove~\eqref{pqp}, it is enough to verify that it holds on each subspace $H^k\subseteq\Ome^q$
for $k\in\cale^q$. Fixing such $k$, and with $\pi\in\cals_m$ satisfying~\eqref{pi}, take any $x\in H^k$
and denote $y=S_\pi x\in H^{k_*}$. Then
using~\eqref{psig} we have that
$$
Px=\sgn(\pi)PS_\pi x=\sgn(\pi)Py=\frac{\sgn(\pi)}{m!}\sum_{\sig\in\cals_m}\sgn(\sig)S_\sig y.
$$
Upon applying the operator $\Pi$, we retain only those terms in the above sum
which lie in $H^{k_*}$, namely, the terms for which $\sig\in\Psi$. Thus
$$
\Pi Px=\frac{\sgn(\pi)}{m!}\sum_{\sig\in\Psi}\sgn(\sig)S_\sig y.
$$
Applying $P$, where we again use~\eqref{psig}, now gives
$$
P\Pi Px=\frac{\sgn(\pi)}{m!}\sum_{\sig\in\Psi}Py=\frac{\sgn(\pi)|\Psi|}{m!}Py=\frac{|\Psi|}{m!}Px.
$$
From this we conclude~\eqref{pqp}, as desired.
It follows directly from~\eqref{pqp} that the map $\Pi$ is one-to-one on the space $P\Ome^q$.
Thus with~\eqref{phom} we conclude that
\be{dform}
\dim(P\Ome^q)=\dim(\Pi PH^{k_*}).
\ee{dform}

Let us now examine the isotropy group $\Psi$ more closely. 
As the spaces $C_i$ in the product~\eqref{hkap} are distinct, it follows
that $\sig\in\Psi$ if and only if $\sig$ permutes only those indices common
to each given term $C_i^{\otimes \kappa_i}$ among themselves without involving
other indices. More precisely,
define sets $\calk_i\subseteq\{1,2,\ldots,m\}$ for $1\le i\le d$ by
$$
\calk_i=\{n\in\ZZ\:|\:\wt\kappa_{i-1}<n\le\wt\kappa_i\},\qquad
\wt\kappa_i=\sum_{j=1}^i\kappa_j,\qquad
\wt\kappa_0=0,
$$
and so $\calk_i$ is the set of indices associated with the factor $C_i^{\otimes\kappa_i}$ in~\eqref{hkap},
and each $n$ in the range $1\le n\le m$ belongs to exactly one $\calk_i$.
Define subgroups $\Psi_i\subseteq\cals_m$, for $1\le i\le d$, by
$$
\Psi_i=\{\sig\in\cals_m\:|\:\sig(n)\in\calk_i\hbox{ if }n\in\calk_i,\hbox{ and }\sig(n)=n\hbox{ if }n\in\calk_j\hbox{ for some }j\ne i\},
$$
consisting of those $\sig$ which permute only the indices in $\calk_i$, leaving all other indices fixed.
Also define operators
$$
P_i=\frac{1}{\kappa_i!}\sum_{\sig\in\Psi_i}\sgn(\sig)S_\sig,\qquad
P_0=P_1P_2\cdots P_d,
$$
for $1\le i\le d$. The one easily sees that
$\Psi$ is precisely the set of elements
of the form
\be{sigi}
\sig=\sig_1\sig_2\cdots\sig_d
\ee{sigi}
with $\sig_i\in\Psi_i$ for $1\le i\le d$, and that the decomposition in~\eqref{sigi} is unique
for each $\sig\in\Psi$. Note the commutativity property, that
$\sig_i\sig_{i'}=\sig_{i'}\sig_i$ if $\sig_i\in\Psi_{i}$ and $\sig_{i'}\in\Psi_{i'}$ with $i\ne i'$.
One sees that operator $P_i^2=P_i$ is a projection on $H^{k_*}$ whose range is the space
$$
C_1^{\otimes \kappa_1}\otimes\cdots\otimes
C_{i-1}^{\otimes \kappa_{i-1}}\otimes C_i^{\wedge \kappa_i}\otimes C_{i+1}^{\otimes \kappa_{i+1}}
\otimes\cdots\otimes C_d^{\otimes \kappa_d},
$$
and using the above-mentioned commutativity, one sees that $P_0^2=P_0$ is also a projection on $H^{k_*}$
whose range is the subspace
$$
C_1^{\wedge\kappa_1}\otimes C_2^{\wedge\kappa_2}\otimes\cdots\otimes C_d^{\wedge\kappa_d}.
$$
We claim that
\be{clm}
\Pi P=\frac{\kappa_1!\kappa_2!\cdots\kappa_d!}{m!}P_0\qquad\hbox{on }H^{k_*},
\ee{clm}
from which it follows directly, with the above remarks, that
\be{qph}
\Pi PH^{k_*}=C_1^{\wedge\kappa_1}\otimes C_2^{\wedge\kappa_2}\otimes\cdots\otimes C_d^{\wedge\kappa_d}.
\ee{qph}
Note that~\eqref{qph}, along with~\eqref{wdim} and~\eqref{dform}, implies our desired result~\eqref{dform2}.
To prove~\eqref{clm}, first observe that for every $x\in H^{k_*}$ we have that
\be{pip}
\Pi Px=\frac{1}{m!}\sum_{\sig\in\cals_m}\sgn(\sig)\Pi S_\sig x =\frac{1}{m!}\sum_{\sig\in\Psi}\sgn(\sig)S_\sig x.
\ee{pip}
Now decomposing $\sig\in\Psi$ as in~\eqref{sigi}, we have that
$$
\begin{array}{lcl}
\ds{\sum_{\sig\in\Psi}\sgn(\sig)S_\sig x}
&\bb = &\bb
\ds{\sum_{\sig_1\in\Psi_1}\sum_{\sig_2\in\Psi_2}\cdots\sum_{\sig_d\in\Psi_d}
\sgn(\sig_1)\sgn(\sig_2)\cdots\sgn(\sig_d)S_{\sig_1} S_{\sig_2}\cdots S_{\sig_d}x}\\
\\
&\bb = &\bb
\ds{(\kappa_1!\kappa_2!\cdots\kappa_d!)P_1P_2\cdots P_dx=(\kappa_1!\kappa_2!\cdots\kappa_d!)P_0x},
\end{array}
$$
which with~\eqref{pip}, proves the claim~\eqref{clm}. With this, the proposition is proved.~\fp
\vv

We now consider the specific case of Banach spaces
$$
X_j=C(\Theta_j)=\{\vph:\Theta_j\to\IR\:|\:\vph\hbox{ is continuous}\}
$$
for $1\le j\le m$, where each $\Theta_j$ is a compact Hausdorff space and where the supremum norm
is taken for $C(\Theta_j)$. As described in~\cite[Chapter~I, Section~4]{df}, one may regard
\be{05}
C(\Theta_0)=C(\Theta_1)\otimes C(\Theta_2)\otimes\cdots\otimes C(\Theta_m),\qquad
\Theta_0=\Theta_1\times\Theta_2\times\cdots\times\Theta_m,
\ee{05}
as follows. First, taking any $\vph_j\in C(\Theta_j)$ for $1\le j\le m$,
define $\vph\in C(\Theta_0)$ by
\be{phi}
\vph(\theta_1,\theta_2,\ldots,\theta_m)=\vph_1(\theta_1)\vph_2(\theta_2)\cdots\vph_m(\theta_m),
\ee{phi}
and identify $\vph_1\otimes\vph_2\otimes\cdots\otimes\vph_m$
with $\vph$. More generally, identify any finite sum
$$
\sum_{i=1}^n\vph_{1,i}\otimes\vph_{2,i}\otimes\cdots\otimes\vph_{m,i}\in
C(\Theta_1)\otimes C(\Theta_2)\otimes\cdots\otimes C(\Theta_m)
$$
where $\vph_{j,i}\in C(\Theta_j)$ for $1\le i\le n$ and $1\le j\le m$, with
$\vph\in C(\Theta_0)$ given by
\be{06}
\vph(\theta_1,\theta_2,\ldots,\theta_m)=
\sum_{i=1}^n\vph_{1,i}(\theta_1)\vph_{2,i}(\theta_2)\cdots\vph_{m,i}(\theta_m).
\ee{06}
One sees that this identification is an isometry, that is,
\be{04}
\|\vph\|=\bigg\|\sum_{i=1}^n\vph_{1,i}\otimes\vph_{2,i}\otimes\cdots\otimes\vph_{m,i}\bigg\|,
\ee{04}
where the norms in~\eqref{04} are those in $C(\Theta_0)$ and in $C(\Theta_1)\otimes C(\Theta_2)\otimes\cdots\otimes C(\Theta_m)$,
respectively.
To prove~\eqref{04}, first take elements $\xi_j\in C(\Theta_j)^*$ of the dual spaces, with $\|\xi_j\|=1$,
for $1\le j\le m$. Each $\xi_j$ is given by integration with respect to a Borel measure
$d\mu_j(\theta_j)$ on $\Theta_j$ with total variation $|\mu_j|(\Theta_j)=1$. Then with~\eqref{06}
we have, following~\eqref{02}, that
\be{meas}
\begin{array}{l}
\ds{\bigg|\sum_{i=1}^n\xi_1(\vph_{1,i})\xi_2(\vph_{2,i})\cdots\xi_m(\vph_{m,i})\bigg|}\\
\\
\qquad\qquad\ds{=\bigg|\int_{\Theta_1}\int_{\Theta_2}\!\!\cdots\int_{\Theta_m}\vph(\theta_1,\theta_2,\ldots,\theta_m)
\:d\mu_m(\theta_m)\cdots d\mu_2(\theta_2)\:d\mu_1(\theta_1)\bigg|
\le\|\vph\|.}
\end{array}
\ee{meas}
Upon taking the supremum over all such $\xi_j$, we have that
\be{measinq}
\bigg\|\sum_{i=1}^n\vph_{1,i}\otimes\vph_{2,i}\otimes\cdots\otimes\vph_{m,i}\bigg\|\le\|\vph\|.
\ee{measinq}
To obtain equality in~\eqref{measinq}, take any point
$(\theta_1^*,\theta_2^*,\ldots,\theta_m^*)\in\Theta_1\times\Theta_2\times\cdots\times\Theta_m$
at which the maximum of $|\vph(\theta_1,\theta_2,\ldots,\theta_m)|$ is achieved, where without loss,
by multiplying $\vph$ by a scalar of norm $+1$, we may assume that
$\vph(\theta_1^*,\theta_2^*,\ldots,\theta_m^*)=\|\vph\|\ge 0$.
Then letting $d\mu_j(\theta_j)$ be the unit point mass at $\theta_j^*$, we see that the integral expression~\eqref{meas}
equals $\|\vph\|$, and thus equality holds in~\eqref{measinq}. This establishes~\eqref{04}.
With~\eqref{06} and~\eqref{04}, it follows that the space $C(\Theta_1)\otimes C(\Theta_2)\otimes\cdots\otimes C(\Theta_m)$ is
isometrically embedded as a subspace of
$C(\Theta_0)$. In fact this subspace is all of $C(\Theta_0)$, that is,
the first equality in~\eqref{05} holds.
This follows directly from the fact that the set of functions $\vph$ of the
form~\eqref{06} is dense in $C(\Theta_0)$, by the Stone-Weierstrass Theorem.

Suppose further that $A_j\in\call(C(\Theta_j))$ for $1\le j\le m$. For $1\le k\le m$ define
an operator $\wt A_k\in\call(C(\Theta_0))$ by
$$
\wt A_k=I\otimes \cdots I\otimes A_k\otimes I\cdots \otimes I,
$$
where the factor $A_k$ occurs in the $k^{\rm th}$ position. Then if $\vph_j\in C(\Theta_j)$ for $1\le j\le m$
and with $\vph$ given by~\eqref{phi},
we have that
$$
(\wt A_k\vph)(\theta_1,\ldots,\theta_m)
=\vph_1(\theta_1)\ldots\vph_{k-1}(\theta_{k-1})[(A_k\vph_k)(\theta_k)]\vph_{k+1}(\theta_{k+1})\cdots\vph_m(\theta_m)
$$
for every $(\theta_1,\ldots,\theta_m)\in\Theta_0$, that is, $A_k$ acts upon the function $\vph_k$ with the
other functions $\vph_j$ for $j\ne k$ untouched. More generally, for any
$\vph \in C(\Theta_0)$ not necessarily of the product form~\eqref{phi}, one has that
\be{cd}
(\wt A_k\vph)(\theta_1,\ldots,\theta_m)=[A_k\vph(\theta_1,\ldots,\theta_{k-1},\,\cdot\,,\theta_{k+1},\ldots,\theta_m)](\theta_k)
\ee{cd}
which is interpreted as follows. Let the points $\theta_j\in \Theta_j$ for $j\ne k$ be held fixed
and regard
$$
\vph(\theta_1,\ldots,\theta_{k-1},\,\cdot\,,\theta_{k+1},\ldots,\theta_m)\in C(\Theta_k)
$$
as a function
of one variable represented by the centered dot ``$\,\cdot\,$''. Apply the operator $A_k$
to this function, and then evaluate the resulting function at the point $\theta_k\in \Theta_k$
to get the right-hand side of~\eqref{cd}. It follows that
to calculate $(A_0\vph)(\theta_1,\theta_2,\ldots,\theta_m)$ where $A_0\in\call(C(\Theta_0))$ is the operator
$$
A_0=A_1\otimes A_2\otimes\cdots\otimes A_m=\wt A_1\wt A_2\cdots\wt A_m,
$$
one successively applies the operators $A_k$ for $1\le k\le m$ with the variable in the $k^{\rm th}$ position free,
while holding the remaining $m-1$ variables fixed. Note that one may apply these operators in any order,
as the operators $\wt A_k$ commute with one another.

In the special case that all the spaces $X_j=X=C(\Theta)$ are the same, then
$\Theta_0=\Theta^m$, the $m$-fold cartesian product, and so we have the identification
$C(\Theta)^{\otimes m}=C(\Theta^m)$. Further, it is clear that $C(\Theta)^{\wedge m}$ is identified
with the subspace of $C(\Theta^m)$ consisting of all anti-symmetric functions, that is,
\be{as}
\begin{array}{lcl}
C(\Theta)^{\wedge m} &\bb = &\bb
\{\vph\in C(\Theta^m)\:|\:
\vph(\theta_{\sig(1)},\theta_{\sig(2)},\ldots,\theta_{\sig(m)})
=\sgn(\sig)\vph(\theta_1,\theta_2,\ldots,\theta_m)\\
\\
&\bb &\bb
\:\:\hbox{for every }(\theta_1,\theta_2,\ldots,\theta_m)\in\Theta^m,\hbox{ and every }\sig\in\cals_m\}.
\end{array}
\ee{as}

As a practical matter, the above observations will be useful in evaluating
tensor products of solution operators of linear delay-differential equations.
In such applications we shall typically work
with the exterior product space $C([-1,0])^{\wedge m}$.

The following basic result will be needed later. Although it is proved in~\cite{ich1},
we provide a proof for completeness.
\vv

\noindent {\bf Proposition~2.5 (Ichinose~\cite[Lemma~3{.}6]{ich1}).} \bxx
Let $X_j$ and $Y_j$ be Banach spaces and $A_j\in\call(X_j,Y_j)$ for $1\le j\le m$.
Assume that each operator $A_j$ is one-to-one. Then the operator
$A_0=A_1\otimes A_2\otimes\cdots\otimes A_m$ is one-to-one from
$X_1\otimes X_2\otimes\cdots\otimes X_m$ to $Y_1\otimes Y_2\otimes\cdots\otimes Y_m$.
\exx
\vv

\noindent {\bf Proof.}
Without loss it is enough to consider the case $m=2$, as the
case of general $m$ can be proved inductively by writing
$A_0=A_{00}\otimes A_m$ where $A_{00}=A_1\otimes A_2\otimes\cdots\otimes A_{m-1}$.
Further, if $m=2$, then by writing
$A_1\otimes A_2=(A_1\otimes I_{Y_2})(I_{X_1}\otimes A_2)$ where $I_{X_1}$ and $I_{Y_2}$
denote the identity operators on $X_1$ and $Y_2$ respectively, we see that it is enough
to prove that both $A_1\otimes I_{Y_2}$ and $I_{X_1}\otimes A_2$ are one-to-one.
In fact, it is enough to prove that the operator $A_1\otimes I_{Y_2}$ is one-to-one.

Therefore, denoting $A=A_1$, $X=X_1$, $Y=Y_1$, and $Z=Y_2$, let us consider an operator $A\in\call(X,Y)$ which is
one-to-one. We must prove that $A\otimes I\in\call(X\otimes Z,Y\otimes Z)$ is also one-to-one,
where $I$ denotes the identity operator on $Z$. Letting $Z^*$ denote the dual space of $Z$,
for any $\zeta\in Z^*$ define an operator $L_X(\zeta)\in\call(X\otimes Z,X)$ by setting
$$
L_X(\zeta)(x\otimes z)=\zeta(z)x
$$
for any $x\in X$ and $z\in Z$, and then extending $L_X(\zeta)$ to all of $X\otimes Z$ first
by linearity and then by continuity. One easily sees that $L_X(\zeta)$ is well-defined, with operator norm
$$
\|L_X(\zeta)\|=\|\zeta\|.
$$
One also easily checks that
\be{yy}
\|u\|=\sup_{{\scriptstyle \zeta\in Z^*}\atop {\scriptstyle \|\zeta\|=1}}\|L_X(\zeta)u\|
\ee{yy}
for every $u\in X\otimes Z$, which in fact follows directly from
the definition~\eqref{02} of the injective norm. We also define the operator
$L_Y(\zeta)\in\call(Y\otimes Z,Y)$ in an analogous fashion. Finally, let us note that
\be{x}
L_Y(\zeta)(A\otimes I)=AL_X(\zeta)
\ee{x}
for every $\zeta\in Z$, which one easily sees by showing that the operators in~\eqref{x}
agree on all elements $x\otimes z\in X\otimes Z$.

Now assume that $(A\otimes I)u=0$ for some $u\in X\otimes Z$. Then from~\eqref{x} we have
for every $\zeta\in Z^*$ that $AL_X(\zeta)u=0$, and hence that $L_X(\zeta)u=0$
as $A$ is one-to-one. But then~\eqref{yy} implies that $\|u\|=0$, thus $u=0$.
We conclude that $A\otimes I$ is one-to-one, as desired.~\fp

\section{Tensor Products of Linear Processes}

Before specializing to the delay-differential equation~\eqref{int1},
we make some general observations about abstract linear processes. These
observations not only apply to~\eqref{int1}, but also to a large class of linear nonautonomous
delay-differential equations as well as to many other systems.

By a {\bf linear process} (sometimes called a {\bf linear evolutionary process}) $U(t,\tau)$ on a Banach space $X$, we mean a family of
bounded linear operators $U(t,\tau)\in\call(X)$, for every $t,\tau\in\IR$ with $t\ge\tau$, for which
\begin{itemize}
\item[{(1)}] $U(\tau,\tau)=I$ for every $\tau\in\IR$;
\item[{(2)}] $U(t,\sig)U(\sig,\tau)=U(t,\tau)$ for every $t,\sig,\tau\in\IR$ with $t\ge\sig\ge\tau$; and
\item[{(3)}] $U(t,\tau)$ is strongly continuous in $t$ and $\tau$, that is, for every $x\in X$ it is the
case that $U(t,\tau)x$ varies continuously in $X$ as a function of $t$ and $\tau$, for $t\ge\tau$.
\end{itemize}
It is easy to check, using the uniform boundedness principle, that there is a bound $\|U(t,\tau)\|\le C$
in the neighborhood of any point $(t_0,\tau_0)$ in the domain of $U(\cdot,\cdot)$, where such $C$ depends on $(t_0,\tau_0)$.

Now fix an integer $m\ge 1$ and consider the $m$-fold tensor product
$X^{\otimes m}$.
For every $k$ satisfying $1\le k\le m$, and with $t$ and $\tau$ as before,
we may define an operator $U_k(t,\tau)\in\call(X^{\otimes m})$ by
\be{uk}
U_k(t,\tau)=I\otimes\cdots\otimes I\otimes U(t,\tau)\otimes I\otimes\cdots\otimes I,
\ee{uk}
where the factor $U(t,\tau)$ occurs in the $k^{\rm th}$ place.
It is easily checked that $U_k(t,\tau)$ is a linear process on $X^{\otimes n}$. Also,
one has from~\eqref{07} that
\be{08}
U_k(t,\tau)U_j(t',\tau')=U_j(t',\tau')U_k(t,\tau)
\ee{08}
for any real numbers $t\ge\tau$ and $t'\ge\tau'$,
with $j\ne k$ in the range $1\le j,k\le m$.
Next define the operator $\UU(t,\tau)\in\call(X^{\otimes m})$ for $t\ge\tau$ by
\be{bigu}
\UU(t,\tau)=U(t,\tau)^{\otimes m}=U_1(t,\tau)U_2(t,\tau)\cdots U_m(t,\tau),
\ee{bigu}
where it does not matter in what order the above product is taken
due to the commutativity~\eqref{08}.
Again, $\UU(t,\tau)$ is a linear process on $X^{\otimes m}$.
It is also clear that the subspace $X^{\wedge m}\subseteq X^{\otimes m}$ is invariant under
$\UU(t,\tau)$, and we shall denote by
\be{bigw}
\WW(t,\tau)=U(t,\tau)^{\wedge m}=\UU(t,\tau)|_{X^{\wedge m}}
\ee{bigw}
the restriction of this linear process to $X^{\wedge m}$.
Certainly, $\WW(t,\tau)\in\call(X^{\wedge m})$ is itself a linear process on the space $X^{\wedge m}$.

It often happens that a linear process $U(t,\tau)$ is periodic, meaning that
there exists some $\gam>0$ such that
$$
U(t+\gam,\tau+\gam)=U(t,\tau)
$$
for every $t$ and $\tau$ with $t\ge\tau$. In this case, for each $\tau\in\IR$ we define
$$
M(\tau)=U(\tau+\gam,\tau),
$$
the so-called {\bf monodromy operator} with initial time $\tau$,
and we note that $M(\tau+\gam)=M(\tau)$.
We refer to the nonzero spectrum of $M(\tau)$ as the {\bf Floquet spectrum} of the linear process,
and we call the set of nonzero elements in the point spectrum of $M(\tau)$ the
set of {\bf Floquet multipliers}.
Let us observe that the Floquet spectrum does not depend on the
initial time $\tau$, that is, $\spec(M(\tau))\setminus\{0\}=\spec(M(\tau'))\setminus\{0\}$
for every $\tau,\tau'\in\IR$. To prove this, without
loss we may take $\tau'=0$ and let $\tau$
lie in the range $0<\tau<\gam$. With $\tau$ so fixed, we have that
\be{decomp}
M(0)=AB,\qquad M(\tau)=BA,\qquad\hbox{where}\qquad A=U(\gam,\tau),\qquad B=U(\tau,0),
\ee{decomp}
and we must show that $\spec(AB)\setminus\{0\}=\spec(BA)\setminus\{0\}$.
In fact, this is a well-known result for any pair of operators. We sketch the proof.
Taking any $\lam\in\IC\setminus\{0\}$ satisfying $\lam\not\in\spec(AB)$,
one easily checks by multiplication that the operator $\lam^{-1}I+\lam^{-1}B(\lam I-AB)^{-1}A$
is the inverse of $\lam I-BA$, and so $\lam\not\in\spec(BA)$. Thus
$\spec(AB)\setminus\{0\}\supseteq\spec(BA)\setminus\{0\}$, with the opposite inclusion
proved similarly.

It is also the case that the nonzero point spectra of $AB$ and $BA$ are the same,
and so the set of Floquet multipliers is independent of the initial time. Indeed, if
$\lam\ne 0$ is in the point spectrum of $M(0)$ then $ABx=\lam x\ne 0$ for
some $x\in X$. Letting $y=Bx\ne 0$, we thus have that $M(\tau)y=BAy=BABx=\lam Bx=\lam y\ne 0$,
and so $\lam$ is in the point spectrum of $M(\tau)$.

It is easily seen from~\eqref{decomp} that if for some $\tau'\in\IR$ and some $n\ge 1$ the operator
$M(\tau')^n$ is compact, then for every $\tau\in\IR$ the operator
$M(\tau)^{n+1}$ is compact. In this case the remarks above imply that
the Floquet spectrum consists entirely of Floquet multipliers.
This will indeed be the case in our studies of delay-differential equations
below.
\vv

\noindent {\bf Remark.}
As noted by Muldowney~\cite{muld},
the above observations in connection with compound (exterior product) systems
have great relevance for the stability of periodic solutions of nonlinear systems.
For example, in the case of an autonomous ordinary differential equation, say
\be{nonl}
\dot x=f(x),\qquad x\in\IR^n,
\ee{nonl}
suppose that $x=\xi(t)$ is a nonconstant periodic solution of minimal period $\gam>0$.
Consider the associated linearized system $\dot y=A(t)y$ where $A(t)=f'(\xi(t))$,
with $U(t,\tau)$ the associated fundamental solution with $U(\tau,\tau)=I$. Let the
Floquet multipliers (the characteristic multipliers) be ordered so that
$|\lam_1|\ge|\lam_2|\ge\cdots\ge|\lam_n|$, with repetitions according to algebraic
multiplicities, and recall that $\lam_k=1$ for
some $k$, the so-called trivial multiplier. Then the periodic solution $\xi(t)$
is exponentially asymptotically stable for the nonlinear system~\eqref{nonl} if and only if
\be{exst}
\lam_1=1>|\lam_2|.
\ee{exst}
Further, consider the second compound linear process $\WW(t,\tau)=U(t,\tau)^{\wedge 2}=U(t,\tau)\wedge U(t,\tau)$,
which acts on the space $\IR^n\wedge\IR^n$ of dimension $\binom{n}{2}=\frac{1}{2}n(n-1)$
and whose Floquet multipliers are precisely the quantities $\mu=\lam_i\lam_j$ for $1\le i<j\le n$.
Then $|\mu|<1$ for every such $\mu$ if and only if~\eqref{exst} holds, that is,
if and only if $\xi(t)$ is exponentially asymptotically stable.

The appropriate generalizations of this conclusion hold for a wide variety of infinite dimensional systems,
including reaction-diffusion equations (see Wang~\cite{wang}), and
a large class of retarded functional differential equations $\dot x(t)=f(x_t)$, not limited
to a single delay. (Here we follow the notation of~\cite{hale}.)

\section{Tensor Products of Delay-Differential Equations}

Consider the linear scalar delay-differential equation
\be{0}
\dot x(t)=-\alp(t)x(t)-\bet(t)x(t-1)
\ee{0}
which is equation~\eqref{int1} from the Introduction. The change of variables
\be{cov}
y=\mu(t)x,\qquad\mu(t)=\exp\bigg(\int_0^t\alp(s)\:ds\bigg),
\ee{cov}
transforms~\eqref{0} into the equivalent equation
\be{1}
\dot y(t)=-b(t)y(t-1),
\ee{1}
where
\be{xyz}
b(t)=\frac{\mu(t)\bet(t)}{\mu(t-1)}=\bet(t)\exp\bigg(\int_{t-1}^t\alp(s)\:ds\bigg).
\ee{xyz}
Note that $\bet(t)$ and $b(t)$ have the same sign for every $t$.

As a standing hypothesis, for the remainder of the paper we shall assume that $\alp,\bet:\IR\to\IR$
and therefore $b:\IR\to\IR$ are locally integrable functions. (Actually, it will sometimes be enough only to assume
these properties on the interval $[\tau,\tau+\eta]$ considered in our results.)
Often we shall work with the simpler equation~\eqref{1} and interpret our results
back for equation~\eqref{0}.

Both equations~\eqref{0} and~\eqref{1} generate linear
processes, which we denote by $U(t,\tau)$ and $\wt U(t,\tau)$, respectively, on the Banach space
$$
X=C([-1,0]),
$$
and where we keep
this notation for the remainder of the paper.
In particular, $U(t,\tau)\in\call(X)$ for any $t,\tau\in\IR$ with $t\ge\tau$
denotes the associated solution operator
on $X$ to equation~\eqref{0} defined as
$$
U(t,\tau)\vph=x_{t}.
$$
Here $x(t)$ satisfies~\eqref{0} for $t\ge\tau$, with $x_t\in X$
defined in the usual fashion~\cite{hale} by $x_t(\theta)=x(t+\theta)$ for $\theta\in[-1,0]$,
and where we take the initial condition
$x_\tau=\vph\in X$, that is, $x(\tau+\theta)=\vph(\theta)$ for $\theta\in[-1,0]$.
Similarly, $\wt U(t,\tau)$ denotes the analogous solution operator for equation~\eqref{1},
and one sees the relation
\be{conj}
\wt U(t,\tau)=\Sig(t)U(t,\tau)\Sig(\tau)^{-1}
\ee{conj}
between these two processes, where $\Sig(t)\in\call(X)$ is defined to be the multiplication
operator
$$
[\Sig(t)\vph](\theta)=\mu(t+\theta)\vph(\theta),\qquad\theta\in[-1,0],
$$
for any $t\in\IR$.

Now fix an integer $m\ge 1$ and recall the identification
$$
X^{\otimes m}=C([-1,0]^m)
$$
of the $m$-fold tensor product
as described in Section~2.
We shall denote the argument of a function
in $X^{\otimes m}$ by $\theta=(\theta_1,\theta_2,\ldots,\theta_m)\in[-1,0]^m$, that is,
we write $\vph(\theta)=\vph(\theta_1,\theta_2,\ldots,\theta_m)$.
Also recall that $X^{\wedge m}$ is identified with the
antisymmetric elements of $C([-1,0]^m)$, namely
$$
\begin{array}{lcl}
X^{\wedge m} &\bb = &\bb
\{\vph\in C([-1,0]^m)\:|\:
\vph(\theta_{\sig(1)},\theta_{\sig(2)},\ldots,\theta_{\sig(m)})
=\sgn(\sig)\vph(\theta_1,\theta_2,\ldots,\theta_m)\\
\\
&\bb &\bb
\:\:\hbox{for every }(\theta_1,\theta_2,\ldots,\theta_m)\in[-1,0]^m,\hbox{ and every }\sig\in\cals_m\},
\end{array}
$$
as in~\eqref{as}.
Note that if
$\vph\in X^{\wedge m}$ then
the values of $\vph(\theta)$ for $\theta\in[-1,0]^m$ are completely
determined by the values for which $\theta\in T_m$, where
\be{tri}
T_m=\{\theta=(\theta_1,\theta_2,\ldots,\theta_m)\in[-1,0]^m\:|\:\theta_1\le\theta_2\le\cdots\le\theta_m\}.
\ee{tri}
Recall further the operators $U_k(t,\tau),\UU(t,\tau)\in\call(X^{\otimes m})$ as in~\eqref{uk} and~\eqref{bigu},
and the restriction $\WW(t,\tau)\in\call(X^{\wedge m})$ of $\UU(t,\tau)$
to the invariant subspace $X^{\wedge m}\subseteq X^{\otimes m}$
as in~\eqref{bigw}, which give linear processes in their respective spaces.
Let $\wt U_k(t,\tau)$, $\wt\UU(t,\tau)$, and $\wwtilde(t,\tau)$ denote the corresponding operators for equation~\eqref{1}.

Our main interest will be positivity properties of the operator $\WW(t,\tau)$ on $X^{\wedge m}$,
where positivity is described in terms of the cone
\be{km}
K_m=\{\vph\in X^{\wedge m}\:|\:\vph(\theta)\ge 0\hbox{ for every }\theta\in T_m\}.
\ee{km}
We shall prove the following theorem, which is one of our main results. To our knowledge, this
result is new even in the case of the constant coefficient version
\be{cc2}
\dot x(t)=-\alp_0x(t)-\bet_0x(t-1)
\ee{cc2}
of equation~\eqref{0}.
\vv

\noindent {\bf Theorem~4.1 (Positivity Theorem).} \bxx
Fix $m\ge 1$, and let $\tau\in\IR$ and $\eta\ge 0$.
Assume that $(-1)^m\beta(t)\ge 0$ for almost every $t\in[\tau,\tau+\eta]$ for the delay coefficient
in equation~\eqref{0}. Then the
operator $\WW(\tau+\eta,\tau)\in\call(X^{\wedge m})$ defined in~\eqref{bigw}
is a positive operator with respect to the cone $K_m$ in~\eqref{km},
that is, $\WW(\tau+\eta,\tau)$ maps $K_m$ into itself.
\exx
\vv

Recall that if $Y$ is any Banach space, then a closed, convex subset $K\subseteq Y$ is called
a {\bf cone} if whenever $u\in K$ then $\sigma u\in K$ for every $\sigma\ge 0$, and if both
$u,-u\in K$ only if $u=0$. If $u,v\in Y$ and $K\subseteq Y$ is a cone, we write $u\le v$
to mean $v-u\in K$. Also, if $A\in\call(Y)$
then we say the operator $A$ is {\bf positive} if it maps $K$ into $K$, that is, $Au\in K$
whenever $u\in K$, and we write $A\ge 0$. (And we write $A\le B$ to mean $B-A\ge 0$.)
We say that a cone $K$ is {\bf reproducing} if $Y=\{u-v\:|\:u,v\in K\}$, and
we say that $K$ is {\bf normal} if there exists a constant $C>0$
such that $\|u\|\le C\|v\|$ whenever $0\le u\le v$.
One sees directly that the cone $K_m$ is both reproducing and normal in $X^{\wedge m}$.

The above Positivity Theorem is a straightforward consequence
of Proposition~4.2 below, which provides more detailed information for the transformed equation~\eqref{1}.
The Positivity Theorem follows from this using the conjugacy~\eqref{conj}
and the fact that $\Sig(t)^{\wedge m}$ and its inverse $[\Sig(t)^{\wedge m}]^{-1}$ are
positive operators on $X^{\wedge m}$ with respect to the cone $K_m$.

From the remarks preceding Proposition~2.5, we have that
$U_k(t,\tau)$ and $\wt U_k(t,\tau)$ are simply
the solution operators to equations~\eqref{0} and~\eqref{1} taken
along the $k^{\rm th}$ coordinate in $[-1,0]^m$, with the remaining
$m-1$ coordinates staying fixed. Therefore, suppose that $x^k(t)$, for $1\le k\le m$, are solutions
of~\eqref{0} for $t\in[\tau,\tau+\eta]$ with initial conditions $x^k_\tau\in C[-1,0]$. For such $t$ define
\be{xi}
\begin{array}{lcl}
\XX_t(\theta) &\bb = &\bb
x^1(t+\theta_1)x^2(t+\theta_2)\cdots x^m(t+\theta_m),\\
\\
\XXI_t(\theta) &\bb = &\bb
\det\left(\begin{array}{ccc}
x^1(t+\theta_1) &\bb \cdots &\bb x^1(t+\theta_m)\\
\vdots &\bb &\bb \vdots\\
x^m(t+\theta_1) &\bb \cdots &\bb x^m(t+\theta_m)
\end{array}\right),
\end{array}
\ee{xi}
for $\theta\in[-1,0]$ and observe that
$\XX_t\in X^{\otimes m}$ and $\XXI_t\in X^{\wedge m}$.
Then it follows immediately that
\be{xi2}
\XX_t=\UU(t,\tau)\XX_\tau,\qquad
\XXI_t=\WW(t,\tau)\XXI_\tau,
\ee{xi2}
for $t\in[\tau,\tau+\eta]$.

If $0\le\eta\le 1$ then explicit formulas can be given. We consider the transformed
equation~\eqref{1}. First, let $\psi\in C[-1,0]$. Then
\be{m1}
[\wt U(\tau+\eta,\tau)\psi](\theta)=
\left\{\begin{array}{lcl}
\psi(\eta+\theta), & \hbox{ for} &\bb -1\le\theta\le -\eta,\\
\\
\ds{\psi(0)-\int^{\eta+\theta}_0b(\tau+s)\psi(s-1)\:ds,} & \hbox{ for} &\bb -\eta\le\theta\le 0,
\end{array}\right.
\ee{m1}
if $0\le\eta\le 1$.
Now fix $k$ in the range $1\le k\le m$ and for any
$\theta=(\theta_1,\theta_2,\ldots,\theta_m)\in[-1,0]^m$ let
$$
\wh\theta=(\theta_1,\ldots,\theta_{k-1},\theta_{k+1},\ldots,\theta_m)\in[-1,0]^{m-1},
$$
which is $\theta$ with the $k^{\rm th}$ coordinate removed.
Then regarding $\wh\theta$ as a fixed parameter and taking $\vph\in X^{\otimes m}$, consider
$\vph(\theta)=\vph(\theta_1,\theta_2,\ldots,\theta_m)$ as a function of $\theta_k$ alone and take this function as the
initial condition for equation~\eqref{1} at initial time~$\tau$, that is, we take
$\psi(\theta_k)=\vph(\theta_k,\wh\theta)$.
(We slightly abuse notation by writing $\vph(\theta_k,\wh\theta)$ for $\vph(\theta)$.)
Denoting the resulting solution by $y(t,\wh\theta)$ for $t\ge\tau$, we have that
$[\wt U_k(t,\tau)\vph](\theta)=y(t+\theta_k,\wh\theta)$, and so by~\eqref{m1} we have
\be{11}
[\wt U_k(\tau+\eta,\tau)\vph](\theta)=
\left\{\begin{array}{lcl}
\vph(\eta+\theta_k,\wh\theta), & \hbox{ for} &\bb -1\le\theta_k\le -\eta,\\
\\
\ds{\vph(0,\wh\theta)
-\int^{\eta+\theta_k}_0b(\tau+s)\vph(s-1,\wh\theta)\:ds,}
& \hbox{ for} &\bb -\eta\le\theta_k\le 0,
\end{array}\right.
\ee{11}
again assuming $0\le\eta\le 1$.

In Proposition~4.2 and below, we denote
\be{brho}
b^\rho(s)=b(\rho+s)
\ee{brho}
for ease of notation. The superscript notation here is formally distinguished from the
subscript notation $x_t(\theta)=x(t+\theta)$ used earlier wherein
the argument $\theta$ was restricted to the interval $\theta\in[-1,0]$ and
$x_t$ was regarded as an element of $C([-1,0])$. No restriction is imposed upon
the argument $s$ of $b^\rho(s)$, and $b^\rho(\cdot)$ is not viewed as an element
of any particular space, but merely as a shorthand notation.
\vv

The following result is the key to proving the Positivity Theorem.
\vv

\noindent {\bf Proposition~4.2.} \bxx
Fix $m\ge 1$, and let $\tau\in\IR$ and $0<\eta\le 1$.
Fix any $\theta=(\theta_1,\theta_2,\ldots,\theta_m)\in T_m$ where $T_m$ is the set~\eqref{tri},
and let $a$ be any integer in the range $0\le a\le m$ satisfying
\be{thet}
\theta_a\le-\eta\le\theta_{a+1},
\ee{thet}
where $a=0$ is allowed in case $-\eta\le\theta_1$,
and $a=m$ is allowed in case $\theta_m\le -\eta$. Then for every $\vph\in X^{\wedge m}$ we have that
\be{uform}
\begin{array}{lcl}
[\wwtilde(\tau+\eta,\tau)\vph](\theta) &\bb = &\bb
\ds{(-1)^{am}\int_{\eta+\theta_{a+1}-1}^{\eta+\theta_{a+2}-1}\!\!\!\!\cdots\int_{\eta+\theta_{m-1}-1}^{\eta+\theta_{m}-1}
b^{\tau+1}(t_{1})\cdots b^{\tau+1}(t_{m-a-1})}\\
\\
&\bb &\bb
\qquad\qquad\times\:\:\vph(t_{1},\ldots,t_{m-a-1},\eta+\theta_1,\ldots,\eta+\theta_a,0)\:dt_{m-a-1}\cdots dt_{1}\\
\\
&\bb &\bb
\ds{+(-1)^{(a+1)m}\int_{-1}^{\eta+\theta_{a+1}-1}\int_{\eta+\theta_{a+1}-1}^{\eta+\theta_{a+2}-1}
\!\!\!\!\cdots\int_{\eta+\theta_{m-1}-1}^{\eta+\theta_{m}-1}
b^{\tau+1}(t_{1})\cdots \makebox[0pt][l]{$b^{\tau+1}(t_{m-a})$}}\\
\\
&\bb &\bb
\qquad\qquad\times\:\:\vph(t_{1},\ldots,t_{m-a},\eta+\theta_1,\ldots,\eta+\theta_a)\:dt_{m-a}\cdots dt_{1},
\end{array}
\ee{uform}
if $0\le a\le m-2$, where the terms $\eta+\theta_j$ in the arguments of $\vph$ are absent if $a=0$.
If $a=m-1$ then
\be{uformx}
\begin{array}{lcl}
[\wwtilde(\tau+\eta,\tau)\vph](\theta) &\bb = &\bb
\vph(\eta+\theta_1,\ldots,\eta+\theta_{m-1},0)\\
\\
&\bb &\bb
\ds{+(-1)^m\int_{-1}^{\eta+\theta_{m}-1}
b^{\tau+1}(t_{1})\vph(t_{1},\eta+\theta_1,\ldots,\eta+\theta_{m-1})\:dt_{1},}
\end{array}
\ee{uformx}
while
\be{uformz}
[\wwtilde(\tau+\eta,\tau)\vph](\theta)=\vph(\eta+\theta_1,\ldots,\eta+\theta_m)
\ee{uformz}
if $a=m$.
\exx
\vv

\noindent {\bf Remark.}
The integer $a$ in the statement of Proposition~4.2 need not be unique; indeed, this is the case if
$\theta_k=-\eta$ for some $k$, where either $a=k-1$ or $a=k$ could be taken.
Indeed, if $a$ is not unique then
any value of $a$ permitted by the statement of the proposition may be taken.
\vv

Proposition~4.2 will be proved in Section~6. However, assuming its validity,
we use it here to prove Theorem~4.1. (In Section~5 we shall also use Proposition~4.2 to prove
Proposition~5.3 and Corollary~5.4, which in turn will be used to prove
Theorem~5.1 and Proposition~5.2.)
\vv

\noindent {\bf Proof of Theorem~4.1.}
We prove the result only for equation~\eqref{1}. The result
for the general equation~\eqref{0} follows directly from the conjugacy~\eqref{conj}
and the positivity of the operators
$\Sig(t)^{\wedge m}$ and $[\Sig(t)^{\wedge m}]^{-1}$.

Due to the fact that $\wwtilde(\tau+\eta,\tau)$ is a linear process, it is enough
to prove the theorem in the case that $0<\eta\le 1$. Taking such $\eta$, we
assume that $(-1)^mb(t)\ge 0$ almost everywhere in $[\tau,\tau+\eta]$.
With $\vph\in K_m\subseteq X^{\wedge m}$ and $\theta\in T_m$ fixed, and with $a$ as in the statement of Proposition~4.2,
consider the formulas~\eqref{uform}, \eqref{uformx}, and~\eqref{uformz} in that
result for $[\wwtilde(\tau+\eta,\tau)\vph](\theta)$.
Note that
$$
(t_{1},\ldots,t_{m-a-1},\eta+\theta_1,\ldots,\eta+\theta_a,0)\in T_m,\qquad
(t_{1},\ldots,t_{m-a},\eta+\theta_1,\ldots,\eta+\theta_a)\in T_m,
$$
both hold for the arguments of $\vph$ in these formulas, in particular because
$\eta+\theta_m-1\le\eta+\theta_1$ and $\eta+\theta_a\le 0$. Therefore $\vph$
evaluated at these points is nonnegative. Thus if $m$ is even, so $b(t)\ge 0$,
it is immediate from these formulas that $[\wwtilde(\tau+\eta,\tau)\vph](\theta)\ge 0$. If $m$ is odd,
so $b(t)\le 0$, the same conclusion holds after noting that
$(-1)^{am}$=$(-1)^{m-a-1}$ and $(-1)^{(a+1)m}$=$(-1)^{m-a}$.
In either case one concludes that $\wwtilde(\tau+\eta,\tau)\vph\in K_m$,
as desired.~\fp
\vv

We believe it is instructive
to verify Proposition~4.2 in the simplest nontrivial case of $m=2$ and $\eta=1$.
First, by equation~\eqref{11} we have that
$$
\begin{array}{lcl}
[\wt U_1(\tau+1,\tau)\vph](\theta) &\bb = &\bb
\ds{\vph(0,\theta_2)-\int_0^{1+\theta_1}b(\tau+s)\vph(s-1,\theta_2)\:ds,}\\
\\
{[}\wt U_2(\tau+1,\tau)\vph](\theta) &\bb = &\bb
\ds{\vph(\theta_1,0)-\int_0^{1+\theta_2}b(\tau+s)\vph(\theta_1,s-1)\:ds,}
\end{array}
$$
for every $\vph\in X^{\otimes 2}=C([-1,0]^2)$, where $\theta=(\theta_1,\theta_2)\in[-1,0]^2$.
We next compose these two formulas as in~\eqref{bigu}, substituting the second
into the first. Denoting
$\psi(\theta)=[\wt U_2(\tau+1,\tau)\vph](\theta)$, we have that
\be{u2}
\begin{array}{lcl}
[\wt\UU(\tau+1,\tau)\vph](\theta) &\bb = &\bb
\ds{\psi(0,\theta_2)-\int_0^{1+\theta_1}b(\tau+s)\psi(s-1,\theta_2)\:ds}\\
\\
&\bb = &\bb
\ds{\vph(0,0)-\int_0^{1+\theta_2}b(\tau+s)\vph(0,s-1)\:ds
-\int_0^{1+\theta_1}b(\tau+s)\psi(s-1,\theta_2)\:ds}\\
\\
&\bb = &\bb
\ds{\vph(0,0)-\int_0^{1+\theta_2}b(\tau+s)\vph(0,s-1)\:ds}\\
\\
&\bb &\bb
\ds{-\int_0^{1+\theta_1}b(\tau+s)
\bigg(\vph(s-1,0)-\int_0^{1+\theta_2}b(\tau+r)\vph(s-1,r-1)\:dr\bigg)\:ds}\\
\\
&\bb = &\bb
\ds{\vph(0,0)-\int_0^{1+\theta_2}b(\tau+s)\vph(0,s-1)\:ds
-\int_0^{1+\theta_1}b(\tau+s)\vph(s-1,0)\:ds}\\
\\
&\bb &\bb
\ds{+\int_0^{1+\theta_1}\int_0^{1+\theta_2}b(\tau+s)b(\tau+r)\vph(s-1,r-1)\:dr\:ds.}
\end{array}
\ee{u2}
The formula occurring after the final equal sign in~\eqref{u2} simplifies in
the anti-symmetric case $\vph\in X^{\wedge 2}$, that is,
where $\vph(\theta_1,\theta_2)\equiv-\vph(\theta_2,\theta_1)$
holds identically. For such $\vph$ we have that $\vph(0,0)=0$. Moreover, we have that
$$
-\int_0^{1+\theta_2}b(\tau+s)\vph(0,s-1)\:ds
-\int_0^{1+\theta_1}b(\tau+s)\vph(s-1,0)\:ds
=\int_{1+\theta_1}^{1+\theta_2}b(\tau+s)\vph(s-1,0)\:ds,
$$
and also that
$$
\begin{array}{l}
\ds{\int_0^{1+\theta_1}\int_0^{1+\theta_2}b(\tau+s)b(\tau+r)\vph(s-1,r-1)\:dr\:ds}\\
\\
\qquad\qquad
\ds{=\int_0^{1+\theta_1}\int_0^{1+\theta_1}b(\tau+s)b(\tau+r)\vph(s-1,r-1)\:dr\:ds}\\
\\
\qquad\qquad\:\:\:\:\:\:
\ds{+\int_0^{1+\theta_1}\int_{1+\theta_1}^{1+\theta_2}b(\tau+s)b(\tau+r)\vph(s-1,r-1)\:dr\:ds}\\
\\
\qquad\qquad
\ds{=\int_0^{1+\theta_1}\int_{1+\theta_1}^{1+\theta_2}b(\tau+s)b(\tau+r)\vph(s-1,r-1)\:dr\:ds,}
\end{array}
$$
where the integral taken over the square $[0,1+\theta_1]^2$ vanishes on account of the anti-symmetry.
With this we obtain from~\eqref{u2} that
\be{uu2}
\begin{array}{lcl}
[\wwtilde(\tau+1,\tau)\vph](\theta) &\bb = &\bb
\ds{\int_{1+\theta_1}^{1+\theta_2}b(\tau+s)\vph(s-1,0)\:ds}\\
\\
&\bb &\bb
\ds{+\int_0^{1+\theta_1}\int_{1+\theta_1}^{1+\theta_2}b(\tau+s)b(\tau+r)\vph(s-1,r-1)\:dr\:ds,}
\end{array}
\ee{uu2}
and one sees this formula coincides with~\eqref{uform}, where $a=0$ is taken.

To check that the operator $\wwtilde(\tau+1,\tau)$ is positive with respect to the cone $K_2$
under the condition $(-1)^mb(t)=b(t)\ge 0$,
take $\vph\in K_2$ and let $\theta\in T_2$, that is, $-1\le\theta_1\le\theta_2\le 0$. Then
$s-1\le \theta_2\le 0$ for the first term in~\eqref{uu2},
thus $(s-1,0)\in T_2$ and so $\vph(s-1,0)\ge 0$.
Similarly, for the second term in~\eqref{uu2} we have that $s-1\le \theta_1\le r-1$, and so $(s-1,r-1)\in T_2$ and
$\vph(s-1,r-1)\ge 0$. It follows that
the expression in~\eqref{uu2} is nonnegative, so
$\wwtilde(\tau+1,\tau)\vph\in K_2$, as desired.

We end this section by describing several properties of the above
linear processes. The first is a well-known compactness property
of $U(\tau+\eta,\tau)$ for $\eta>0$.
Its significance is that in the case of a periodic process, some
power of the monodromy operator is compact, and so the Floquet spectrum
consists entirely of Floquet multipliers, and these are isolated
values $\lam\in\IC\setminus\{0\}$ of finite multiplicity which can only cluster
at $\lam=0$.
\vv

\noindent {\bf Proposition~4.3.} \bxx
If $\eta\ge 1$ then the solution operator $U(\tau+\eta,\tau)$ for equation~\eqref{0} is compact.
More generally, if $\eta\ge\frac{1}{n}$ for some $n\ge 1$ then the $n^{\xx th}$ power
$U(\tau+\eta,\tau)^n$ is compact.
\exx
\vv

\noindent {\bf Proof.}
It is enough to work with the transformed equation~\eqref{1} and consider the
operator $\wt U(\tau+\eta,\tau)$.
The operator $\wt U(\tau+1,\tau)$ is easily seen from~\eqref{m1} to be compact, being the sum of the rank-one operator $\vph(0)$
and an integral operator, and thus if $\eta\ge 1$ then the operator
$\wt U(\tau+\eta,\tau)=\wt U(\tau+\eta-1,\tau+1)\wt U(\tau+1,\tau)$
is also compact.

Now suppose that $\eta\ge\frac{1}{n}$ for some $n\ge 1$. Define a new function $\wh b(t)$ by
setting $\wh b(t)=b(t)$ for $\tau\le t<\tau+\eta$, and extending it periodically so that $\wh b(t+\eta)=\wh b(t)$ for every
$t\in \IR$. Let $\wh U(\sig+\zeta,\sig)$ denote the linear process
associated to equation~\eqref{1} but with $\wh b(t)$ replacing $b(t)$, where
$\sig\in\IR$ and $\zeta\ge 0$ are general arguments.
Note that $\wh U(\tau+\eta,\tau)=\wt U(\tau+\eta,\tau)$
for our specific $\tau$ and $\eta$, as these operators only involve the range where $\wh b(t)$ and $b(t)$ agree.
Also note that $\wh U(\sig+\zeta+\eta,\sig+\eta)=\wh U(\sig+\zeta,\sig)$
for every $\sig$ and $\zeta$, due to the $\eta$-periodicity of  $\wh b(t)$. From this it follows that
$\wt U(\tau+\eta,\tau)^n=\wh U(\tau+\eta,\tau)^n=\wh U(\tau+n\eta,\tau)$, which is a compact operator
by our earlier remarks as $n\eta\ge 1$.~\fp
\vv

The next result concerns one-to-oneness of the above linear processes.
\vv

\noindent {\bf Proposition~4.4.} \bxx
Assume, for some $\tau\in\IR$
and $\eta>0$, that $\beta(t)\ne 0$ for almost every $t\in[\tau,\tau+\eta]$ for the
delay coefficient in equation~\eqref{0}. Let $m\ge 1$.
Then the operators $\UU(\tau+\eta,\tau)\in\call(X^{\otimes m})$ and  thus $\WW(\tau+\eta,\tau)\in\call(X^{\wedge m})$
are one-to-one.
\exx
\vv

\noindent {\bf Proof.}
It is enough to consider equation~\eqref{1}, where $b(t)\ne 0$ for almost every $t\in[\tau,\tau+\eta]$, and by
Proposition~2.5 it is enough to show that $\wt U(\tau+\eta,\tau)\in\call(X)$ is one-to-one. Also,
we may assume without loss that $0<\eta\le 1$ due to the fact that $\wt U(t,\tau)$ is a linear process.
Assume for some such $\eta$ that $\wt U(\tau+\eta,\tau)\vph=0$ for  some $\vph\in X$. Then from~\eqref{m1}
we have that $\vph(\theta)=0$ for every $\theta\in[\eta-1,0]$ and that
$$
\int_0^\theta b(\tau+s)\vph(s-1)\:ds=0\qquad\mbox{for every }\theta\in[0,\eta].
$$
Differentiating the above integral shows that $b(\tau+s)\vph(s-1)=0$ for almost every
$s\in[0,\eta]$, and as $b(\tau+s)\ne 0$ for almost every such $s$, we conclude that
$\vph(s-1)=0$ for $s\in[0,\eta]$. Thus $\vph(\theta)=0$ for every $\theta\in[-1,0]$,
and this gives the result.~\fp

\section{Positivity and Floquet Theory}

In this section we describe some basic consequences of the Positivity Theorem
in the context of Floquet theory.
One main result of this section, Theorem~5.1, provides additional structure
to the Floquet multipliers and their associated eigenfunctions in the case
that the feedback coefficients are periodic and
the delay coefficient $\beta(t)$ is of constant sign. As noted, these results confirm conjectures and
supplement results of Sell and one of the authors~\cite{sell} which provided
partial information on the multipliers. We note that the results of Theorem~5.1 have already
been established for a class of cyclic systems of ordinary differential equations
with a signed feedback; see~\cite{smith}.

Assume the coefficients in equation~\eqref{0}
are $\gam$-periodic for some $\gam>0$, that is
\be{gper}
\alp(t+\gam)=\alp(t),\qquad
\bet(t+\gam)=\bet(t),
\ee{gper}
for almost every $t$, and recall the monodromy operator $M(0)=U(\gam,0)$
where $U(t,\tau)$ denotes the linear process on $X=C([-1,0])$ associated 
to equation~\eqref{0}. In what follows we always take the initial time as $\tau=0$,
so we write $M=M(0)$ for simplicity.
As some power of $M$ is compact, the Floquet
spectrum $\spec(M)\setminus\{0\}$ consists entirely of point spectrum
(Floquet multipliers), of which there are at most countably many, and each is isolated and of finite multiplicity.
We have the following result.
\vv

\noindent {\bf Theorem~5.1.} \bxx
Consider equation~\eqref{0} where $\alp:\IR\to\IR$ and $\bet:\IR\to\IR$
are locally integrable and $\gam$-periodic for some $\gam>0$, and so satisfy~\eqref{gper}
for almost every $t$. Also assume that
\be{pari}
(-1)^m\bet(t)\ge(-1)^m\bet_0>0
\ee{pari}
for almost every $t$, for some integer $m$ and some quantity $\bet_0$.
Then there are countably infinitely many Floquet multipliers $\{\lam_k\}_{k=1}^\infty$,
that is, nonzero elements of the spectrum of the monodromy operator. Further, if the multipliers are labelled
so that
\be{mults}
|\lam_1|\ge|\lam_2|\ge|\lam_3|\ge\cdots,
\ee{mults}
with repetitions according to algebraic multiplicity, then it is the case that
the strict inequality
\be{gap}
|\lam_k|>|\lam_{k+1}|
\ee{gap}
holds whenever $k-m$ is even. More precisely, in the notation of Proposition~5.6 below one has that
\be{in1}
\calj(|\lam_{k-1}|)=\calj(|\lam_k|)=k-1,\qquad
\lam_{k-1}\lam_k>0,
\ee{in1}
whenever $k-m$ is even with $k\ge 2$, and
\be{in2}
\calj(|\lam_1|)=0,\qquad
\lam_1>0,
\ee{in2}
when $k=1$ with $m$ odd.
In particular, the inequalities~\eqref{jdim} and~\eqref{jdim1} below are in fact equalities.
\exx
\vv

It will be useful to introduce some notation. With $M$ denoting
the monodromy operator for equation~\eqref{0} as above, for any $\lambda\in\IC\setminus\{0\}$ denote
$$
G_\lambda=\bigcup_{k=0}^\infty\ker(M-\lambda I)^k,
$$
which is the generalized eigenspace for $M$, and which is a
finite-dimensional subspaces of $X$ (and is just $\{0\}$ for all but countably many $\lambda$). Define for any $\rho>0$
\be{calg}
\calg_\rho=\re\bigg(\bigoplus_{|\lambda|=\rho}G_\lambda\bigg),\qquad
\calh_\rho=\re\bigg(\bigoplus_{|\lambda|\ge\rho}G_\lambda\bigg),
\ee{calg}
which are the real parts of the spans of the generalized eigenspaces of eigenvalues
with norms as indicated, and which are both finite-dimensional. Also define
\be{p}
P=\{\rho>0\:|\:\rho=|\lambda|\hbox{ for some }\lambda\in\spec(M)\}.
\ee{p}
It is easy to see that due to the gap~\eqref{gap}, Theorem~5.1 implies that the space $\calh_{|\lambda_m|}$ has dimension $m$
for $m$ of the parity as in~\eqref{pari}.

The next result establishes positivity of the eigenfunction of the leading eigenvalue
of $M^{\wedge m}=\WW(\tau,0)$. In the constant coefficient case~\eqref{cc2} the solutions may be taken as
simply $x^j(t)=e^{\zeta_jt}$, where $\{\zeta_j\}_{j=1}^\infty$ are the roots of the characteristic
equation $\zeta=-\alpha_0-\beta_0e^{-\zeta}$ ordered so that
\be{rez}
\re\zeta_1\ge\re\zeta_2\ge\re\zeta_3\ge\cdots.
\ee{rez}
But even in this case we know of no elementary proof of this result.
\vv

\noindent {\bf Proposition~5.2.} \bxx
In the setting of Theorem~5.1, let $x^j(t)$, for $1\le j\le m$, be linearly independent solutions
of equation~\eqref{0} with
$x^j_0\in\calh_{|\lambda_m|}$ for each $j$. Then there exists a constant $C$ such that
\be{gl}
C\det\left(\begin{array}{ccc}
x^1(t+\theta_1) &\bb \cdots &\bb x^1(t+\theta_m)\\
\vdots &\bb &\bb \vdots\\
x^m(t+\theta_1) &\bb \cdots &\bb x^m(t+\theta_m)
\end{array}\right)\ge 0,
\ee{gl}
for every $t\in\IR$ and $\theta=(\theta_1,\theta_2,\ldots,\theta_m)\in T_m$, and where
for every $t$ the left-hand side of~\eqref{gl} is not identically zero as a function of $\theta$.
\exx
\vv

We shall prove Theorem~5.1 and Proposition~5.2 at the end of this section, following the proofs of Proposition~5.3
and Corollary~5.4 below, on which they rely. All of these results ultimately depend on Proposition~4.2,
which will be proved later, in Section~6.

Proposition~5.3 compares the Floquet multipliers of two
systems whose coefficients are appropriately ordered, and Corollary~5.4
provides computable lower bounds for the modulus of the Floquet
multipliers.
\vv

\noindent {\bf Proposition~5.3.} \bxx
Consider equation~\eqref{0} and also the equation
\be{5a}
\dot x(t)=-\alp(t)x(t)-\wh\bet(t)x(t-1)
\ee{5a}
with the same coefficient $\alp(t)$ of $x(t)$ but a different delay coefficient $\wh\bet(t)$.
Assume that $\alp,\bet,\wh\bet:\IR\to\IR$ are locally integrable and $\gam$-periodic
for some $\gam>0$ (with the same period $\gam$ for all coefficients). Also assume that
$$
(-1)^m\bet(t)\ge(-1)^m\wh\bet(t)\ge 0
$$
for almost every $t$, for some integer $m$.
Let $\{\lam_k\}_{k=1}^\infty$ and $\{\wh\lam_k\}_{k=1}^\infty$ denote the Floquet multipliers
of equations~\eqref{0} and~\eqref{5a}, respectively, ordered as in~\eqref{mults}, with the convention that $\lam_k=0$
for $k>n$ if there are only a finite number, $n$, of multipliers, and similarly for $\wh\lam_k$.
Then
\be{linq2}
|\lam_1\lam_2\cdots\lam_k|\ge|\wh\lam_1\wh\lam_2\cdots\wh\lam_k|
\ee{linq2}
for every $k\ge 1$ for which $k-m$ is even.
\exx
\vv

\noindent {\bf Corollary~5.4.} \bxx
Consider equation~\eqref{0} where the coefficients
$\alp,\bet:\IR\to\IR$ are locally integrable and $\gam$-periodic
for some $\gam>0$. Let $b(t)$ be given by~\eqref{xyz} and suppose that
$$
(-1)^mb(t)\ge(-1)^mb_0>0
$$
for almost every $t$, for some integer $m$.
Let $\{\lam_k\}_{k=1}^\infty$ denote the Floquet multipliers of equation~\eqref{0}, ordered as in~\eqref{mults},
and let $\{\zeta_k\}_{k=1}^\infty$ be the
roots of the characteristic equation $\zeta=-b_0e^{-\zeta}$, ordered as in~\eqref{rez}.
Also let
\be{rez2}
\alpha_0=\frac{1}{\gamma}\int_0^\gamma\alpha(s)\:ds,\qquad
Q=\exp\bigg(\int_0^\gamma |b(s)|\:ds\bigg).
\ee{rez2}
Then
$$
|\lambda_k|\ge Q^{-(k-1)}\exp\bigg(-\gamma\alpha_0+\gamma\sum_{j=1}^k\re\zeta_j\bigg)
$$
if $k-m$ is even. In particular, there are infinitely many Floquet multipliers $\lambda_k\ne 0$.
\exx
\vv

The theory of so-called {\bf lap numbers} will play a crucial role in the
proof of Theorem~5.1. We follow the treatment of~\cite{sell} in which
the lap number of a solution is the value
of the function $V^\pm$ defined below. Roughly, the
lap number counts the number of oscillations of a solution in a given delay interval.
More precisely, define a function $\sch:X\to\{0,1,2,\ldots,\infty\}$,
where as always $X=C[-1,0]$, by
$$
\begin{array}{lcl}
\sch(\vph) &\bb = &\bb
\sup\{k\ge 1\:|\:\hbox{there exist }-1<\theta_0<\theta_1<\cdots<\theta_k<0\\
\\
&\bb &\bb
\hbox{ such that }\vph(\theta_{i-1})\vph(\theta_i)<0\hbox{ for }1\le i\le k\},
\end{array}
$$
which is the number of sign changes of $\vph$ in $[-1,0]$. By convention, $\sch(\vph)=0$
if either $\vph(\theta)\ge 0$ or $\vph(\theta)\le 0$ for every $\theta\in[-1,0]$.
Now define two function, $V^-:X\setminus\{0\}\to\{1,3,5,\ldots,\infty\}$ and
$V^+:X\setminus\{0\}\to\{0,2,4,\ldots,\infty\}$, by
$$
\begin{array}{lcl}
V^-(\vph) &\bb = &\bb \left\{\begin{array}{lcl}
\sch(\vph), & \hbox{ for} &\bb \sch(\vph)\hbox{ odd or infinite},\\
\\
\sch(\vph)+1, & \hbox{ for} &\bb \sch(\vph)\hbox{ even},
\end{array}\right.\\
\\
V^+(\vph) &\bb = &\bb \left\{\begin{array}{lcl}
\sch(\vph), & \hbox{ for} &\bb \sch(\vph)\hbox{ even or infinite},\\
\\
\sch(\vph)+1, & \hbox{ for} &\bb \sch(\vph)\hbox{ odd}.
\end{array}\right.
\end{array}
$$
The significance of these functions is that they are Lyapunov functions for
the system~\eqref{0}. In particular,
the function $V^-$ is used in the case of negative feedback ($m$ even in the statement of Theorem~5.1)
while the function $V^+$ is used in the case of positive feedback ($m$ odd).
The following result from~\cite{sell} (specifically, for the case $N=0$
in the notation of~\cite{sell}) states some basic properties.
\vv

\noindent {\bf Proposition~5.5 (see~\cite[Theorems~2{.}1 and~2{.}2]{sell}).} \bxx
Consider equation~\eqref{0} where $\alp,\bet:\IR\to\IR$ are locally integrable and
$(-1)^m\bet(t)\ge 0$ for almost every $t\in[\tau,\tau+\eta]$, for some $\tau\in\IR$ and $\eta>0$.
Suppose that $x(t)$ is a solution of~\eqref{0} for $t\ge\tau$. Then the following are true
where the sign~$\pm$ is such that $\pm(-1)^m=-1$.
\begin{itemize}
\item[{(1)}] $V^\pm(x_t)$ is nonincreasing for $t\in[\tau,\tau+\eta]$ provided that $x_{\tau+\eta}\ne 0$;
that is, $V^\pm(x_{t_1})\ge V^\pm(x_{t_2})$ for $t_1,t_2\in[\tau,\tau+\eta]$ with $t_1\le t_2$.
\item[{(2)}] Suppose in addition that the strict inequality
$(-1)^m\beta(t)>0$ holds for almost every $t\in[\tau,\tau+\eta]$ and that $\eta\ge 3$, where
again $x_{\tau+\eta}\ne 0$. Also suppose that $V^\pm(x_t)$ is constant and finite throughout $[\tau,\tau+\eta]$, that is,
$V^\pm(x_t)=V^\pm(x_\tau)<\infty$ for every $t\in[\tau,\tau+\eta]$. Then
\be{r2}
(x(t),x(t-1))\ne (0,0)\in\IR^2
\ee{r2}
for every $t\in[\tau+3,\tau+\eta]$.
\end{itemize}
\exx
\vv

\noindent {\bf Remark.}
If the strict inequality $(-1)^m\beta(t)>0$ holds for almost every $t\in\IR$, and $x(t)$ is a solution with $x_t\ne 0$
and $V^\pm(x_t)=J<\infty$ for some $J$, for every $t\in\IR$, then~\eqref{r2} holds for every $t\in\IR$.
\vv

In the case of periodic coefficients, there is an intimate connection between lap numbers and
Floquet solutions, as the following result shows.
\vv

\noindent {\bf Proposition~5.6 (see~\cite[Theorem~3{.}1]{sell}).} \bxx
Assume the conditions of Theorem~5.1, and recall the set $P$ in~\eqref{p}, where $M$ is the monodromy operator
associated with equation~\eqref{0}. Also assume that $m$ is even (negative feedback).
Then there exists a function $\calj:P\to \{1,3,5,\ldots\}$ such that if $\rho\in P$, then
$$
V^-(\vph)=\calj(\rho)\qquad\hbox{for every }\vph\in\calg_\rho\setminus\{0\}.
$$
Further,
$$
\calj(\rho_1)\ge\calj(\rho_2)\qquad\hbox{if }\rho_1\le\rho_2
$$
for $\rho_1,\rho_2\in P$. Finally,
\be{jdim}
\dim\bigg(\bigoplus_{\calj(\rho)=J} \calg_\rho\bigg)\le 2
\ee{jdim}
for any $J\in\{1,3,5,\ldots\}$.

In the case that $m$ is odd (positive feedback), the analogous results hold
with $V^+$ in place of $V^-$, except that $\calj:P\to \{0,2,4,\ldots\}$ and
$J\in\{0,2,4,\ldots\}$ in~\eqref{jdim}, and that
\be{jdim1}
\dim\bigg(\bigoplus_{\calj(\rho)=0} \calg_\rho\bigg)\le 1
\ee{jdim1}
in the case $J=0$.
\exx
\vv

\noindent {\bf Remark.}
Proposition~5.6 implies that every Floquet multiplier $\lambda\in\spec(M)\setminus\{0\}$
is associated to an integer $\calj(|\lambda|)$ of a specific parity (odd or even), which gives
the lap number of all corresponding generalized eigensolutions.
Note that $\calj(|\lambda|)$ depends only on the modulus of $\lambda$. Now
suppose that the Floquet multipliers
are enumerated as in~\eqref{mults} in the statement of Theorem~5.1. Then
the left-hand side of the inequality~\eqref{jdim}, namely the dimension
of the indicated subspace, is precisely the number of integers $k$
for which $\calj(|\lambda_k|)=J$, and similarly for~\eqref{jdim1} for $J=0$.
One sees easily from this fact that in order to prove Theorem~5.1, it suffices to prove
that the inequalities in~\eqref{jdim} and~\eqref{jdim1} are in fact equalities, and as well to establish
the inequalities in~\eqref{in1} and~\eqref{in2}.

In the case of constant coefficients, namely the equation
$\dot x(t)=-\alpha_0x(t)-\beta_0x(t-1)$,
it is known (see for example~\cite[Theorem~3{.}2]{sell}
and also~\cite{morse}) that~\eqref{jdim} and~\eqref{jdim1} are equalities. It is also the case
that~\eqref{in1} and~\eqref{in2} hold; this follows easily from~\cite[Corollary~3{.}3]{sell} and the
fact that the multipliers are given by $\lambda=e^{\gamma\zeta}$ for roots $\zeta=-\alpha_0-\beta_0e^{-\zeta}$
of the characteristic equation.
\vv

\noindent {\bf Remark.}
Proving that~\eqref{jdim} and~\eqref{jdim1}
are equalities is equivalent to proving that
$\calj(|\lambda_k|)=k-1$ for every $k$ for which $k-m$ is even.
(Note that it is always the case that $\calj(|\lambda_k|)\ge k-1$ for such $k$.)
Indeed, if $\calj(|\lambda_{k_0}|)=k_0-1$ holds for some $k_0$ for which $k_0-m$ is even, then
\be{con}
\begin{array}{lcl}
\calj(|\lam_{1}|) = 0 &\bb \qquad &\bb \mbox{if }m\mbox{ is odd,}\\
\\
\calj(|\lam_{k-1}|)=\calj(|\lam_{k}|) = k-1 &\bb &\bb \mbox{for } 2\le k\le k_0\mbox{ with }k-m\mbox{ even,}\\
\\
\calj(|\lam_{k_0+1}|)\ge k_0+1.
\end{array}
\ee{con}
It follows immediately from~\eqref{con} that
$$
\begin{array}{ll}
|\lam_{k_0-2}|>|\lam_{k_0-1}|\ge|\lam_{k_0}|>|\lam_{k_0+1}| \qquad &\bb \hbox{if }k_0\ge 3,\\
\\
|\lam_{1}|\ge|\lam_{2}|>|\lam_{3}| \qquad &\bb \hbox{if }k_0=2,\\
\\
|\lam_{1}|>|\lam_{2}| \qquad &\bb \hbox{if }k_0=1.
\end{array}
$$
In any case the quantity $\lambda_{k_0-1}\lambda_{k_0}$ (if $k_0\ge 2$) must be real,
as either both $\lambda_{k_0-1}$ and $\lambda_{k_0}$ are real, or else both
are complex with $\lambda_{k_0-1}$=$\overline\lambda_{k_0}$. Similarly $\lambda_1$
must be real if $k_0=1$.

These observations, and those of the previous remark,
follow in an elementary fashion from Proposition~5.6, essentially by a pigeon-hole argument.
\vv

Before proving Proposition~5.3 and Corollary~5.4 we need the following result. Although it is
relatively well-known, we provide a proof for the convenience of the reader.
Recall here the definitions of reproducing and normal cones from Section~4.
\vv

\noindent {\bf Proposition~5.7.} \bxx
Let $Y$ be a Banach space and $K\subseteq Y$ a cone which is both
reproducing and normal. Suppose $A,B\in\call(Y)$ satisfy
$0\le A\le B$, that is, both $A$ and $B-A$ are positive operators. Then
$$
r(A)\le r(B)
$$
for the spectral radii of these operators.
\exx
\vv

\noindent {\bf Proof.}
First define
$$
\tr y\tr=\inf_{y=u-v\atop u,v\in K}(\|u\|+\|v\|)
$$
for any $y\in K$. This quantity is well-defined as $K$ is reproducing.
It is known~(see, for example,~\cite{deim} and~\cite{schwlf}) that $\tr\cdot\tr$ is a norm on $Y$ which is
equivalent to $\|\cdot\|$, and for which $\tr u\tr=\|u\|$ for every $u\in K$.
Next, for any operator $A\in\call(Y)$ which is positive, $A\ge 0$, define
$$
\tr A\tr_K=\sup_{y\in K\atop \tr y\tr=1}\tr Ay\tr.
$$
(See in particular~\cite{kr} for more information about this and related quantities.)

We claim that $\tr A\tr_K=\tr A\tr$ where $\tr A\tr$ denotes the usual operator
norm taken with respect to the new norm $\tr\cdot\tr$ on $Y$. Clearly
$\tr A\tr_K\le\tr A\tr$ holds. On the other hand, given any $y\in Y$ and $\eps>0$,
there exist $u,v\in K$ such that $y=u-v$ and
$$
\tr y\tr\ge\|u\|+\|v\|-\eps=\tr u\tr+\tr v\tr-\eps.
$$
As $Au,Av\in K$ it follows that
$$
\tr Ay\tr\le\|Au\|+\|Av\|=\tr Au\tr+\tr Av\tr\le\tr A\tr_K(\tr u\tr+\tr v\tr)\le\tr A\tr_K(\tr y\tr+\eps),
$$
and thus $\tr Ay\tr\le\tr A\tr_K\tr y\tr$ as $\eps$ is arbitrary. From this one concludes
directly that $\tr A\tr\le\tr A\tr_K$, to give the claim.

Now assume that $A$ and $B$ are as in the statement of the proposition. Then for any $y\in K$
with $\tr y\tr=1$ and integer $n\ge 0$, we have that $0\le A^ny\le B^ny$ and thus
$$
\tr A^ny\tr=\|A^ny\|\le C\|B^ny\|=C\tr B^ny\tr\le C\tr B^n\tr,
$$
where the normality of $K$ is used. Taking the supremum over all such $y$ gives
$$
\tr A^n\tr=\tr A^n\tr_K=\sup_{y\in K\atop \tr y\tr=1}\tr A^ny\tr\le C\tr B^n\tr.
$$
Thus $\tr A^n\tr^{1/n}\le C^{1/n}\tr B^n\tr^{1/n}$, and letting $n\to\infty$ gives
$r(A)\le r(B)$, as desired.~\fp
\vv

\noindent {\bf Proof of Proposition~5.3.}
We again work with the transformed equations, noting that the same transformation~\eqref{cov} applies
to both equations~\eqref{0} and~\eqref{5a}.
The transformed equations have the form~\eqref{1}
with a coefficient $b(t)$ and a similar equation with a coefficient denoted $\wh b(t)$, where
\be{5b}
(-1)^mb(t)\ge(-1)^m\wh b(t)\ge 0.
\ee{5b}
Note that the transformation~\eqref{cov} alters the Floquet multipliers; to be precise, if $\{\lam_k\}_{k=1}^\infty$ denotes
the multipliers of~\eqref{0}, then the multipliers $\{\wt\lam_k\}_{k=1}^\infty$ of the transformed system~\eqref{1}
are given by $\wt\lam_k=e^{\gam\alp_0}\lam_k$, with $\alp_0$ as in~\eqref{rez2}.
This follows directly from the fact that $\mu(t+\gam)=e^{\gam\alp_0}\mu(t)$ in~\eqref{cov}. However,
as the same transformation~\eqref{cov} is also applied to~\eqref{5a}, it suffices to prove the result for
the two transformed systems with coefficients satisfying~\eqref{5b}.

In this proof we use a tilde~$\wt{{\ }{\ }}$ to denote objects associated to the equation with $b(t)$,
and a hat~$\wh{{\ }{\ }}$ to denote objects associated to the equation with $\wh b(t)$.
Also, without loss we may assume that $k=m$.

By~\eqref{bigw} we have $\wwtilde(\gam,0)=\wt U(\gam,0)^{\wedge m}=\mtilde^{\wedge m}$, and by
by Proposition~4.3 some power of $\mtilde$ is compact. Thus all points of $\spec(\mtilde)\setminus\{0\}$
(the Floquet multipliers) are isolated points with finite multiplicity.
By Proposition~2.4 and the remarks following it, the eigenvalues of $\mtilde^{\wedge m}$ have
the form~\eqref{lj} where $1\le j_1<j_2<\cdots<j_m$ with $\wt\lambda_k$ denoting the eigenvalues of $\mtilde$, and thus the quantity
$|\wt\lam_1\wt\lam_2\cdots\wt\lam_m|$ equals the spectral radius $r(\mtilde^{\wedge m})$.
Therefore, in order to establish~\eqref{linq2} we must show that $r(\mtilde^{\wedge m})\ge r(\mhat^{\wedge m})$.
We claim in fact that
$$
r(\wwtilde(\tau+\eta,\tau))\ge r(\wwhat(\tau+\eta,\tau))
$$
holds for every $\tau\in\IR$ and $\eta\ge 0$.

By Theorem~4.1 the operators $\wwtilde(\tau+\eta,\tau)$ and $\wwhat(\tau+\eta,\tau)$
are positive with respect to the cone $K_m$, which we have noted is both
reproducing and normal. Thus by Proposition~5.7 it suffices to show that
$$
\wwtilde(\tau+\eta,\tau)\ge \wwhat(\tau+\eta,\tau).
$$
The proof of this is very similar to the proof of Theorem~4.1, using the formulas~\eqref{uform},
\eqref{uformx}, and~\eqref{uformz} with $b(t)$ and with $\wh b(t)$,
and relying on the inequality~\eqref{5b}. We omit the details.~\fp
\vv

\noindent {\bf Proof of Corollary~5.4.}
The change of variables~\eqref{cov} converts equation~\eqref{0} into equation~\eqref{1},
however, the Floquet multipliers are changed. As noted in the proof above,
with $\{\lambda_k\}_{k=1}^\infty$ denoting the multipliers
of~\eqref{0} and $\{\wt\lambda_k\}_{k=1}^\infty$ denoting the multipliers of~\eqref{1}, we have $\wt\lambda_k=e^{\gamma\alp_0}\lambda_k$
with $\alp_0$ as in~\eqref{rez2}.

Now applying Proposition~5.3 to equation~\eqref{1}, we obtain
$$
|\wt\lambda_1\wt\lambda_2\cdots\wt\lambda_k|\ge\exp\bigg(\gamma\sum_{j=1}^k\re\zeta_j\bigg)
$$
provided $k-m$ is even, since the Floquet multipliers $\{\wh\lambda_j\}_{j=1}^\infty$ of the equation $\dot y(t)=-b_0y(t-1)$
are simply $\wh\lambda_j=e^{\gamma\zeta_j}$. Further, we have for the norm of the monodromy operator
for equation~\eqref{1} that $\|\mtilde\|\le Q$ by a simple Gronwall argument, and thus $|\wt\lambda_j|\le Q$ for every $j$.
Therefore
$$
|\lambda_k|=e^{-\gamma\alpha_0}|\wt\lambda_k|\ge e^{-\gamma\alpha_0}Q^{-(k-1)}|\wt\lambda_1\wt\lambda_2\cdots\wt\lambda_k|
\ge Q^{-(k-1)}\exp\bigg(-\gamma\alpha_0+\gamma\sum_{j=1}^k\re\zeta_j\bigg),
$$
for $k-m$ even, as desired.~\fp
\vv

In the proof of Theorem~5.1 we follow
the approach of~\cite{smith}, where a cyclic system of ordinary differential equations was
considered. In particular, the proof will involve a homotopy wherein we consider the system
\be{hom}
\begin{array}{l}
\dot x(t)=-\alp_\kap(t) x(t)-\bet_\kappa(t)x(t-1),\\
\\
\alp_\kap(t)=\kap\alp(t),\qquad
\beta_\kappa(t)=\kappa \bet(t)+(1-\kappa)\bet_0,
\end{array}
\ee{hom}
where $0\le\kappa\le 1$. Denoting the monodromy operator for equation~\eqref{hom} by $M_\kappa=U_\kappa(\gamma,0)$,
where $U_\kappa(t,\tau)$ denotes the linear process for equation~\eqref{hom},
in a standard fashion one easily obtains a common bound
\be{gron}
\|M_\kappa\|\le\exp\bigg(\int_0^\gamma |\alp_\kap(t)|+|\bet_\kap(t)|\:dt\bigg)
\le\exp\bigg(\int_0^\gamma |\alp(t)|+|\bet(t)|\:dt\bigg),
\ee{gron}
independent of $\kappa\in[0,1]$, via Gronwall's inequality, since $|\bet_\kappa(t)|\le|\bet(t)|$ holds.
We shall need to obtain bounds for
the norm of the inverse of $M_\kappa$ restricted to certain subspaces (eigenspaces), and to
this end several lemmas are required.

In the following lemma
the space ${{E}}$ is of course isomorphic to $\IR^n$, but there is no assumption made about the norm $\|\cdot\|$.
The crucial point is that the constant $C$ depends only on the integer $n$, and on $Q$ and $B$, but not on the norm
with which ${{E}}$ is endowed.
\vv

\noindent {\bf Lemma~5.8.} \bxx
Given a positive integer $n$ and positive quantities $Q$ and $B$, there exists a constant $C$ such that the following holds.
If $({{E}},\|\cdot\|)$ is any $n$-dimensional normed linear space and $L:{{E}}\to {{E}}$ is an invertible linear
transformation with norm satisfying $\|L\|\le Q$ and all of whose eigenvalues $\lambda$ satisfy $|\lambda|\ge B$,
then $\|L^{-1}\|\le C$ for the norm of the inverse. Here the norms of the operators are those inherited in the usual fashion
from the norm on ${{E}}$.
\exx
\vv

\noindent {\bf Proof.}
Let $({{E}},\|\cdot\|)$ be any $n$-dimensional normed linear space.
We first claim there exists a basis $\{v_k\}_{k=1}^n$ for ${{E}}$ with the property that
\be{vprop}
\|v_k\|=1,\qquad
\dist(v_k,{{E}}_{k-1})=1,
\ee{vprop}
for $1\le k\le n$, where ${{E}}_k=\spann\{v_1,v_2,\ldots,v_k\}$ for every $k$ (with ${{E}}_0=\{0\}$).
Indeed, such a basis can be constructed inductively. To begin, let $v_1\in {{E}}$ be any vector with $\|v_1\|=1$.
Assuming that vectors $v_j$ have been constructed for $1\le j\le k$
so that~\eqref{vprop} holds, where $k<n$, we construct $v_{k+1}$ as follows. Take any vector $x\in {{E}}\setminus {{E}}_k$
and let $y\in {{E}}_k$ be any point in ${{E}}_k$ which minimizes the distance $\|x-y\|$, that is, for which
$$
\|x-y\|=\dist(x,{{E}}_k).
$$
Let $v_{k+1}=(x-y)/\|x-y\|$. Clearly $\|v_{k+1}\|=1$, and one easily checks that $\dist(v_{k+1},{{E}}_k)=1$. This completes the induction.

With the basis $\{v_k\}_{k=1}^n$ fixed as above, we construct for each $k$ a linear functional $\xi_k:{{E}}\to\IR$
as follows. First set $\xi_k(v_j)=0$ for $1\le j\le k-1$ and $\xi_k(v_k)=1$, then extend
$\xi_k$ linearly to all of ${{E}}_k$. One checks using~\eqref{vprop} that $\|\xi_k\|=1$ for the norm of this functional. Then
by the Hahn-Banach Theorem, extend $\xi_k$ as a linear functional to all of ${{E}}$ with the same norm $\|\xi_k\|=1$.

We next claim that
\be{eqnrm}
\sum_{j=1}^k|a_j|\le (2^k-1)\bigg\|\sum_{j=1}^ka_jv_j\bigg\|
\ee{eqnrm}
for $1\le k\le n$ and any choice of real numbers $a_j$. We prove~\eqref{eqnrm} by induction on $k$, first noting
that it holds trivially for $k=1$. Assume therefore that~\eqref{eqnrm} holds for some $k$ where $1\le k\le n-1$. Then
taking any $x\in {{E}}_{k+1}$, denote
$$
x=\sum_{j=1}^{k+1}a_jv_j,\qquad y=\sum_{j=1}^ka_jv_j,
$$
and observe that
$$
|a_{k+1}|=\dist(a_{k+1}v_{k+1},{{E}}_k)\le\|a_{k+1}v_{k+1}+y\|=\|x\|
$$
and that
$$
\|y\|=\|x-a_{k+1}v_{k+1}\|=\|x-\xi_{k+1}(x)v_{k+1}\|\le\|x\|+\|\xi_{k+1}\|\|x\|\|v_{k+1}\|=2\|x\|.
$$
Therefore,
$$
\sum_{j=1}^{k+1}|a_j|\le\sum_{j=1}^k|a_j|+\|x\|\le(2^k-1)\|y\|+\|x\|\le2(2^k-1)\|x\|+\|x\|=(2^{k+1}-1)\|x\|,
$$
where the induction hypothesis~\eqref{eqnrm} was used in the second inequality. This completes the induction
and establishes~\eqref{eqnrm} for all $k$.

Now let $L:{{E}}\to {{E}}$ be any linear map and let $A$ be the matrix representation of $L$ with respect to
the above basis $\{v_k\}_{k=1}^n$. Let the space $\IR^n$ be endowed with the $l^1$ norm, that is,
$$
\|a\|=\sum_{j=1}^n|a_j|,\qquad a=(a_1,a_2,\ldots,a_n)\in\IR^n.
$$
For any vector $x\in {{E}}$ let $y=Lx$, and let $a,b\in\IR^n$ be such that
$$
x=\sum_{j=1}^na_jv_j,\qquad
y=Lx=\sum_{j=1}^nb_jv_j,\qquad
a=(a_1,a_2,\ldots,a_n),\qquad
b=Aa=(b_1,b_2,\ldots,b_n).
$$
Then
$$
(2^n-1)^{-1}\|Aa\|=
(2^n-1)^{-1}\sum_{j=1}^n|b_j|\le\|y\|=\|Lx\|\le\|L\|\|x\|\le\|L\|\sum_{j=1}^n|a_j|=\|L\|\|a\|
$$
where~\eqref{eqnrm} has been used, and therefore
\be{al}
\|A\|\le (2^n-1)\|L\|.
\ee{al}
Here $\|L\|$ and $\|A\|$ denote the operator norms inherited from the norms on ${{E}}$ and $\IR^n$, respectively.
Again using~\eqref{eqnrm}, we also have that
$$
\|Lx\|=\|y\|\le\sum_{j=1}^n|b_j|=\|Aa\|\le\|A\|\|a\|=\|A\|\sum_{j=1}^n|a_j|\le(2^n-1)\|A\|\|x\|
$$
and it follows that
\be{la}
\|L\|\le(2^n-1)\|A\|.
\ee{la}

To complete the proof of the lemma, let $Q$ and $B$ be positive constants and suppose
that $L$ is as in the statement of the lemma. With the basis $\{v_k\}_{k=1}^n$ and the matrix $A$
as above, we have from~\eqref{al} that
$\|A\|\le (2^n-1)Q$. Also, let $A^*$ denote the adjugate matrix of $A$, namely the transpose of the matrix of cofactors
of $A$. Each entry of $A^*$ is the determinant of an $(n-1)\times (n-1)$ submatrix of $A$.
Due to the choice of the $l^1$ norm on $\IR^n$, this submatrix has norm at most $\|A\|\le(2^n-1)Q$,
and thus each of its eigenvalues is at most $(2^n-1)Q$ in magnitude.
Thus each entry of $A^*$ has magnitude no greater than $(2^n-1)^{n-1}Q^{n-1}$, and so
$\|A^*\|\le n(2^n-1)^{n-1}Q^{n-1}$. Now $|\det A|\ge B^n$ since
all the eigenvalues $\lambda$ of $A$ satisfy
$|\lambda|\ge B$, and as $A^{-1}=(\det A)^{-1}A^*$, we have that
$\|A^{-1}\|\le n(2^n-1)^{n-1}Q^{n-1}B^{-n}$. Therefore, by~\eqref{la},
$$
\|L^{-1}\|\le(2^n-1)\|A^{-1}\|\le n(2^n-1)^nQ^{n-1}B^{-n}=C
$$
where the above formula serves as the definition of $C$. This completes the proof.~\fp
\vv

\noindent {\bf Lemma~5.9.} \bxx
Consider equation~\eqref{0} where $\alp,\bet:\IR\to\IR$ are $\gamma$-periodic
and where~\eqref{pari} holds for some $\bet_0$, for almost every $t$.
Also consider the homotopy system~\eqref{hom}, and let $M_\kappa=U_\kappa(\gamma,0)$
denote the monodromy operator, where $U_\kappa(t,\tau)$ denotes the linear process on $X=C[-1,0]$
associate to this system. Suppose for some sequence $\kappa_i\to\kappa_0$, where $\kappa_i,\kappa_0\in[0,1]$,
and for some integer $n$, that there exists a
sequence of $n$-dimensional subspaces ${E}_i\subseteq X$ for $i\ne 0$ for which
$$
M_{\kappa_i}{E}_i={E}_i
$$
for every $i$. Further suppose there is a positive lower bound
$$
\min\{|\lambda|\:|\:\lambda\in\spec(M_{\kappa_i}|_{{E}_i})\}\ge B>0
$$
for the spectra of the monodromy operators $M_{\kappa_i}$ restricted to ${E}_i$. Take any sequence $\vph^i\in{E}_i$
for which the norms $\|\vph^i\|$ are bounded, and let $x^i(t)$ denote the solution of equation~\eqref{hom}
with $x^i_0=\vph^i$ for $t\in\IR$, with $\kappa=\kappa_i$. Then there exists a subsequence $x^{i_j}(t)$ such that
$$
x^{i_j}(t)\to x^0(t)
$$
uniformly on compact intervals, where $x^0(t)$ is a solution of equation~\eqref{hom} with $\kappa=\kappa_0$.
\exx
\vv

\noindent {\bf Proof.}
By Lemma~5.8, and using the bound~\eqref{gron}, there exists a bound
$$
\|(M_{\kappa_i}|_{{E}_i})^{-1}\|\le C
$$
for the inverse of $M_{\kappa_i}$ restricted to ${E}_i$, which is independent of $i$.
It follows that for each $i$ the solution $x^i(t)$ enjoys the bounds
$$
\|x^i_{-k\gamma}\|=\|(M_{\kappa_i}|_{{E}_i})^{-k}\vph^i\|\le C^k\|\vph^i\|
$$
for integers $k\ge 0$. Therefore, by a Gronwall argument with initial condition
$x^i_{-k\gam}$ and initial time $-k\gam$, one has the uniform boundedness property
$$
|x^i(t)|\le Q^{2k}\|x^i_{-k\gam}\|\le Q^{2k}C^k\|\vph^i\|\qquad \hbox{for every }t\in[-k\gamma,k\gamma],
$$
for integers $k\ge 0$. Here $Q$ denotes the quantity in the right-hand side of~\eqref{gron}.
Further, the sequence of functions $x^i(t)$
is equicontinuous on $[-k\gam+1,k\gam]$, as one sees directly from the differential equation~\eqref{hom}.
From this, with an application of Ascoli's Theorem, the result follows.~\fp
\vv

The functions $V^\pm$, being integer-valued, have discontinuities. However, they exhibit the following
pseudocontinuity property for certain solutions of interest; see also~\cite[Lemma~8{.}1]{morse}.
\vv

\noindent {\bf Lemma~5.10.} \bxx
Consider the setting of Lemma~5.9 and the notation therein.
Suppose for some sequence $\kappa_i\to\kappa_0$, where $\kappa_i,\kappa_0\in[0,1]$,
that $x^i(t)$ is a solution of equation~\eqref{hom} for $t\in\IR$, with $\kappa=\kappa_i$.
Suppose further that $x^i(t)\to x^0(t)$ as $i\to\infty$ uniformly on compact intervals and that
$x^0(t)$ does not vanish identically. Finally suppose that there exists an integer $J$ such that
$V^\pm(x^0_t)=J$ for all $t\in\IR$, where $V^+$ is taken if $m$ is odd (positive feedback) and $V^-$
is taken if $m$ is even (negative feedback).
Then given any bounded interval $I$, it is the case that $V^\pm(x^i_t)=J$ for all $t\in I$, for all
sufficiently large $i$.
\exx
\vv

\noindent {\bf Proof.}
Certainly $x^0(t)$ is a solution of~\eqref{hom} for $\kap=\kap_0$.
The fact that $V^\pm(x^0_t)$ is finite and constant in $t$ implies, by~(2) of Proposition~5.5,
that $(x^0(t),x^0(t-1))\ne(0,0)$ for all $t\in\IR$. That is, whenever $x^0(t_0)=0$ for some $t_0$,
then $x^0(t_0-1)\ne 0$.
This implies, from the differential equation~\eqref{hom}, that
if $x^0(t_0)=0$ for some $t_0$, then for some $\delta\in\{-1,1\}$ and
some $C>0$ and $\eps>0$ we have
$$
\del\frac{d}{dt}\bigg(x^0(t)\exp\bigg(\int_0^t\alp(s)\:ds\bigg)\bigg)
=-\del\bet(t)x^0(t-1)\exp\bigg(\int_0^t\alp(s)\:ds\bigg)\ge C
$$
for almost every $t\in[t_0-\eps,t_0+\eps]$.
(Note that we have used~\eqref{pari}.)
In turn, this implies that $t_0$ must be an isolated zero of $x^0(t)$. In fact, one has for
an appropriate (possibly smaller) $C$ and $\eps$, that for every sufficiently large $i$
$$
\del\frac{d}{dt}\bigg(x^i(t)\exp\bigg(\int_0^t\alp(s)\:ds\bigg)\bigg)
=-\del\bet(t)x^i(t-1)\exp\bigg(\int_0^t\alp(s)\:ds\bigg)\ge C
$$
for almost every $t\in[t_0-\eps,t_0+\eps]$,
and that $x^i(t)$ has exactly one zero in $[t_0-\eps,t_0+\eps]$,
which is located in the interior of this interval, and at which $x^i(t)$ changes sign.

Now let $\tau\in\IR$ be such that $x^0(\tau)\ne 0$ and $x^0(\tau-1)\ne 0$. (Certainly, the set of such $\tau$ is
dense.) Then it follows from the above observations and the compact-uniform convergence $x^i(t)\to x^0(t)$ that
$$
\sch(x^i_\tau)=\sch(x^0_\tau)\qquad\hbox{hence}\qquad V^\pm(x^i_\tau)=V^\pm(x^0_\tau)=J,
$$
for all large $i$. To complete the proof of the lemma, let $I\subseteq\IR$ be any bounded
interval and fix $\tau_1,\tau_2\in\IR$ such that $\tau_1\le t\le \tau_2$ for every $t\in I$,
and where $x^0(\tau_j)\ne 0$ and $x^0(\tau_j-1)\ne 0$ for $j=1,2$. Then
$V^\pm(x^i_{\tau_1})=V^\pm(x^i_{\tau_2})=J$ for all large $i$, and thus for such $i$ we have that
$V^\pm(x^i_t)=J$ for every $t\in[\tau_1,\tau_2]$ by the monotonicity property~(1) of Proposition~5.5.
In particular, this holds for all $t\in I$, as desired.~\fp
\vv

\noindent {\bf Proof of Theorem~5.1.}
Consider the homotoped equation~\eqref{hom}, with $\kappa\in[0,1]$, and
let $U_\kap(t,\tau)$ denote the associated linear process on $X=C([-1,0])$.
Let $M_\kap=U_\kap(\gam,0)$ denote the associated monodromy operator as before, and
denote the Floquet multipliers, that is, the nonzero spectra of $M_\kap$,
by $\{\lam_{\kap,k}\}_{k=1}^\infty$, ordered so that~\eqref{mults} holds.
(No {\it a priori} assumption of continuity or any regularity in $\lam_{\kap,k}$ as
$\kap$ varies is required or assumed.)
By Corollary~5.4 we have $\lam_{\kap,k}\ne 0$ for every $\kap\in[0,1]$ and every $k\ge 1$.

Now let $k_0$ be a positive integer for which $k_0-m$ is even.
Motivated by the remarks following the statement of Proposition~5.6,
we shall say that a value $\kap\in[0,1]$ is {\boldmath $k_0$}{\bf -regular} if
$$
\calj_\kap(|\lam_{\kap,k_0}|)=k_0-1\qquad\hbox{and}\qquad
\left\{\begin{array}{lcr}
\lambda_{\kap,k_0-1}\lambda_{\kap,k_0}>0 & \hbox{ if} &\bb k_0\ge 2,\\
\\
\lambda_{\kap,1}>0 & \hbox{ if} &\bb k_0=1.
\end{array}\right.
$$
Here $\calj_\kap(\rho)$ is as in Proposition~5.6, but for the homotoped equation~\eqref{hom}.
As remarked, in order to prove the theorem it suffices
to show that every $\kap\in[0,1]$ is $k_0$-regular, for every positive $k_0$ for which
$k_0-m$ is even. Let such $k_0$ be fixed, and define the set
$$
S=\{\kap\in[0,1]\:|\:\kap\mbox{ is }k_0\mbox{-regular}\}.
$$
Also as remarked, we have $0\in S$, namely $k_0$-regularity for the constant coefficient equation.
Thus to prove that $S=[0,1]$, it is enough to show that $S$ is both open and closed in the relative
topology of $[0,1]$.

To prove closedness, take a sequence $\kap_i\in S$ with a limit $\kap_i\to\kap_0$.
Without loss we may assume that
\be{ilim}
\lim_{i\to\infty}\lambda_{\kap_i,k}=\lambda_{\kap_0,k}
\ee{ilim}
for every $k$ by taking a subsequence of $i\to\infty$ and possibly reordering the
indexing of the limits (but still maintaining~\eqref{mults}).
Denoting the subspace $\calg_{\kap,\rho}$ as in~\eqref{calg} for the homotoped equation~\eqref{hom}, consider
for $i\ne 0$ the subspace
$$
{E}_i=\calg_{\kap_i,\rho_i},\qquad \rho_i=|\lam_{\kap_i,k_0}|.
$$
Without loss all the ${E}_i$
have the same dimension, either $1$ or $2$, by again passing to a subsequence.
Then $M_{\kap_i}{E}_i={E}_i$ and the eigenvalues of $M_{\kap_i}$ restricted to ${E}_i$
are $\lam_{\kap_i,k_0}$ and possibly also $\lam_{\kap_i,k_0\pm 1}$. Then in light of the
limits~\eqref{ilim} and the fact that $\lam_{\kap_0,k}\ne 0$,
Lemma~5.9 applies. Thus take any $\vph^i\in{E}_i$ with $\|\vph^i\|=1$, and let $x^i(t)$ denote the
corresponding solution of~\eqref{hom} for $t\in\IR$. Therefore, by passing to a subsequence,
we have that $x^i(t)\to x^0(t)$ uniformly on compact intervals, where $x^0(t)$ satisfies~\eqref{hom} with $\kap=\kap_0$.
Thus $\vph^i\to\vph^0=x^0_0$ where $\vph^0\in\calg_{\kap_0,\rho_0}$, where
$$
\rho_0=\lim_{i\to\infty}\rho_i=\lim_{i\to\infty}|\lam_{\kap_i,k_0}|=|\lam_{\kap_0,k_0}|.
$$
Therefore
$$
V^\pm(x^i_t)=\calj_{\kap_i}(|\lam_{\kap_i,k_0}|)=k_0-1,\qquad
V^\pm(x^0_t)=\calj_{\kap_0}(|\lam_{\kap_0,k_0}|),
$$
for all $t\in\IR$, and
it follows by Lemma~5.10 that $\calj_{\kap_0}(|\lam_{\kap_0,k_0}|)=k_0-1$, as desired.
Finally,
\be{posl}
\lam_{\kap_0,k_0-1}\lam_{\kap_0,k_0}=\lim_{i\to\infty}\lam_{\kap_i,k_0-1}\lam_{\kap_i,k_0}\ge 0,
\ee{posl}
and as $\lam_{\kap_0,k}\ne0$ for every $k$, the above limit must be strictly positive
(with the obvious modification for $k_0=1$).
This proves that $\kap_0$ is $k_0$-regular, that is, $\kap_0\in S$. Thus $S$ is closed.

We now prove that $S$ is (relatively) open in $[0,1]$. Let
$\kap_0\in S$ and consider any sequence $\kap_i\in[0,1]$ with $\kap_i\to\kap_0$. We must show that
$\kap_i\in S$ for all large $i$.
Assume to the contrary that $\kap_i\not\in S$ for infinitely many $i$. In fact, without loss,
by taking a subsequence we may assume that $\kap_i\not\in S$ for all large $i$. We also again
may assume that~\eqref{ilim} holds. Now denote $J_i=\calj_{\kap_i}(|\lambda_{\kap_i,k_0}|)$.
Arguing much as in the above paragraph where closedness is proved, we obtain a sequence
of solutions $x^i(t)\to x^0(t)$, although now with
$$
V^\pm(x^i_t)=\calj_{\kap_i}(|\lam_{\kap_i,k_0}|)=J_i,\qquad
V^\pm(x^0_t)=\calj_{\kap_0}(|\lam_{\kap_0,k_0}|)=k_0-1,
$$
for all $t\in\IR$. Lemma~5.10 again applies, and implies that $J_i=k_0-1$ for all large $i$.
Finally, we note again the limit in~\eqref{posl} holds, with strict
positivity of the limit, hence $\lam_{\kap_i,k_0-1}\lam_{\kap_i,k_0}>0$ for all large $i$.
Thus $\kap_i\in S$ for all large $i$, a contradiction. This proves the openness of $S$ and completes the proof.~\fp
\vv

\noindent {\bf Proof of Proposition~5.2.}
For any choice of solutions $x^j(t)$ as in the statement of the proposition, let $\XXI_t\in X^{\wedge m}$
denote the determinant as in~\eqref{xi}, and thus $\XXI_t=\WW(t,\tau)\XXI_\tau$ for $t\ge\tau$, as in~\eqref{xi2}.
Further, $x^i(t+\gamma)$ for each $i$ is a linear combination of the solutions $x^j(t)$ for $1\le j\le m$, say
$$
x^i(t+\gamma)=\sum_{j=1}^mq_{i,j}x^j(t)
$$
for $t\in\IR$. The eigenvalues of the matrix $Q=\{q_{i,j}\}$ are simply the eigenvalues $\lambda_k$, for $1\le k\le m$,
of the monodromy operatory $M$ of equation~\eqref{0} restricted to $\calh_{|\lambda_m|}$, and so $\det Q=\lambda_0$ where
$\lambda_0=\lambda_1\lambda_2\cdots\lambda_m$. Therefore for every $t\in\IR$
\be{fl}
\XXI_{t+\gamma}=\lambda_0\XXI_t,\qquad\mbox{and thus}\qquad
M^{\wedge m}\XXI_0=\lambda_0\XXI_0.
\ee{fl}
That is, $\XXI_0$ is an eigenvector of $M^{\wedge m}$ for the eigenvalue $\lambda_0$.

On the other hand, by Theorem~5.1 and Proposition~2.4, and the remarks following Proposition~2.4, the quantity $\lambda_0$
is a positive eigenvalue of $M^{\wedge m}$ of simple algebraic multiplicity, and with $|\lambda|<\lambda_0$
for all other spectral points of $M^{\wedge m}$. As noted in the Introduction, by a generalization~\cite{31} of
the Kre\u\i n-Rutman Theorem the operator $M^{\wedge m}$ possesses an eigenvector in the positive cone $K_m$
for the eigenvalue $\lambda_0$. Thus by the simplicity $\lambda_0$, this eigenvector is a multiple of the
eigenvector $\XXI_0$, and this gives~\eqref{gl} for some $C$, at least for $t=0$. Positivity for $t\ge 0$ then
holds because of the positivity of the operator $\WW(t,0)$, and positivity for all $t\in\IR$ holds
because of the Floquet property~\eqref{fl}.~\fp

\section{The Proof of Proposition~4.2}

Here we prove Proposition~4.2 for general $m$. Recall that this result was used in
the proofs of Theorem~4.1 (the Positivity Theorem), and as well required for
Theorem~5.1, Propositions~5.2 and~5.3, and Corollary~5.4.
Our approach to proving Proposition~4.2 follows that of the special case with $m=2$ and $\eta=1$
outlined in Section~4. Namely, with
$\eta\le 1$ in Proposition~4.2, we have an explicit
expression~\eqref{11} for each $\wt U_k(\tau+\eta,\tau)$, and composing these expressions will provide an
explicit formula for $\wt\UU(\tau+\eta,\tau)\vph$, for any $\vph\in X^{\otimes m}$.
Then, assuming that $\vph\in X^{\wedge m}$, namely that $\vph$ is anti-symmetric, we will observe significant
cancellations in this formula, and this will yield a much simpler formula for $\wwtilde(\tau+\eta,\tau)\vph$.
\vv

\noindent {\bf Proof of Proposition~4.2.}
Fix $m$, $\tau$, $\eta$, and $b(t)$, as in the statement of the proposition.
Recall the notation~\eqref{brho} and define operators
$Z_k,B_k\in\call(X^{\otimes m})$ for $1\le k\le m$ by
$$
\begin{array}{lcl}
(Z_k\psi)(\theta) &\bb = &\bb
\psi(\theta_1,\ldots,\theta_{k-1},\ome(\theta_k),\theta_{k+1},\ldots,\theta_m),\\
\\
(B_k\psi)(\theta) &\bb = &\bb
\ds{\int_{-1}^{\Ome(\theta_k)}b^{\tau+1}(s)\psi(\theta_1,\ldots,\theta_{k-1},s,\theta_{k+1},\ldots,\theta_m)\:ds,}
\end{array}
$$
where we denote
$$
\ome(s)=\min\{\eta+s,0\},\qquad \Ome(s)=\max\{\eta+s,0\}-1.
$$
Then from~\eqref{11} one sees that
$$
\wt U_k(\tau+\eta,\tau)=Z_k-B_k,
$$
and so
\be{bigu2}
\wt\UU(\tau+\eta,\tau)=(Z_1-B_1)(Z_2-B_2)\cdots(Z_m-B_m)
\ee{bigu2}
by equation~\eqref{bigu}. Next
observe the commutativity properties
\be{com}
Z_jZ_k=Z_kZ_j,\qquad B_jB_k=B_kB_j,\qquad Z_jB_k=B_kZ_j,\qquad \hbox{provided }j\ne k.
\ee{com}
Concerning symmetries, let us define the swap operators $S_{j,k}\in\call(X^{\otimes m})$ for $1\le j,k\le m$ by
$$
(S_{j,k}\psi)(\theta)=\psi(\overline\theta),\qquad
\overline\theta_i=\left\{\begin{array}{lcl}
\theta_k, & \hbox{ if} &\bb i=j,\\
\\
\theta_j, & \hbox{ if} &\bb i=k,\\
\\
\theta_i, & \hbox{ if} &\bb i\ne j\hbox{ and }i\ne k.
\end{array}\right.
$$
Note that $S_{j,k}=S_{\sig_{j,k}}$ in the notation~\eqref{cals},
where $\sig_{j,k}\in\cals_m$ is the permutation satisfying
$\sig_{j,k}(j)=k$ and $\sig_{j,k}(k)=j$, with $\sig_{j,k}(i)=i$ if $i\ne j$ and $i\ne k$.
We write $S_{j,k}$ rather than $S_{\sig_{j,k}}$ for simplicity of notation.
Each $S_{j,k}$ is an isometry on the space $X^{\otimes m}$, and
of course $S_{j,k}=S_{k,j}=S_{j,k}^{-1}$. One easily checks that
\be{inter}
\begin{array}{lcl}
S_{j,k}Z_k=Z_jS_{j,k}, &\bb \qquad &\bb S_{j,k}B_k=B_jS_{j,k},\\
\\
S_{j,k}Z_i=Z_iS_{j,k}, &\bb \qquad &\bb S_{j,k}B_i=B_iS_{j,k},\\
\\
\makebox[0pt][l]{in the last two cases provided $i\ne j$ and $i\ne k$.}
\end{array}
\ee{inter}
Let us note that $S_{j,k}\wt\UU(\tau+\eta,\tau)=\wt\UU(\tau+\eta,\tau)S_{j,k}$ by~\eqref{com} and~\eqref{inter},
and so the space $X^{\wedge m}$ of anti-symmetric functions is invariant under
the operator $\wt\UU(\tau+\eta,\tau)$. However, also note that $X^{\wedge m}$ is not in general invariant under
either of the operators $Z_k$ or $B_k$.

Fix $\vph\in X^{\wedge m}$ and let the right-hand side of~\eqref{bigu2} be expanded and act on $\vph$. We obtain
a sum of all terms of the form
\be{cprod}
\pm C_1C_2\cdots C_m\vph,\qquad C_k=Z_k\hbox{ or }B_k,
\ee{cprod}
where the sign $\pm$ is $(-1)^i$, where $i$ is the number of $B_k$ appearing in the product~\eqref{cprod}.
Now fix $\theta\in T_m$, along with the integer $a$ satisfying~\eqref{thet} as in the statement of the proposition.
Both $\theta$ and $a$ will stay fixed for the remainder of this proof.
Here we make two crucial observations. First,
suppose that $C_k=B_k$ for some $k$ satisfying $1\le k\le a$. Note that $\theta_k\le\theta_a\le-\eta$ so $\Ome(\theta_k)=-1$,
and thus $(B_k\psi)(\theta)=0$ for every $\psi\in X^{\otimes m}$, for our chosen $\theta$. (We are not claiming
that $B_k\psi$ is the zero function, however, but simply that it vanishes
at this particular $\theta$.) In particular, taking
$\psi=C_1\cdots C_{k-1}\wh C_k C_{k+1}\cdots C_m\vph$ where the hat~$\wh{{\ }{\ }}$ indicates
the term is omitted, and using the commutativity properties~\eqref{com}, we have that
$C_1C_2\cdots C_m\vph=B_k\psi$ and thus
\be{czer}
(C_1C_2\cdots C_m\vph)(\theta)=0.
\ee{czer}
Secondly, suppose that $C_{k_1}=Z_{k_1}$ and $C_{k_2}=Z_{k_2}$ for two values $k_1,k_2$
satisfying $a<k_1<k_2\le m$. Then $\ome(\theta_{k_1})=\ome(\theta_{k_2})=0$ and so
$$
(Z_{k_1}Z_{k_2}\psi)(\theta)=
\psi(\theta_1,\ldots,\theta_{k_1-1},0,\theta_{k_1+1},\ldots,\theta_{k_2-1},0,\theta_{k_2+1},\ldots,\theta_m)
$$
for every $\psi\in X^{\otimes m}$. If it is further the case that $\psi$ satisfies
\be{k1k2}
S_{k_1,k_2}\psi=-\psi
\ee{k1k2}
identically as functions, then in fact $(Z_{k_1}Z_{k_2}\psi)(\theta)=0$.
(Again, this is for our chosen $\theta$, and there is no claim that $Z_{k_1}Z_{k_2}\psi$ is the zero function.)
Taking
$$
\psi=C_1\cdots C_{k_1-1}\wh C_{k_1} C_{k_1+1}\cdots C_{k_2-1}\wh C_{k_2} C_{k_2+1}\cdots C_m\vph
$$
and recalling that $\vph\in X^{\wedge m}$,
we see from the anti-symmetry of $\vph$ and the properties~\eqref{inter} that~\eqref{k1k2} holds.
Thus $C_1C_2\cdots C_m\vph=Z_{k_1}Z_{k_2}\psi$ and again~\eqref{czer} holds.

Following the above two crucial observation, we see that only a few select terms survive
in the expansion of ~\eqref{bigu2} when applied to $\vph\in X^{\wedge m}$ and evaluated at our chosen $\theta$.
These are the terms as in~\eqref{cprod} for which the following both hold:
$C_k=Z_k$ for every $k$ satisfying $1\le k\le a$;
and $C_k=Z_k$ for at most one $k$ in the range $a<k\le m$. We conclude that for every $\vph\in X^{\wedge m}$,
we have that
\be{gam}
[\wwtilde(\tau+\eta,\tau)\vph](\theta)=(-1)^{m-a}[(\Gamma_0\vph)(\theta)]
+(-1)^{m-a+1}\sum_{k=1}^{m-a}(\Gamma_k\vph)(\theta),
\ee{gam}
where
\be{gammas}
\begin{array}{lcl}
\Gamma_0 &\bb = &\bb Z_1\cdots Z_aB_{a+1}\cdots B_m,\\
\\
\Gamma_k &\bb = &\bb Z_1\cdots Z_aB_{a+1}\cdots B_{a+k-1}Z_{a+k}B_{a+k+1}\cdots B_m,
\end{array}
\ee{gammas}
for $1\le k\le m-a$. The signs $(-1)^{m-a}$ and $(-1)^{m-a+1}$ occurring in the above formulas count
the number of $B_j$ terms.

(We remark on these formulas in the extreme cases $a=0$ and $a=m$. If $a=0$ then the factors $Z_1\cdots Z_a$
in~\eqref{gammas} are simply absent. If $a=m$ then the summation in~\eqref{gam} is empty, so no $\Gamma_k$
are defined for $k>0$. Also, if $a=m$ then $\Gamma_0=Z_1\cdots Z_m$. We leave the verification of these facts to the reader.)

We note a peculiarity of the formula~\eqref{gam}, namely
that the operators $\Gamma_0$ and $\Gamma_k$ depend on $a$, which in
turn depends on $\theta$. Thus it is not the case, in general, that the operator $\wwtilde(\tau+\eta,\tau)$
is the sum of the operators $\Gamma_0$ and $\Gamma_k$ with the indicated signs as in~\eqref{gam}.
Rather, equation~\eqref{gam}
is only valid pointwise for those $\theta$ which satisfy~\eqref{thet}. We emphasize that
for this reason, we work with a fixed $\theta\in T_m$.

Let us now evaluate the terms in~\eqref{gam}. We first consider the term involving $\Gamma_k$ with $1\le k\le m-a$.
As noted, the case $a=m$ is vacuous. If $a=m-1$, then $k=1$ and $\Gamma_1=Z_1\cdots Z_{m-1}Z_m$, which immediately gives
\be{ma1}
(\Gamma_1\vph)(\theta)=\vph(\eta+\theta_1,\ldots,\eta+\theta_{m-1},0),
\ee{ma1}
as we note that $\ome(\theta_j)=\eta+\theta_j$ for $1\le j\le m-1$ but
$\ome(\theta_m)=0$. Now let us assume that $0\le a\le m-2$.
We have, from the above formula~\eqref{gammas} for $\Gamma_k$, and also from the formulas for $Z_k$ and $B_k$, that
\be{gk}
\begin{array}{lcl}
(\Gamma_k\vph)(\theta) &\bb = &\bb
\ds{\int_{-1}^{\Ome(\theta_{a+1})}\!\!\!\!\cdots\wh{\int_{-1}^{{\Ome(\theta_{a+k})}}}\!\!\!\!\cdots\int_{-1}^{\Ome(\theta_{m})}
b^{\tau+1}(s_{1})\cdots \wh{b^{\tau+1}(s_k)}\cdots b^{\tau+1}(s_{m-a})}\\
\\
&\bb &\bb
\times\:\:\vph(\eta+\theta_1,\ldots,\eta+\theta_a, s_{1},\ldots,s_{k-1},0,s_{k+1},\ldots,s_{m-a})\:ds_{m-a}\cdots \wh{ds_{k}}\cdots ds_{1},
\end{array}
\ee{gk}
where here again (and below) the hat~$\wh{{\ }{\ }}$ denotes that the indicated expression is omitted,
and where we note that $\ome(\theta_j)=\eta+\theta_j$ for $1\le j\le a$ but $\ome(\theta_{a+k})=0$.
(If $a=0$, then the terms $\eta+\theta_j$ in the integral~\eqref{gk} are simply absent.)
Next, by permuting the arguments of $\vph$ in~\eqref{gk} and using
the anti-symmetry of $\vph$, we see that
\be{perm}
\begin{array}{lcl}
\makebox[0pt][l]{$\vph(\eta+\theta_1,\ldots,\eta+\theta_a, s_{1},\ldots,s_{k-1},0,s_{k+1},\ldots,s_{m-a})$}\\
\\
&\bb = &\bb
(-1)^{a+k+(a+1)m}\vph(t_{1},\ldots,t_{k-1},t_{k},\ldots,t_{m-a-1},\eta+\theta_1,\ldots,\eta+\theta_a,0),\\
\\
\hbox{where } t_j &\bb = &\bb
\left\{\begin{array}{lcl}
s_j, & \hbox{ for} &\bb 1\le j\le k-1,\\
\\
s_{j+1}, & \hbox{ for} &\bb k\le j\le m-a-1.
\end{array}\right.
\end{array}
\ee{perm}
The explanation for the term $(-1)^{a+k+(a+1)m}$, arising from the anti-symmetry of $\vph$, is as follows.
Each term $s_j$, for $1\le j\le k-1$,
is moved leftward $a$ places by means of swaps with adjacent terms $\eta+\theta_i$
for $1\le i\le a$. This is a total of $a(k-1)$
swaps of such terms. Each term $s_j$, for $k+1\le j\le m-a$, is moved leftward $a+1$ places by means
of swaps with adjacent terms $0$ and then $\eta+\theta_i$. This is an additional $(a+1)(m-a-k)$ swaps.
The total number of swaps is thus $a(k-1)+(a+1)(m-a-k)$, and one sees
easily that this number has the same parity as $a+k+(a+1)m$.

Let us now introduce
some notation which will simplify our calculations. We define
a bounded linear operator $E\in\call(X^{\wedge m},X^{\wedge (m-a-1)})$ by
\be{q}
(E\psi)(t_1,\ldots,t_{m-a-1})=\psi(t_{1},\ldots,t_{m-a-1},\eta+\theta_1,\ldots,\eta+\theta_a,0).
\ee{q}
We also define operators
$\iop{i}{j},\jop{i}{j}\in\call(X^{\otimes (m-a-1)})$ by
\be{iop}
\begin{array}{lcr}
(\iop{i}{j}\psi)(t_1,\ldots,t_{m-a-1}) &\bb =&\bb
\ds{\int_{\Ome(\theta_{a+i-1})}^{\Ome(\theta_{a+i})}b^{\tau+1}(s)\psi(t_{1},\ldots,t_{j-1},s,t_{j+1},\ldots,t_{m-a-1})\:ds,}\\
\\
(\jop{i}{j}\psi)(t_1,\ldots,t_{m-a-1}) &\bb =&\bb
\ds{\int_{-1}^{\Ome(\theta_{a+i})}b^{\tau+1}(s)\psi(t_{1},\ldots,t_{j-1},s,t_{j+1},\ldots,t_{m-a-1})\:ds,}
\end{array}
\ee{iop}
where $1\le i\le m-a$ and $1\le j\le m-a-1$, and where we note that
$\Ome(\theta_{a})=-1$. (If $a=0$, then by convention we set $\Ome(\theta_0)=-1$.)
Observe that
$$
\jop{i}{j}=\iop{1}{j}+\iop{2}{j}+\cdots+\iop{i}{j}
$$
holds. With this notation,
and upon inserting the formula~\eqref{perm} into~\eqref{gk}, one sees that~\eqref{gk} takes the form
\be{jform}
(\Gamma_k\vph)(\theta)=(-1)^{a+k+(a+1)m}\jop{1}{1}\jop{2}{2}\cdots \jop{k-1}{k-1}\jop{k+1}{k}\cdots \jop{m-a}{m-a-1}E\vph.
\ee{jform}
We remind the reader again, that $\theta$ has been fixed and does not serve as the
argument of the functions $E\psi$, $\iop{i}{j}\psi$, and $\jop{i}{j}\psi$ above. Rather, these are functions
of the variables $(t_1,\ldots,t_{m-a-1})$. The functions $\iop{i}{j}\psi$ and $\jop{i}{j}\psi$, in particular,
are constant in the variable $t_j$, as the right-hand sides in~\eqref{iop}
are independent of $t_j$. Thus
the right-hand side of~\eqref{jform} is formally a function
of $(t_1,\ldots,t_{m-a-1})$, and in fact is a constant function
of those variables. That constant value is the value of the function $\Gamma_k\vph$ evaluated at the point~$\theta$.

The operators $\iop{i}{j}$ and $\jop{i}{j}$ act on the full tensor product, and not just the exterior product. That
is, no symmetry assumption is made on the argument function $\psi\in X^{\otimes (m-a-1)}$ in~\eqref{iop}.
Observe that
\be{asym}
\begin{array}{l}
S_{j_1,j_2}\iop{i}{j}=\iop{i}{j}S_{j_1,j_2},\qquad
S_{j_1,j_2}\jop{i}{j}=\jop{i}{j}S_{j_1,j_2},\\
\\
\hbox{in both cases provided }j\ne j_1\hbox{ and }j\ne j_2.
\end{array}
\ee{asym}
Thus if it is the case that $\psi$ is anti-symmetric in $t_{j_1}$ and $t_{j_2}$, meaning that
$S_{j_1,j_2}\psi=-\psi$, then
$\iop{i}{j}\psi$ and $\jop{i}{j}\psi$ are also anti-symmetric in these variables as long as $j\ne j_1,j_2$.

It is easily seen that
\be{zer2}
\begin{array}{l}
\jop{i}{j_1}\jop{i}{j_2}\psi=0,\qquad
\jop{i}{j_1}\jop{i+1}{j_2}\psi=\jop{i}{j_1}\iop{i+1}{j_2}\psi,\qquad
\jop{i}{j_1}\jop{i+2}{j_2}\psi=\jop{i}{j_1}(\iop{i+1}{j_2}+\iop{i+2}{j_2})\psi,\\
\\
\hbox{in every case provided }j_1\ne j_2\hbox{ and }S_{j_1,j_2}\psi=-\psi.
\end{array}
\ee{zer2}
Indeed, the first equation in~\eqref{zer2} holds as it is simply the integral over
the square $[-1,\Ome(\theta_{a+i})]^2$ of an anti-symmetric function of $(t_{j_1},t_{j_2})$. The
second and third equations in~\eqref{zer2} follow from the first equation
because $\jop{i+1}{j}=\jop{i}{j}+\iop{i+1}{j}$ and $\jop{i+2}{j}=\jop{i}{j}+\iop{i+1}{j}+\iop{i+2}{j}$.

We now use the identities~\eqref{zer2} to obtain
a simplification of equation~\eqref{jform} when
the function $\vph$ is anti-symmetric.
We begin with the rightmost pair of $J$-operators in~\eqref{jform}, namely $\jop{m-a-1}{m-a-2}\jop{m-a}{m-a-1}$
and move to the left. The result is that each factor $\jop{i+1}{i}$ is replaced with $\iop{i+1}{i}$, and each $\jop{i}{i}$
is replaced with $\iop{i}{i}$, except for the factor $\jop{k+1}{k}$ which is replaced with $\iop{k}{k}+\iop{k+1}{k}$.
At each stage we observe, using~\eqref{asym}, that the relevant function is anti-symmetric in the appropriate
variables, as required by~\eqref{zer2}. Finally noting
that $\jop{1}{1}=\iop{1}{1}$, we conclude directly that
\be{iform}
(\Gamma_k\vph)(\theta)=(-1)^{a+k+(a+1)m}\iop{1}{1}\iop{2}{2}\cdots \iop{k-1}{k-1}(\iop{k}{k}+\iop{k+1}{k})\iop{k+2}{k+1}\cdots \iop{m-a}{m-a-1}E\vph.
\ee{iform}

The reader can verify that the formula~\eqref{iform} degenerates in the extreme cases
of the indices, as follows. If $m-a\ge 4$ then
$$
\begin{array}{rcl}
(\Gamma_1\vph)(\theta) &\bb = &\bb (-1)^{a+1+(a+1)m}(\iop{1}{1}+\iop{2}{1})\iop{3}{2}\cdots \iop{m-a}{m-a-1}E\vph,\\
\\
(\Gamma_2\vph)(\theta) &\bb = &\bb (-1)^{a+2+(a+1)m}\iop{1}{1}(\iop{2}{2}+\iop{3}{2})\iop{4}{3}\cdots \iop{m-a}{m-a-1}E\vph,\\
\\
(\Gamma_{m-a-1}\vph)(\theta) &\bb = &\bb (-1)^{m-1+(a+1)m}\iop{1}{1}\iop{2}{2}\cdots \iop{m-a-2}{m-a-2}(\iop{m-a-1}{m-a-1}+\iop{m-a}{m-a-1})E\vph,\\
\\
(\Gamma_{m-a}\vph)(\theta) &\bb = &\bb (-1)^{m+(a+1)m}\iop{1}{1}\iop{2}{2}\cdots \iop{m-a-1}{m-a-1}E\vph.
\end{array}
$$
If $m-a=3$ then we have
$$
\begin{array}{lcl}
(\Gamma_1\vph)(\theta) &\bb = &\bb (-1)^{a+1+(a+1)m}(\iop{1}{1}+\iop{2}{1})\iop{3}{2}E\vph,\\
\\
(\Gamma_2\vph)(\theta) &\bb = &\bb (-1)^{a+2+(a+1)m}\iop{1}{1}(\iop{2}{2}+\iop{3}{2})E\vph,\\
\\
(\Gamma_{3}\vph)(\theta) &\bb = &\bb (-1)^{a+3+(a+1)m}\iop{1}{1}\iop{2}{2}E\vph,
\end{array}
$$
while if $m-a=2$ then we have
$$
\begin{array}{lcl}
(\Gamma_1\vph)(\theta) &\bb = &\bb (-1)^{a+1+(a+1)m}(\iop{1}{1}+\iop{2}{1})E\vph,\\
\\
(\Gamma_2\vph)(\theta) &\bb = &\bb (-1)^{a+2+(a+1)m}\iop{1}{1}E\vph.
\end{array}
$$
(If $m-a=1$ we have simply the formula~\eqref{ma1}.)
These formulas follow easily from~\eqref{jform} using the identities~\eqref{zer2}.
Now define the operators
$$
\begin{array}{l}
R_1 = \iop{2}{1}\iop{3}{2}\cdots \iop{m-a}{m-a-1},\\
\\
R_k = \iop{1}{1}\cdots \iop{k-1}{k-1} \iop{k+1}{k}\cdots \iop{m-a}{m-a-1},\qquad 2\le k\le m-a-1,\\
\\
R_{m-a} = \iop{1}{1} \iop{2}{2}\cdots \iop{m-a-1}{m-a-1},\\
\\
R_{m-a+1} = 0.
\end{array}
$$
Then~\eqref{iform} can be rewritten as
$$
(\Gamma_k\vph)(\theta)=(-1)^{a+k+(a+1)m}(R_{k+1}+R_{k})E\vph,
$$
and we see that the formula is valid for all values $1\le k\le m-a$ with $m-a\ge 2$,
including the extreme cases above.
It follows immediately that the summation in~\eqref{gam} telescopes to give
$$
\hspace*{-.1em}
\begin{array}{l}
\ds{(-1)^{a+1+(a+1)m}\sum_{k=1}^{m-a}(\Gamma_k\vph)(\theta)=R_1E\vph}\\
\\
\qquad\qquad
\ds{=\int_{\Ome(\theta_{a+1})}^{\Ome(\theta_{a+2})}\!\!\!\!\cdots\int_{\Ome(\theta_{m-1})}^{\Ome(\theta_{m})}
b^{\tau+1}(t_{1})\cdots b^{\tau+1}(t_{m-a-1})
[(E\vph)(t_{1},\ldots,t_{m-a-1})]\:dt_{m-a-1}\cdots \makebox[0pt][l]{$dt_1.$}}
\end{array}
$$
Upon multiplying the above formula by $(-1)^{am}$ and noting that
$(-1)^{am}(-1)^{a+1+(a+1)m}=(-1)^{m-a+1}$, we obtain
\be{x1}
\begin{array}{lcl}
\ds{(-1)^{m-a+1}\sum_{k=1}^{m-a}(\Gamma_k\vph)(\theta)} &\bb = &\bb
\ds{(-1)^{am}\int_{\eta+\theta_{a+1}-1}^{\eta+\theta_{a+2}-1}\!\!\!\!\cdots\int_{\eta+\theta_{m-1}-1}^{\eta+\theta_{m}-1}
b^{\tau+1}(t_{1})\cdots b^{\tau+1}(t_{m-a-1})}\\
\\
&\bb &\bb
\qquad\qquad
\times\:\:\vph(t_{1},\ldots,t_{m-a-1},\eta+\theta_1,\ldots,\eta+\theta_a,0)\:dt_{m-a-1}\cdots \makebox[0pt][l]{$dt_1$,}
\end{array}
\ee{x1}
where the formula~\eqref{q} for $E$ is used along with the fact that
$\Ome(\theta_j)=\eta+\theta_j-1$ for $a+1\le j\le m$.

The calculation of $(\Gamma_0\vph)(\theta)$ is handled in a similar fashion,
and in fact is slightly simpler. If $a=m$ then $\Gamma_0=Z_1Z_2\cdots Z_m$, and so
\be{gzz}
(\Gamma_0\vph)(\theta)=\vph(\eta+\theta_1,\ldots,\eta+\theta_{m}).
\ee{gzz}
Let us therefore assume that $0\le a\le m-1$.
We have first that
\be{gz}
\begin{array}{lcl}
(\Gamma_0\vph)(\theta) &\bb = &\bb
\ds{\int_{-1}^{\Ome(\theta_{a+1})}\!\!\!\!\cdots\int_{-1}^{\Ome(\theta_{m})}
b^{\tau+1}(s_{1})\cdots b^{\tau+1}(s_{m-a})}\\
\\
&\bb &\bb
\qquad\qquad\times\:\:\vph(\eta+\theta_1,\ldots,\eta+\theta_a,s_{1},\ldots,s_{m-a})\:ds_{m-a}\cdots ds_{1},
\end{array}
\ee{gz}
and also, using the anti-symmetry of $\vph$, that
$$
\vph(\eta+\theta_1,\ldots,\eta+\theta_a, s_{1},\ldots,s_{m-a})
=(-1)^{a(m-a)}\vph(s_{1},\ldots,s_{m-a},\eta+\theta_1,\ldots,\eta+\theta_a).
$$
Introducing the operator $E_0\in\call(X^{\wedge m},X^{\wedge (m-a)})$ given by
\be{q0}
(E_0\psi)(t_1,t_2,\ldots,t_{m-a})=\psi(t_{1},\ldots,t_{m-a},\eta+\theta_1,\ldots,\eta+\theta_a),
\ee{q0}
we have from~\eqref{gz} that
$$
(\Gamma_0\vph)(\theta)=(-1)^{a(m-a)}\jop{1}{1}\jop{2}{2}\cdots \jop{m-a}{m-a}E_0\vph.
$$
Here the operators $\jop{i}{j}$, and $\iop{i}{j}$ are as before, except they now operate on
functions of $m-a$ variables rather than $m-a-1$ variables.
Using~\eqref{asym} and~\eqref{zer2} as before, we obtain
$$
\begin{array}{l}
(-1)^{a(m-a)}[(\Gamma_0\vph)(\theta)]=\iop{1}{1}\iop{2}{2}\cdots \iop{m-a}{m-a}E_0\vph\\
\\
\qquad\qquad
\ds{=\int_{\Ome(\theta_{a})}^{\Ome(\theta_{a+1})}\!\!\!\!\cdots\int_{\Ome(\theta_{m-1})}^{\Ome(\theta_{m})}
b^{\tau+1}(t_{1})\cdots b^{\tau+1}(t_{m-a})
[(E_0\vph)(t_{1},\ldots,t_{m-a})]\:dt_{m-a}\cdots dt_{1}.}
\end{array}
$$
Upon multiplying the above formula by $(-1)^{(a+1)m}$ and noting that $(-1)^{(a+1)m}(-1)^{a(m-a)}=(-1)^{m-a}$,
we obtain
\be{x2}
\begin{array}{lcl}
(-1)^{m-a}[(\Gamma_0\vph)(\theta)] &\bb = &\bb
\ds{(-1)^{(a+1)m}\int_{-1}^{\eta+\theta_{a+1}-1}\int_{\eta+\theta_{a+1}-1}^{\eta+\theta_{a+2}-1}
\!\!\!\!\cdots\int_{\eta+\theta_{m-1}-1}^{\eta+\theta_{m}-1}
b^{\tau+1}(t_{1})\cdots b^{\tau+1}(t_{m-a})}\\
\\
&\bb &\bb
\qquad\qquad\times\:\:\vph(t_{1},\ldots,t_{m-a},\eta+\theta_1,\ldots,\eta+\theta_a)\:dt_{m-a}\cdots dt_{1},
\end{array}
\ee{x2}
where the formula~\eqref{q0} for $E_0$ is used, and where we have that $\Ome(\theta_a)=-1$. Adding the two
equations~\eqref{x1} and~\eqref{x2} and using~\eqref{gam} gives~\eqref{uform}, as desired, at least in
the case that $0\le a\le m-2$. If $a=m-1$ then~\eqref{x1} must be replaced by~\eqref{ma1} to give
the desired formula~\eqref{uformx}.
Finally, if $a=m$ then the term corresponding to~\eqref{x1} is absent, while~\eqref{x2} is given by~\eqref{gzz}
to give~\eqref{uformz}, again as desired. With this, the proposition is proved.~\fp

\section{{\boldmath $u_0$}-Positivity}

Here we consider the question of $u_0$-positivity of the linear process $\WW(t,\tau)\in\call(X^{\wedge m})$.
We assume a uniform sign condition
\be{usc}
0<\bet_1\le(-1)^m\beta(t)\le \bet_2
\ee{usc}
on bounded intervals, which is slightly stronger than the assumption $(-1)^m\beta(t)\ge 0$ of Theorem~4.1.
We maintain the same notation as in the previous sections, with $X=C([-1,0])$
and the cone $K_m\subseteq X^{\wedge m}$ given by~\eqref{km},
with the set $T_m\subseteq[-1,0]$ given by~\eqref{tri}, and where $\WW(t,\tau)\in\call(X^{\wedge  m})$
is the $m$-fold exterior product of the linear process associated to the delay-differential equation~\eqref{0}.

Generally, if $Y$ is a Banach space and $K\subseteq Y$ is a cone,
then we say two elements $u,v\in K\setminus\{0\}$ are {\bf comparable} in case
there exist quantities $C_2\ge C_1>0$ such that
$$
C_1v\le u\le C_2v,
$$
where $\le$ denotes the ordering with respect to $K$. We denote this relation by
$$
u\sim v.
$$
Clearly, $\sim$ is an equivalence relation on $K$.
The following definition is classical.
\vv

\noindent {\bf Definition.}
Suppose that $A\in\call(Y)$ is a positive operator with respect
to a cone $K\subseteq Y$ in a Banach space $Y$,
that is, $Au\in K$ whenever $u\in K$. Let $u_0\in K\setminus\{0\}$.
Then we say that the operator $A$ is {\boldmath $u_0$}{\bf -positive} in case
there exists an integer $k_0\ge 1$ such that $A^{k_0}u\sim u_0$
for every $u\in K\setminus\{0\}$.
(We note that this condition implies further that $A^ku\sim u_0$ for every $k\ge k_0$.)
\vv

In the case of a linear process, we make the following related definition which accounts
for continuous rather than discrete time. Here the constants $C_1$ and $C_2$
are uniform with respect to compact time-intervals.
\vv

\noindent {\bf Definition.}
Suppose that $U(t,\tau)\in\call(Y)$, for $t\ge \tau$, is a linear process
on a Banach space $Y$. Suppose also that $U(t,\tau)$ is a positive
operator with respect to a cone $K\subseteq Y$, for every $t,\tau\in\IR$ with
$t\ge \tau$. Let $u_0\in K\setminus\{0\}$. Then we say that the process $U(t,\tau)$
is {\boldmath $u_0$}{\bf -positive} in case there exists $\eta_0>0$ such
that the following holds. Given any $u\in K\setminus\{0\}$, and given $\tau_1\le\tau_2$ and $\eta_*\ge\eta_0$,
then there exist $C_1>0$ and $C_2>0$
such that
$$
C_1u_0\le U(\tau+\eta,\tau)u\le C_2u_0
$$
for every $\tau\in[\tau_1,\tau_2]$ and $\eta\in[\eta_0,\eta_*]$.
\vv

As was noted in the Introduction, $u_0$-positivity of an operator $A$
implies additional and precise information about the
spectrum and the eigenvectors of $A$, at least when the cone $K$ is reproducing. In particular, if $Au=\lambda u$ with
$u\in K\setminus\{0\}$ and $\lambda>0$, then it is necessarily the case that $u\sim u_0$. In the context of Floquet
theory for delay-differential equations, this refines Proposition~5.2 to give sharper bounds, in
Proposition~7.2 below, on the primary Floquet eigenfunction of the linear process $\WW(t,\tau)$.
\vv

For each $m\ge 2$ let us define a function
\be{u0}
u_m(\theta)=\bigg(\prodform{\xxa}{1\le i<j\le m}{(\theta_j-\theta_i)}\bigg)
\bigg(\prodform{\xxb}{1\le i<j\le m-1}{(1+\theta_i-\theta_j)}\bigg)\qquad
\hbox{for }\theta\in T_m,
\ee{u0}
and extend $u_m$ to all of $[-1,0]^m$ as an anti-symmetric function,
and so $(S_\sig u_m)(\theta)=\sgn(\sig)u_m(\theta)$ for every $\sig\in\cals_m$
and every $\theta\in[-1,0]^m$.
Note that for $m=2$ the range $1\le i<j\le m-1$ of the indices in the second factor of~\eqref{u0} is empty.
In this and other such cases, here and below, we interpret such an empty product to be equal to $+1$ identically.
Also note that the extended function $u_m$ is continuous throughout $[-1,0]^m$, that is $u_m\in X^{\wedge m}$.
This holds because $u_m(\theta)=0$ whenever $\theta\in T_m$ is such that $\theta_j=\theta_{j+1}$ for some $j$
with $1\le j\le m-1$. Of course the polynomial formula~\eqref{u0} is not generally valid
for $\theta\in[-1,0]^m\setminus T_m$.

The following theorem is the main result of this section. It is followed by
a conjecture concerning the natural generalization of this result.
\vv

\noindent {\bf Theorem~7.1.} \bxx
Let $m=2$ or $m=3$ be fixed. Assume that the coefficients
$\alpha,\beta:\IR\to\IR$ in equation~\eqref{0} are locally integrable and that
for every compact interval $[t_1,t_2]\subseteq\IR$ there exist $\bet_1$ and $\bet_2$
such that~\eqref{usc} holds for almost every $t\in[t_1,t_2]$.
Then the linear process $\WW(t,\tau)$ on $X^{\wedge m}$ is $u_m$-positive
with respect to the cone $K_m$, with the above function $u_m$.
Moreover, we have that $\eta_0=3$ if $m=2$ and $\eta_0=5$ if $m=3$ for the
quantity $\eta_0$ in the definition of a $u_0$-positive linear process.
\exx
\vv

\noindent {\bf Conjecture~A\@.}
Let $m\ge 4$. Then Theorem~7.1 holds as stated, but instead for this $m$, and with $\eta_0=2m-1$.
\vv

The following result, which extends Proposition~5.2, is a consequence of Theorem~7.1 and Conjecture~A\@.
One easily sees that the bounds~\eqref{72bnd} there are in fact valid for $t$ in
any compact interval, where $C_1$ and $C_2$ depend on the interval.
Even for the constant coefficient problem~\eqref{cc2}, Conjecture~A and the bounds~\eqref{72bnd}
are unknown for $m\ge 4$.
\vv

\noindent {\bf Proposition~7.2.} \bxx
Let $m\ge 2$ be such that Theorem~7.1 holds as stated, but instead for this $m$.
Consider the $\gamma$-periodic problem in the setting of Proposition~5.2,
but with the additional assumption~\eqref{usc} for some $\bet_1$ and $\bet_2$, for almost every $t\in\IR$.
Then with $x^j(t)$ as in the statement of Proposition~5.2, there exists a constant $C$
and constants $C_2\ge C_1>0$ such that
\be{72bnd}
C_1u_m(\theta)\le C\det\left(\begin{array}{ccc}
x^1(t+\theta_1) &\bb \cdots &\bb x^1(t+\theta_m)\\
\vdots &\bb &\bb \vdots\\
x^m(t+\theta_1) &\bb \cdots &\bb x^m(t+\theta_m)
\end{array}\right)\le C_2u_m(\theta),
\ee{72bnd}
for every $t\in[0,\gam]$ and $\theta=(\theta_1,\theta_2,\ldots,\theta_m)\in T_m$.
\exx
\vv

\noindent {\bf Proof.}
Let $\XXI_t\in X^{\wedge m}$ denote the left-hand side of~\eqref{gl} in Proposition~5.2, with $C$
as indicated. Then this result states that
$\XXI_t\in K_m$ for every $t\in\IR$. Moreover, $\XXI_t=\WW(t,\tau)\XXI_\tau$ holds for every $t\ge\tau$.
In particular, if we take $\tau=-\eta_0$ where $\eta_0$ is as in the definition of a $u_0$-positive process,
there exist $C_2\ge C_1>0$ such that
$$
C_1u_m\le \XXI_t=\WW(t,-\eta_0)\XXI_{-\eta_0}\le C_2u_m
$$
for every $t\in[0,\gam]$. This proves the result.~\fp
\vv

To prove Theorem~7.1 it will be enough to consider the transformed equation~\eqref{1}. Further,
it will be sufficient for our purpose to take $\eta=1$ in the formula~\eqref{uform},
and so we may take $a=0$ in that formula, as per the statement of Proposition~4.2.
Therefore assume that $(-1)^mb(t)\ge 0$ for almost every $t$ satisfying $\tau\le t\le\tau+1$. Then~\eqref{uform} gives
\be{eta1x}
\begin{array}{lcl}
[\wwtilde(\tau+1,\tau)\vph](\theta) &\bb = &\bb
\ds{\int_{\theta_{1}}^{\theta_{2}}\!\!\!\!\cdots\int_{\theta_{m-1}}^{\theta_{m}}
|b^{\tau+1}(t_{1})\cdots b^{\tau+1}(t_{m-1})|
\vph(t_{1},\ldots,t_{m-1},0)\:dt_{m-1}\cdots dt_{1}}\\
\\
&\bb &\bb
\ds{+\int_{-1}^{\theta_{1}}\int_{\theta_{1}}^{\theta_{2}}
\!\!\!\!\cdots\int_{\theta_{m-1}}^{\theta_{m}}
|b^{\tau+1}(t_{0})\cdots b^{\tau+1}(t_{m-1})|\vph(t_{0},\ldots,t_{m-1})\:dt_{m-1}\cdots dt_{0}}
\end{array}
\ee{eta1x}
for every $\vph\in X^{\wedge m}$ provided that $\theta\in T_m$.
(For convenience later, we have reindexed the variables $t_j$ in the final term of~\eqref{eta1x}.)
Note that in the case of odd $m$, where $b(t)\le 0$, the identities
$(-1)^{am}$=$(-1)^{m-a-1}$ and $(-1)^{(a+1)m}$=$(-1)^{m-a}$
are used in taking the absolute values of $b(t)$.

We shall need both positive upper and lower bounds for the operator $\wwtilde(\tau+1,\tau)$,
which is why we assume in Theorem~7.1 that
there are uniform positive upper and lower bounds for $|\beta(t)|$, and therefore for $|b(t)|$, on compact intervals.
The bounds for $\wwtilde(\tau+1,\tau)$
will then be given by appropriate multiples of the operator $\AA=\AA_0+\AA_1$, where
\be{aa}
\begin{array}{lcl}
(\AA_0\vph)(\theta) &\bb = &\bb
\ds{\int_{\theta_{1}}^{\theta_{2}}\!\!\!\!\cdots\int_{\theta_{m-1}}^{\theta_{m}}
\vph(t_{1},\ldots,t_{m-1},0)\:dt_{m-1}\cdots dt_{1},}\\
\\
(\AA_1\vph)(\theta) &\bb = &\bb
\ds{\int_{-1}^{\theta_{1}}\int_{\theta_{1}}^{\theta_{2}}
\!\!\!\!\cdots\int_{\theta_{m-1}}^{\theta_{m}}
\vph(t_{0},\ldots,t_{m-1})\:dt_{m-1}\cdots dt_{0}.}
\end{array}
\ee{aa}
The operator $\AA$ is a central object of study below. Although we prove some general results (Propositions~7.6 and~7.7)
which are valid for every $m\ge 1$, our focus is ultimately on the cases $m=2$ and $m=3$, as
in Theorem~7.1.

We shall consider $\AA$ as acting on the space $C(T_m)$ of all
continuous functions $\vph:T_m\to\IR$, which is in contrast to
earlier sections where we worked with the space $X^{\wedge m}$.
However, note that $X^{\wedge m}$ is isometrically isomorphic to the subspace
$$
X^{\wedge m}_{\rm r}=\{\vph\in C(T_m)\:|\:\vph(\theta)=0\hbox{ whenever }\theta_j=\theta_{j+1}\hbox{ for some }
j\hbox{ satisfying }1\le j\le m-1\}
$$
of $C(T_m)$ consisting of all restrictions $\vph|_{T_m}$ of functions $\vph\in X^{\wedge m}$
to $T_m\subseteq[-1,0]^m$.
As such, we shall freely regard the function $u_m$ in~\eqref{u0} to be
an element of $C(T_m)$, in fact, $u_m\in X^{\wedge m}_{\rm r}\subseteq C(T_m)$. Also, without loss,
we may regard $\wwtilde(t,\tau)$ to be an operator on $X^{\wedge m}_{\rm r}$ rather than on $X^{\wedge m}$,
as we wish to compare $\wwtilde(t,\tau)$ with powers of the operator $\AA$.
Note that the ranges of $\AA_0$ and $\AA_1$ on $C(T_m)$ lie
in the subspace $X^{\wedge m}_{\rm r}$, and so the subspace $X^{\wedge m}_{\rm r}\subseteq C(T_m)$
is invariant under these operators.

Let us also denote the positive cone in $C(T_m)$ by
\be{pm}
C(T_m)^+=\{\vph\in C(T_m)\:|\:\vph(\theta)\ge 0\hbox{ for every }\theta\in T_m\}.
\ee{pm}
The crucial part in proving Theorem~7.1 is to show that
the operator $\AA$ is $u_m$-positive
with respect to $C(T_m)^+$. Indeed, we have the following result.
\vv

\noindent {\bf Proposition~7.3.} \bxx
Let $m\ge 2$ and suppose it is the case that the operator $\AA\in\call(C(T_m))$ given above is a $u_m$-positive
operator with respect to the cone $C(T_m)^+$, with $u_m$ as in~\eqref{u0}.
Assume that $\alpha,\beta:\IR\to\IR$ are as in the statement of Theorem~7.1. Then the conclusions
of Theorem~7.1 (and thus also of Proposition~7.2) hold, but with the value of $m$ chosen here. Further, we have that
$\eta_0=k_0$ for the quantities in the above definitions of $u_0$-positive operator
and $u_0$-positive linear process, corresponding to the operator $\AA$ and to the
linear process $\WW(t,\tau)$.
\exx
\vv

Proposition~7.3 implies that in order to prove Theorem~7.1, it is sufficient
to prove the following result.
\vv

\noindent {\bf Theorem~7.4.} \bxx
Let $m=2$ or $m=3$. Then the operator $\AA$ acting on $C(T_m)$ is $u_m$-positive
with respect to the cone $C(T_m)^+$, with $u_m$
as in~\eqref{u0}.
In fact, if $\vph\in C(T_m)^+\setminus\{0\}$ then $\AA^{k_0}\vph\sim u_m$ for $k_0=3$ if $m=2$,
and for $k_0=5$ if $m=3$.
\exx
\vv

\noindent {\bf Remark.} 
The fact that $u_m\in X^{\wedge m}_{\rm r}$, along with the invariance of $X^{\wedge m}_{\rm r}$
under $\AA$, implies for Theorem~7.4 that $\AA$, as an operator
on $X^{\wedge m}_{\rm r}$, is also $u_m$-positive for that space with respect to the cone
$C(T_m)^+\cap X^{\wedge m}_{\rm r}$.
\vv

\noindent {\bf Conjecture~B\@.}
Let $m\ge 4$. Then the conclusions of Theorem~7.4 hold for this $m$, but where $\AA^{k_0}\vph\sim u_m$
for $k_0=2m-1$.
\vv

It is clear from Proposition~7.3 that if Conjecture~B holds, then so does Conjecture~A\@.
\vv

\noindent {\bf Proof of Proposition~7.3.}
It is enough to prove only the conclusions of Theorem~7.1 pertaining
to equation~\eqref{1}, that is, for the linear process $\wwtilde(t,\tau)$.
As usual, the corresponding conclusions for equation~\eqref{0} can be obtained from these
using the conjugacy~\eqref{conj}.

By assumption, there exists $k_0\ge 1$ such that $\AA^{k_0}\vph\sim u_m$
for every $\vph\in C(T_m)^+\setminus\{0\}$. Let $\tau_1\le\tau_2$
and let $\eta_*\ge k_0$ be given. Then there exist $b_1$ and $b_2$ such that
\be{bbnd}
0<b_1\le (-1)^mb(t)\le b_2\qquad\hbox{for almost every }t\in[\tau_1,\tau_2+\eta_*],
\ee{bbnd}
from the bounds~\eqref{usc} and the formula~\eqref{xyz}. We shall work in the space $X^{\wedge m}_{\rm r}\subseteq C(T_m)$,
and we note that $C(T_m)^+\cap X^{\wedge m}_{\rm r}$ is a cone in that space. First
note that if $\vph\in C(T_m)^+\cap X^{\wedge m}_{\rm r}$, then from the formulas~\eqref{eta1x} and~\eqref{aa},
and the bounds~\eqref{bbnd}, we have that if $[\sig,\sig+1]\subseteq[\tau_1,\tau_2+\eta_*]$ then
$$
B_1\AA\vph\le (b_1^{m-1}\AA_0+b_1^m\AA_1)\vph\le \wwtilde(\sig+1,\sig)\vph\le(b_2^{m-1}\AA_0+b_2^m\AA_1)\vph\le B_2\AA\vph
$$
where
$$
B_1=\min\{b_1^{m-1},\,b_1^m\},\qquad
B_2=\max\{b_2^{m-1},\,b_2^m\},
$$
with $\le$ denoting the ordering in the cone $C(T_m)^+\cap X^{\wedge m}_{\rm r}$. It follows
by iteration that if $[\sig,\sig+k]\subseteq[\tau_1,\tau_2+\eta_*]$ for some integer $k\ge 1$, then
\be{kbnd}
B_1^k\AA^k\vph\le \wwtilde(\sig+k,\sig)\vph\le B_2^k\AA^k\vph.
\ee{kbnd}

Now let $\vph\in [C(T_m)^+\cap X^{\wedge m}_{\rm r}]\setminus\{0\}$ be given; we shall keep $\vph$
fixed for the remainder of the proof. For any $\eps>0$ denote
$$
\calu(\tau,\eta,\eps)=\{(\tau',\eta')\in[\tau_1,\tau_2]\times[k_0,\eta_*]\:|\:|\tau'-\tau|<\eps\mbox{ and }|\eta'-\eta|<\eps\}.
$$
Then given any $(\tau,\eta)\in[\tau_1,\tau_2]\times [k_0,\eta_*]$, let
$\psi=\wwtilde(\tau+\eta-k_0,\tau)\vph$, and note that $\psi\in [C(T_m)^+\cap X^{\wedge m}_{\rm r}]\setminus\{0\}$
where Proposition~4.4 is used. Thus there exists $\eps=\eps(\tau,\eta)>0$ such that if we define
$$
\psi_1(\theta)=\inf_{(\tau',\eta')\in\calu(\tau,\eta,\eps)}
[\wwtilde(\tau'+\eta'-k_0,\tau')\vph](\theta),
$$
for $\theta\in T_m$, then $\psi_1(\theta)>0$ for some $\theta$ and so $\psi_1\in [C(T_m)^+\cap X^{\wedge m}_{\rm r}]\setminus\{0\}$.
Fix such $\eps$ and let
$$
\psi_2(\theta)=\sup_{(\tau',\eta')\in\calu(\tau,\eta,\eps)}
[\wwtilde(\tau'+\eta'-k_0,\tau')\vph](\theta),
$$
and so also $\psi_2\in [C(T_m)^+\cap X^{\wedge m}_{\rm r}]\setminus\{0\}$.
Then for any $(\tau',\eta')\in\calu(\tau,\eta,\eps)$ we have that
\be{ord}
\psi_1\le \wwtilde(\tau'+\eta'-k_0,\tau')\vph\le\psi_2,
\ee{ord}
and upon applying the positive operator $\wwtilde(\tau'+\eta',\tau'+\eta'-k_0)$ to~\eqref{ord}, one obtains
$$
\wwtilde(\tau'+\eta',\tau'+\eta'-k_0)\psi_1\le \wwtilde(\tau'+\eta',\tau')\vph\le \wwtilde(\tau'+\eta',\tau'+\eta'-k_0)\psi_2.
$$
It follows, by~\eqref{kbnd}, that
\be{ord4}
Q_1u_m\le B_1^{k_0}\AA^{k_0}\psi_1\le \wwtilde(\tau'+\eta',\tau')\vph\le B_2^{k_0}\AA^{k_0}\psi_2\le Q_2u_m
\ee{ord4}
for some $Q_2\ge Q_1>0$, where the existence of $Q_1$ and $Q_2$
follows directly from the assumption that $\AA$ is $u_m$-positive, specifically, that
$\AA^{k_0}\psi_-\sim u_m$ and $\AA^{k_0}\psi_+\sim u_m$.

To complete the proof of the proposition, let us denote the constants in~\eqref{ord4} by $Q_{j,(\tau,\eta)}$, for $j=1,2$
for any $(\tau,\eta)\in[\tau_1,\tau_2]\times [k_0,\eta_*]$. We observe that the sets $\calu(\tau,\eta,\eps(\tau,\eta))$
for such $(\tau,\eta)$
form a (relatively) open cover of $[\tau_1,\tau_2]\times [k_0,\eta_*]$, so we may extract a finite subcover, corresponding to points
$(\tau_i,\eta_i)$ for $1\le i\le p$. Then upon setting
$C_1=\ds{\min_{1\le i\le p}}\{Q_{1,(\tau_i,\eta_i)}\}$ and
$C_2=\ds{\max_{1\le i\le p}}\{Q_{2,(\tau_i,\eta_i)}\}$,
we see that
$$
C_1u_m\le \wwtilde(\tau'+\eta',\tau')\vph\le C_2u_m
$$
for every $(\tau',\eta')\in[\tau_1,\tau_2]\times [k_0,\eta_*]$. With this, the proof is complete.~\fp
\vv

Moving toward the proof of Theorem~7.4,
we shall first obtain a pointwise upper bound for $|(\AA^k\vph)(\theta)|$ in Proposition~7.6 below,
and in fact we shall obtain such for every $m\ge 2$.
To this end we define functions $u_m^q\in C(T_m)$ by
\be{umq}
u_m^q(\theta)
=\bigg(\prodformtwo{\xxc}{1\le i,j\le m}{1\le j-i\le q}{(\theta_j-\theta_i)}\bigg)
\bigg(\prodformtwo{\xxd}{1\le i,j\le m-1}{j-i\ge m-q}{(1+\theta_i-\theta_j)}\bigg),
\ee{umq}
for $0\le q\le m-1$, recalling that an empty product takes the value $+1$. Here and below we shall always
assume $q$ is in this range, although we shall sometimes impose
additional restrictions on $q$. We assume that $m\ge 2$ and that $\theta\in T_m$.
Note that $u_m^0(\theta)\equiv 1$ identically, and that $u_m^{m-1}(\theta)=u_m(\theta)$ as in~\eqref{u0}.
Generally, as $q$ increases the polynomial~\eqref{umq} has more factors. Observe that
all factors in~\eqref{umq} are nonnegative and bounded above by $+1$.
Now define
$$
w_m(\theta)=
\prod_{j=1}^{m-1}(\theta_{j+1}-\theta_j),\qquad
\wt w_m(\theta)=
(1+\theta_1-\theta_{m-1})w_m(\theta_1,\ldots,\theta_m),
$$
and let polynomials $v_m^q$ and $\wt v_m^q$ be defined by
$$
\begin{array}{lclcl}
u_m^q(\theta) &\bb = &\bb v_m^q(\theta)w_m(\theta), &\bb \qquad &\bb \mbox{for }1\le q\le m-1,\\
\\
u_m^q(\theta) &\bb = &\bb \wt v_m^q(\theta)\wt w_m(\theta), &\bb \qquad &\bb \mbox{for }2\le q\le m-1,\\
\\
\wt v_m^1(\theta) &\bb \equiv &\bb 1 \mbox{ identically.}
\end{array}
$$
It is easy to check that $v_m^q$ and $\wt v_m^q$ are well-defined polynomials, as every factor of 
of $w_m$ and $\wt w_m$ occurs as a factor of the polynomial $u_m^q$ for
the indicated ranges of $q$.

Now let us take quantities $t_j$ for $0\le j\le m-1$ satisfying
\be{tthet}
t_0\in[-1,\theta_1],\qquad t_j\in[\theta_j,\theta_{j+1}]\qquad\hbox{for }1\le j\le m-1,
\ee{tthet}
as in the integrands of~\eqref{aa}. The following lemma provides a crucial estimate
needed for the proof of Proposition~7.6.
\vv

\noindent {\bf Lemma~7.5.} \bxx
With $m\ge 2$, let $\theta\in T_m$ and let $t_j$ for $0\le j\le m-1$ satisfy~\eqref{tthet}. Then
\be{uinq}
0\le u_m^q(t_{1},\ldots,t_{m-1},0)\le v_m^{q+1}(\theta),\qquad
0\le u_m^q(t_{0},\ldots,t_{m-1})\le \wt v_m^{q+1}(\theta),
\ee{uinq}
for $0\le q\le m-2$. Further,
\be{uinq2}
0\le u_m^{m-1}(t_{1},\ldots,t_{m-1},0)\le v_m^{m-1}(\theta),\qquad
0\le u_m^{m-1}(t_{0},\ldots,t_{m-1})\le \wt v_m^{m-1}(\theta),
\ee{uinq2}
holds.
\exx
\vv

\noindent {\bf Proof.}
In this proof care must be taken to ensure the correct ranges of the
indices $i$ and $j$, and it will be helpful to note that $i<j$ in many places.

Assume that $\theta\in T_m$ and that~\eqref{tthet} holds.
We begin by observing that
\be{tinq}
\begin{array}{rclcl}
0\le t_j-t_i
&\bb \le &\bb
\theta_{j+1}-\theta_i\le 1,
& \hbox{ for} &\bb 1\le i<j\le m-1,\\
\\
0\le t_j-t_0
&\bb \le &\bb
1+\theta_{j+1}-\theta_{m-1}\le 1,
\quad
& \hbox{ for} &\bb 1\le j\le m-2,\\
\\
0\le -t_i
&\bb \le &\bb
1+\theta_1-\theta_i\le 1,
& \hbox{ for} &\bb 1\le i\le m-1,\\
\\
0\le 1+t_i-t_j
&\bb \le &\bb
1+\theta_{i+1}-\theta_j\le 1,
& \hbox{ for} &\bb 0\le i<j\le m-1.
\end{array}
\ee{tinq}
We first establish~\eqref{uinq}. Suppose that $0\le q\le m-2$. Then from~\eqref{umq} we have that
\be{prod0}
u_m^q(t_1,\ldots,t_{m-1},0)
=\bigg(\prodformtwo{\xxe}{1\le i,j\le m-1}{1\le j-i\le q}{(t_j-t_i)}\bigg)
\bigg(\prodform{\xxf}{m-q\le i\le m-1}{(-t_i)}\bigg)
\bigg(\prodformtwo{\xxg}{1\le i,j\le m-1}{j-i\ge m-q}{(1+t_i-t_j)}\bigg).
\ee{prod0}
Using~\eqref{tinq} we see that
\be{first}
\prodformtwo{\xxh}{1\le i,j\le m-1}{1\le j-i\le q}{(t_j-t_i)}
\:\:\le\!\!\!\!
\prodformtwo{\xxi}{1\le i,j\le m-1}{1\le j-i\le q}{(\theta_{j+1}-\theta_i)}
\:\:=\!\!\!\!
\prodformtwo{\xxj}{1\le i,j\le m}{2\le j-i\le q+1}{(\theta_j-\theta_i)}
\ee{first}
for the first product in~\eqref{prod0}. We also have that
\be{third}
\prodformtwo{\xxk}{1\le i,j\le m-1}{j-i\ge m-q}{(1+t_i-t_j)}
\:\:\le\!\!\!\!
\prodformtwo{\xxl}{1\le i,j\le m-1}{j-i\ge m-q}{(1+\theta_{i+1}-\theta_j)}
\:\:=\!\!\!\!
\prodformtwo{\xxm}{2\le i,j\le m-1}{j-i\ge m-q-1}{(1+\theta_i-\theta_j)}
\ee{third}
for the third product in~\eqref{prod0}.
Using the inequality $-t_i\le 1+\theta_1-\theta_i$ from~\eqref{tinq} in the second product
in~\eqref{prod0}, and reindexing using $j$ instead of $i$,
we may combine this with~\eqref{third} to obtain
$$
\bigg(\prodform{\xxn}{m-q\le i\le m-1}{(-t_i)}\bigg)
\bigg(\prodformtwo{\xxo}{1\le i,j\le m-1}{j-i\ge m-q}{(1+t_i-t_j)}\bigg)
\:\:\le\!\!\!\!
\prodformtwo{\xxp}{1\le i,j\le m-1}{j-i\ge m-q-1}{(1+\theta_i-\theta_j)}.
$$
Combining this further with~\eqref{first}, we see that with~\eqref{prod0} this gives
$$
u_m^q(t_1,\ldots,t_{m-1},0)\leq
\bigg(\prodformtwo{\xxq}{1\le i,j\le m}{2\le j-i\le q+1}{(\theta_j-\theta_i)}\bigg)
\bigg(\prodformtwo{\xxr}{1\le i,j\le m-1}{j-i\ge m-q-1}{(1+\theta_i-\theta_j)}\bigg)
=v_m^{q+1}(\theta),
$$
to give the first half of~\eqref{uinq}.

Next observe that
\be{prod1}
u_m^q(t_0,\ldots,t_{m-1})
=\bigg(\prodformtwo{\xxs}{0\le i,j\le m-1}{1\le j-i\le q}{(t_j-t_i)}\bigg)
\bigg(\prodformtwo{\xxt}{0\le i,j\le m-2}{j-i\ge m-q}{(1+t_i-t_j)}\bigg).
\ee{prod1}
For the first product in~\eqref{prod1} we have, again using~\eqref{tinq}, that
\be{first2}
\begin{array}{lcl}
\ds{\prodformtwo{\xxu}{0\le i,j\le m-1}{1\le j-i\le q}{(t_j-t_i)}}
&\bb \le &\bb
\ds{\bigg(\prodformtwo{\xxv}{1\le i,j\le m-1}{1\le j-i\le q}{(\theta_{j+1}-\theta_i)}\bigg)
\bigg(\prodform{\xxw}{1\le j\le q}{(1+\theta_{j+1}-\theta_{m-1})}\bigg)}\\
\\
&\bb \le &\bb
\ds{\bigg(\prodformtwo{\xxx}{1\le i,j\le m}{2\le j-i\le q+1}{(\theta_j-\theta_i)}\bigg)
\bigg(\prodform{\xxy}{2\le i\le q}{(1+\theta_i-\theta_{m-1})}\bigg).}
\end{array}
\ee{first2}
Note that in the second inequality of~\eqref{first2}, we have used the estimate
$1+\theta_{q+1}-\theta_{m-1}\le 1$, which holds because
$q\le m-2$, and which allows us to drop the term $1+\theta_{q+1}-\theta_{m-1}$.
Now for the second product in~\eqref{prod1} we have that
\be{second2}
\prodformtwo{\xxz}{0\le i,j\le m-2}{j-i\ge m-q}{(1+t_i-t_j)}
\:\:\le\!\!\!\! \prodformtwo{\xxaa}{0\le i,j\le m-2}{j-i\ge m-q}{(1+\theta_{i+1}-\theta_j)}
\:\:=\!\!\!\!\prodformtwo{\xxbb}{1\le i,j\le m-2}{j-i\ge m-q-1}{(1+\theta_i-\theta_j)}.
\ee{second2}
If $q\ge 1$ then combining~\eqref{first2} and~\eqref{second2} gives
$$
u_m^q(t_0,\ldots,t_{m-1})\le
\bigg(\prodformtwo{\xxcc}{1\le i,j\le m}{2\le j-i\le q+1}{(\theta_j-\theta_i)\bigg)}
\bigg(\prodformtwo{\xxdd}{1\le i,j\le m-1}{m-3\ge j-i\ge m-q-1}{(1+\theta_i-\theta_j)}\bigg)
=\wt v_m^{q+1}(\theta),
$$
while if $q=0$ we have directly that
$$
u_m^0(t_0,\ldots,t_{m-1})=1=\wt v_m^1(\theta).
$$
This establishes the second half of~\eqref{uinq}.

We now prove~\eqref{uinq2}. This follows directly by noting that
$$
\begin{array}{l}
0\le u_m^{m-1}(t_{1},\ldots,t_{m-1},0)\le u_m^{m-2}(t_{1},\ldots,t_{m-1},0),\\
\\
0\le u_m^{m-1}(t_{0},\ldots,t_{m-1})\le u_m^{m-2}(t_{0},\ldots,t_{m-1}),
\end{array}
$$
and then applying~\eqref{uinq} for $q=m-2$.
With this the proof is complete.~\fp
\vv

\noindent {\bf Proposition~7.6.} \bxx
Let $m\ge 2$. Then we have the pointwise bounds
\be{ubnd}
0\le (\AA_iu_m^q)(\theta)\le u_m^{q+1}(\theta),\qquad
0\le (\AA_iu_m)(\theta)\le u_m(\theta),
\ee{ubnd}
for $0\le q\le m-2$ and $i=0,1$, for $\theta\in T_m$.
Thus for every $\vph\in C(T_m)$ we have the pointwise bound
\be{b3}
|(\AA^k\vph)(\theta)|\le 2^ku_m(\theta)\|\vph\|,
\ee{b3}
for $k\ge m-1$ and $\theta\in T_m$.
\exx
\vv

\noindent {\bf Proof.}
Let $\vph=u_m^q$ in~\eqref{aa} where $0\le q\le m-2$. The using~\eqref{uinq} in
Lemma~7.5, we have for every $\theta\in T_m$ that
$$
\begin{array}{lcl}
0\le(\AA_0u_m^q)(\theta)
&\bb = &\bb
\ds{\int_{\theta_{1}}^{\theta_{2}}\!\!\!\!\cdots\int_{\theta_{m-1}}^{\theta_{m}}
u_m^q(t_1,\ldots,t_{m-1},0)\:dt_{m-1}\cdots dt_1}\\
\\
&\bb \le &\bb
\ds{\int_{\theta_{1}}^{\theta_{2}}\!\!\!\!\cdots\int_{\theta_{m-1}}^{\theta_{m}}
v_m^{q+1}(\theta)\:dt_{m-1}\cdots dt_1=w_m(\theta)v_m^{q+1}(\theta)=u_m^{q+1}(\theta).}
\end{array}
$$
Similarly, if $1\le q\le m-2$ we have that
\be{a1q}
\begin{array}{lcl}
0\le(\AA_1u_m^q)(\theta)
&\bb = &\bb
\ds{\int_{-1}^{\theta_{1}}\int_{\theta_{1}}^{\theta_{2}}
\!\!\!\!\cdots\int_{\theta_{m-1}}^{\theta_{m}}
u_m^q(t_0,\ldots,t_{m-1})\:dt_{m-1}\cdots dt_0}\\
\\
&\bb \le &\bb
\ds{\int_{-1}^{\theta_{1}}\int_{\theta_{1}}^{\theta_{2}}
\!\!\!\!\cdots\int_{\theta_{m-1}}^{\theta_{m}}
\wt v_m^{q+1}(\theta)\:dt_{m-1}\cdots dt_0}\\
\\
&\bb = &\bb
\ds{\bigg(\frac{1+\theta_1}{1+\theta_1-\theta_{m-1}}\bigg)
\wt w_m(\theta)\wt v_m^{q+1}(\theta)\le \wt w_m(\theta)\wt v_m^{q+1}(\theta)=u_m^{q+1}(\theta).}
\end{array}
\ee{a1q}
If $q=0$ we again have~\eqref{a1q} except with an inequality~$\le$ in place of the final equal sign, as
$$
\wt w_m(\theta)\wt v_m^1(\theta)=\wt w_m(\theta)\le w_m(\theta)=u_m^1(\theta).
$$
This gives the first half of~\eqref{ubnd}. For the second half of~\eqref{ubnd},
involving $(\AA_i u_m)(\theta)$, one argues similarly except using~\eqref{uinq2} instead of~\eqref{uinq},
where we recall that $u_m=u_m^{m-1}$. We omit the details.

It follows that $0\le(\AA u_m^q)(\theta)\le 2u_m^{q+1}(\theta)$ if $0\le q\le m-2$, while
$0\le(\AA u_m^{m-1})(\theta)\le 2u_m^{m-1}(\theta)$, for every $\theta\in T_m$. Thus
\be{b1}
0\le(\AA^ku_m^0)(\theta)\le 2^ku_m^{\gam(k)}(\theta),\qquad\gam(k)=\min\{k,\,m-1\},
\ee{b1}
for every $k\ge 1$. Also, as $\AA$ is a positive operator with respect to the cone $C(T_m)^+$
in~\eqref{pm}, we have the pointwise bound
\be{b2}
|(\AA^k\vph)(\theta)|\le(\AA^k|\vph|)(\theta)\le[(\AA^ku_m^0)(\theta)]\|\vph\|
\ee{b2}
for every $\vph\in C(T_m)$, where we recall that $u_m^0(\theta)\equiv 1$ identically.
Combining~\eqref{b1} and~\eqref{b2} gives~\eqref{b3}, as desired.~\fp
\vv

Related to the operator $\AA$ is the operator $\BB$, which we define as $\BB=\BB_0+\BB_1$, where
\be{bformx}
\begin{array}{lcl}
(\BB_0\vph)(\theta) &\bb = &\bb
\ds{\frac{1}{u_m(\theta)}\int_{\theta_{1}}^{\theta_{2}}\!\!\!\!\cdots\int_{\theta_{m-1}}^{\theta_{m}}
u_m(t_{1},\ldots,t_{m-1},0)
\vph(t_{1},\ldots,t_{m-1},0)\:dt_{m-1}\cdots dt_{1},}\\
\\
(\BB_1\vph)(\theta) &\bb = &\bb
\ds{\frac{1}{u_m(\theta)}\int_{-1}^{\theta_{1}}\int_{\theta_{1}}^{\theta_{2}}
\!\!\!\!\cdots\int_{\theta_{m-1}}^{\theta_{m}}
u_m(t_{0},\ldots,t_{m-1})
\vph(t_{0},\ldots,t_{m-1})\:dt_{m-1}\cdots dt_{0}.}
\end{array}
\ee{bformx}
Formally, $\BB$ is conjugate to $\AA$ via the operator given by multiplication by $u_m$,
and $\BB$ will play a significant role in proving Theorem~7.4.
In particular, obtaining the required equivalence $\AA^{k_0}\vph\sim u_m$
is essentially the same as showing that $\BB^{k_0}\psi\sim 1$ where $\vph=u_m\psi$.
However, this conjugacy between $\AA$ and $\BB$ is only formal, as multiplication by $u_m$
is not an isomorphism on $C(T_m)$. Indeed, as described in the next section, and in particular
in Theorem~8.2, $\BB$ need not be a bounded operator on $C(T_m)$: It can happen that the function
$\BB\vph$ has discontinuities in $T_m$, so is not an element of $C(T_m)$, even when $\vph\in C(T_m)$.
On the other hand, it is always the case that if $\vph\in C(T_m)$ then the
function $\BB\vph$ is continuous almost everywhere on $T_m$, specifically, it is continuous
at each point $\theta\in O_m$ where
\be{om}
O_m=\{\theta\in T_m\:|\:\theta_j<\theta_{j+1}\mbox{ for every }j\mbox{ satisfying }1\le j\le m-1\}.
\ee{om}
In light of the estimate~\eqref{b3}, one might therefore wish to consider $\BB$ acting on the space $L^\infty(T_m)$
of bounded measurable functions; but $\BB\vph$ is not well-defined for general $\vph\in L^\infty(T_m)$ due to
the zero entry in the final argument of $\vph$ in the formula~\eqref{bformx} for $\BB_0\vph$.
However, if we define
\be{wm}
W_m=\{\vph\in L^\infty(T_m)\:|\:\vph\mbox{ is continuous at every point }\theta\in O_m\},
\ee{wm}
then $W_m\subseteq L^\infty(T_m)$ is a closed subspace, and
one easily sees that $\BB$ is indeed well-defined as an operator on $W_m$ with range in $W_m$, that is, $\BB\in\call(W_m)$.
We have have the following result.
\vv

\noindent {\bf Proposition~7.7.} \bxx
Let $m\ge 2$. Also let $W_m$ be defined by~\eqref{wm} with \eqref{om},
with the norm inherited from $L^\infty(T_m)$.
Then $\BB_0$ and $\BB_1$ in~\eqref{bformx}, and thus $\BB=\BB_0+\BB_1$, define bounded linear operators
on the space $W_m$, with $\|\BB_0\|,\|\BB_1\|\le 1$ and $\|\BB\|\le 2$.
\exx
\vv

\noindent {\bf Proof.}
The proof is very similar to the proof of Proposition~7.6, and entails
using the bounds~\eqref{uinq2} of Lemma~7.5 to estimate the integrals~\eqref{bformx}
just as before. We omit the details.~\fp
\vv

It is natural to ask what is the minimal closed invariant subspace $Y\subseteq W_m$
for the operator $\BB$ which contains $C(T_m)$. Such $Y$ would, in a sense, be the ``natural''
space on which $\BB$ acts and would be given by
\be{y}
Y=\overline{\bigcup_{n=0}^\infty Y_n},\qquad\mbox{where}\qquad
Y_0=C(T_m),\qquad
Y_{n+1}=\BB Y_n+C(T_m)\qquad\mbox{for }n\ge 0.
\ee{y}
Note here that $Y_0\subseteq Y_1\subseteq Y_2\subseteq\cdots\subseteq Y\subseteq W_m$.
In the next section, as part of our efforts to prove Theorem~7.4, we show that if $m=2$
then $Y=Y_0=C(T_2)$, while if $m=3$ then $Y=Y_1=C(T_3)\oplus V$ where $V$ is a certain
two-dimensional subspace of $W_m$.

\section{{\boldmath $u_0$}-Positivity for {\boldmath $m=2$} and {\boldmath $m=3$}}

Let us now specialize to the cases $m=2$ and $m=3$, as in Theorem~7.4. We retain all the conventions
and notation of the previous section. In working toward the proof of Theorem~7.4, our analysis here is largely
concerned with the operator $\BB$.

For convenience of notation, we denote the so-called average integral by
$$
\tint_a^bf(x)\:dx=\frac{1}{b-a}\int_a^bf(x)\:dx,\qquad
\hbox{with}\qquad\tint_a^af(x)\:dx=f(a).
$$
Note that for locally integrable $f$, the average integral is continuous as a function of $a$ and $b$ for $a\ne b$.
It is also continuous where $a=b$ provided that $f$ is continuous at this point.

We first consider the case $m=2$. Then
$u_2(\theta_1,\theta_2)=\theta_2-\theta_1$, and so
\be{71}
(\BB_0\vph)(\theta)=-\:\tint_{\theta_{1}}^{\theta_{2}}t_1\vph(t_{1},0)\:dt_{1},\qquad
(\BB_1\vph)(\theta)=\int_{-1}^{\theta_{1}}\tint_{\theta_{1}}^{\theta_{2}}
(t_1-t_0)\vph(t_{0},t_{1})\:dt_{1}\:dt_{0},
\ee{71}
as in~\eqref{bformx}. It is immediate that $\BB_0$ and $\BB_1$, and thus $\BB$, are bounded linear operators on the space $C(T_2)$,
and so $Y=Y_0=C(T_2)$ for the spaces $Y$ and $Y_0$ in~\eqref{y}.
One sees moreover that $\BB_0$, $\BB_1$, and $\BB$ are positive operators on $C(T_2)$ with respect to the cone $C(T_2)^+$.
(Keep in mind that $t_1\le 0$ for the integrand
in the formula for $\BB_0$.)

We have the following result.
\vv

\noindent {\bf Proposition~8.1.} \bxx
Let $\vph\in C(T_2)^+\setminus\{0\}$. Then for every $k\ge 3$ there
exists $B_k>0$ such that $(\BB^k\vph)(\theta)\ge B_k$ for every
$\theta\in T_2$. Thus the operator $\BB$
with $m=2$ and acting on $C(T_2)$ is $u_0$-positive with respect to
the cone $C(T_2)^+$,
where $u_0(\theta)\equiv 1$ identically on $T_2$.
\exx
\vv

\noindent {\bf Proof.}
Clearly $u_0$-positivity follows from the existence  of the lower bound $B_3$,
as we have the pointwise upper bounds $|(\BB^3\vph)(\theta)|\le\|\BB^3\vph\|$.
Also, it is sufficient to prove the existence only of $B_3$, as the constants $B_k$ for $k\ge 4$ follow
directly by induction, using the positivity of $\BB$. Indeed, having obtained $B_j$ for $3\le j\le k$,
we obtain a lower bound $B_{k+1}$ for $|(\BB^{k+1})(\theta)|$ by applying $\BB^3$ to
the function $\BB^{k-2}\vph\in C(T_2)^+\setminus\{0\}$.

To show that $B_3$ exists, it is enough, due to the continuity of $\BB^3\vph$, to show that if
$\vph\in C(T_2)^+\setminus\{0\}$ then we have strict positivity $(\BB^3\vph)(\theta)>0$ for every $\theta\in T_2$.
To this end let
$$
L=\{\theta\in T_2\:|\:\theta_2=0\}=[-1,0]\times\{0\},
$$
which is the upper boundary of the set $T_2\subseteq \IR^2$. Then it is enough to prove that the following three facts
hold for every $\vph\in C(T_2)^+$.
\begin{itemize}
\item[{(1)}] If $\vph(\theta)>0$ for some $\theta\in T_2$, then $(\BB\vph)(\wt\theta)>0$ for some $\wt\theta\in L$;
\item[{(2)}] if $\vph(\theta)>0$ for some $\theta\in L$, then $(\BB\vph)(\wt\theta)>0$ for every $\wt\theta\in L$; and
\item[{(3)}] if $\vph(\theta)>0$ for every $\theta\in L$, then $(\BB\vph)(\wt\theta)>0$ for every $\wt\theta\in T_2$.
\end{itemize}
The proofs of properties~(1) through~(3) follow easily from the formulas~\eqref{71}.
With $\vph\in C(T_2)^+$, if $\vph(\theta)>0$ for some $\theta=(\theta_1,\theta_2)\in T_2$
we may assume without loss that $-1<\theta_1<\theta_2\le 0$.
Then $(\BB_1\vph)(\theta_1,0)>0$ holds, in particular because the integrand $(t_1-t_0)\vph(t_0,t_1)$ in~\eqref{71},
which is nonnegative throughout the range $-1\le t_0\le\theta_1\le t_1\le 0$, is strictly positive at $(t_0,t_1)=(\theta_1,\theta_2)$.
With this,~(1) is established.

Now suppose that $\vph(\theta)>0$ for some $\theta=(\theta_1,0)\in L$. Then $(\BB_0\vph)(\wt\theta_1,0)>0$
for every $\wt\theta_1$ satisfying $-1\le\wt\theta_1\le\theta_1$ and $\wt\theta_1\ne 0$, and $(\BB_1\vph)(\wt\theta_1,0)>0$
for every $\wt\theta_1$ satisfying $\theta_1\le\wt\theta_1\le 0$ and $\wt\theta_1\ne -1$. In any case, $(\BB\vph)(\wt\theta)>0$
for every $\wt\theta=(\wt\theta_1,0)\in L$. This establishes~(2).

Finally suppose that $\vph(\theta)>0$ for every $\theta\in L$. Then $(\BB_0\vph)(\wt\theta)>0$ for every
$\wt\theta\in T_2$ except $\wt\theta=(0,0)$. However, $(\BB_1\vph)(0,0)>0$ for this point. With this,~(3)
is established and the result is proved.~\fp
\vv

Let us now consider the case $m=3$, so $\theta=(\theta_1,\theta_2,\theta_3)\in T_3$. Here
$$
u_3(\theta)=(\theta_2-\theta_1)(\theta_3-\theta_2)(\theta_3-\theta_1)(1+\theta_1-\theta_2).
$$
We introduce the functions
\be{nu}
\nu_0(\theta)=\frac{-(\theta_2+\theta_3)}{1+\theta_1-\theta_2},\qquad
\nu_1(\theta)=\frac{1+\theta_1}{1+\theta_1-\theta_2},
\ee{nu}
which will play a key role in our analysis.
Observe that due to the ordering of the $\theta_j$ in the definition~\eqref{tri} of $T_3$,
the functions $\nu_0$ and $\nu_1$ are well-defined and continuous everywhere in $T_3$ except
at the point $\theta=(-1,0,0)$. Further,
we have the bounds $0\le\nu_0(\theta_1,\theta_2,\theta_3)\le 2$ and
$0\le\nu_1(\theta_1,\theta_2,\theta_3)\le 1$ throughout $T_3\setminus\{(-1,0,0)\}$,
so $\nu_0,\nu_1\in W_3$, where we recall the definition of $W_m$ in~\eqref{wm}.

After a short calculation one sees from~\eqref{bformx} that
\be{bm3}
\begin{array}{lcl}
(\BB_0\vph)(\theta) &\bb = &\bb
\ds{\tint_{\theta_{1}}^{\theta_{2}}\tint_{\theta_{2}}^{\theta_{3}}
\Phi_0(t_1,t_2,\theta_1,\theta_2,\theta_3)
\vph(t_{1},t_{2},0)\:dt_{2}\:dt_{1},}\\
\\
(\BB_1\vph)(\theta) &\bb = &\bb \nu_1(\theta)(\wt\BB_1\vph)(\theta),\qquad\hbox{where}\\
\\
(\wt\BB_1\vph)(\theta) &\bb = &\bb
\ds{\tint_{-1}^{\theta_{1}}\tint_{\theta_{1}}^{\theta_{2}}\tint_{\theta_{2}}^{\theta_{3}}
\Phi_1(t_0,t_1,t_2,\theta_1,\theta_3)
\vph(t_{0},t_{1},t_{2})\:dt_{2}\:dt_1\:dt_{0},}
\end{array}
\ee{bm3}
for any $\vph\in W_3$, and where the kernels $\Phi_0$ and $\Phi_1$ are given by
\be{phiker}
\begin{array}{lcl}
\Phi_0(t_1,t_2,\theta_1,\theta_2,\theta_3) &\bb = &\bb
\ds{\bigg(\frac{t_2-t_1}{\theta_3-\theta_1}\bigg)t_1\bigg(\frac{t_2}{1+\theta_1-\theta_2}\bigg)(1+t_1-t_2),}\\
\\
\Phi_1(t_0,t_1,t_2,\theta_1,\theta_3) &\bb = &\bb
\ds{(t_1-t_0)\bigg(\frac{t_2-t_1}{\theta_3-\theta_1}\bigg)(t_2-t_0)(1+t_0-t_1).}
\end{array}
\ee{phiker}
Note that we have grouped like terms in the kernels~\eqref{phiker}, so that each ratio in these
formulas is at most $+1$ in absolute value. In particular, we have that
$$
0\le t_2-t_1\le\theta_3-\theta_1,\qquad 0\le -t_2\le 1+\theta_1-\theta_2,
$$
and so
\be{phibnd}
0\le\Phi_0(t_1,t_2,\theta_1,\theta_2,\theta_3)\le 1,\qquad
0\le \Phi_1(t_0,t_1,t_2,\theta_1,\theta_3)\le 1,
\ee{phibnd}
as long as $-1\le t_0\le \theta_1\le t_1\le \theta_2\le t_2\le\theta_3\le 0$.
This confirms the conclusion of Proposition~7.7 in the case $m=3$,
in particular that $\BB_0$, $\BB_1$, and $\BB$
define bounded linear operators on $W_3$. However, in contrast
to the case $m=2$ above, we shall see that here $\BB$ is not an operator on $C(T_3)$, as $\BB\vph$ is
not in general a continuous function on $T_3$ even if $\vph$ is continuous there.
Instead, the following result holds.
\vv

\noindent {\bf Theorem~8.2.} \bxx
Let $V\subseteq W_3$ denote the two-dimensional vector space spanned by the functions
$\nu_0$ and $\nu_1$ in~\eqref{nu},
and let $C_{0,V}\subseteq C_V\subseteq W_3$ be defined as
$$
C_V=C(T_3)\oplus V,\qquad
C_{0,V}=C_0(T_3)\oplus V,\qquad
C_0(T_3)=\{\vph\in C(T_3)\:|\:\vph(-1,0,0)=0\}.
$$
Then the space $C_V$ is invariant under the operators $\BB_0$, $\BB_1$, and thus $\BB$, and moreover,
the ranges of these operators on $C_V$ are contained in $C_{0,V}$.
In particular, if $\vph\in C_V$ then
\be{49}
\begin{array}{l}
(\BB\vph)(\theta)=Q_0\nu_0(\theta)+Q_1\nu_1(\theta)+\psi(\theta),\\
\\
\ds{Q_0=\frac{1}{2}\:\int_{-1}^{0}t^2(1+t)\vph(t,0,0)\:dt,\qquad
Q_1=\int_{-1}^{0}t^2(1+t)\vph(-1,t,0)\:dt,}
\end{array}
\ee{49}
where $\psi\in C_0(T_3)$.
Further, we have $Y=Y_1=C_V$ for the spaces $Y$ and $Y_1$ defined in~\eqref{y}.
\exx
\vv

\noindent {\bf Remark.}
The above theorem implies that although $\BB\vph$ need not be continuous even if $\vph$ is continuous,
the discontinuities of $\BB\vph$ can only be of a special form and located at the specific point $(-1,0,0)$ on
the boundary of $T_3$.
The analogous issue for $m\ge 4$, namely a description or classification
of the possible discontinuities that can arise for iterates
$\BB^k\vph$ where $\vph\in C(T_m)$, or more generally a characterization of the space $Y$,
should be relevant to the conjectures stated earlier, as well as being
an interesting question in its own right.
\vv

Now define the set
$$
C_V^+=\{\vph\in C_V\:|\:\vph(\theta)\ge 0\mbox{ for every }\theta\in T_3\setminus\{(-1,0,0)\}\},
$$
which is a cone in $C_V$. The following result is the analog of Proposition~8.1 for $m=3$.
\vv

\noindent {\bf Proposition~8.3.} \bxx
Let $\vph\in C_V^+\setminus\{0\}$. Then for every $k\ge 5$ there
exists $B_k>0$ such that $(\BB^k\vph)(\theta)\ge B_k$ for every
$\theta\in T_3\setminus\{(-1,0,0)\}$. Thus the operator $\BB$
with $m=3$ and acting on $C_V$ is $u_0$-positive with respect to
the cone $C_V^+$,
where $u_0(\theta)\equiv 1$ identically on $T_3$.
\exx
\vv

A number of preliminary results are needed before proving Theorem~8.2 and Proposition~8.3. We begin
by examining the continuity properties of $\BB_0\vph$ and $\wt\BB_1\vph$ in $T_3$
for $\vph\in C_V$.
For such $\vph$, it is clear from the form~\eqref{phiker} of the kernels $\Phi_i$ and from the formulas~\eqref{nu}
for $\nu_0$ and $\nu_1$
that the only possible points $\theta\in T_3$
at which $\BB_0\vph$ is discontinuous are where either $\theta_3-\theta_1=0$
or where $1+\theta_1-\theta_2=0$, and that the only possible points of discontinuity of $\wt\BB_1\vph$ are where $\theta_3-\theta_1=0$.
Note that $\theta_3-\theta_1=0$ for $\theta\in T_3$ if and only if $\theta=(\theta_*,\theta_*,\theta_*)$ for some $\theta_*\in[-1,0]$.
Also observe that $1+\theta_1-\theta_2=0$ for $\theta\in T_3$ if and only if $\theta=(-1,0,0)$.
The following lemma describes these continuity properties of $\BB_0\vph$ and $\wt\BB_1\vph$ at these points.
\vv

\noindent {\bf Lemma~8.4.} \bxx
Let $\vph\in C_V$. Then
the only possible point $\theta\in T_3$ of discontinuity of $\BB_0\vph$ is $\theta=(-1,0,0)$,
while $\wt\BB_1\vph$ is continuous throughout $T_3$, that is, $\wt\BB_1\vph\in C(T_3)$.
Further,
\be{btrip}
\begin{array}{lcl}
(\BB_0\vph)(\theta_*,\theta_*,\theta_*) &\bb = &\bb
\ds{\frac{\theta_*^2\vph(\theta_*,\theta_*,0)}{2},}\\
\\
(\wt\BB_1\vph)(\theta_*,\theta_*,\theta_*) &\bb = &\bb
\ds{\frac{1}{2}\:\tint_{-1}^{\theta_*}(\theta_*-t)^2(1+t-\theta_*)\vph(t,\theta_*,\theta_*)\:dt,}
\end{array}
\ee{btrip}
for every $\theta_*\in[-1,0]$.
\exx
\vv

\noindent {\bf Proof.}
From the remarks preceding the statement of the lemma, all that is necessary is to prove continuity of
$\BB_0\vph$ and $\wt\BB_1\vph$ at each point of the form
$(\theta_*,\theta_*,\theta_*)$ in $T_3$ and to establish the formulas~\eqref{btrip}.
We present only the proof for $\BB_0$, as the proof for $\wt\BB_1$
is similar. With $\vph\in C_V$ fixed, let $\gam_1=\theta_2-\theta_1$ and $\gam_2=\theta_3-\theta_2$,
which are nonnegative quantities for $\theta\in T_m$.
Making the change of variables $t_1=\theta_1+\tau_1\gam_1$ and $t_2=\theta_1+\gam_1+\tau_2\gam_2$
in~\eqref{bm3}, we obtain
$$
\begin{array}{lcl}
(\BB_0\vph)(\theta) &\bb = &\bb
\ds{\int_{0}^{1}\int_{0}^{1}
\Phi_0(\theta_1+\tau_1\gam_1,\:\theta_1+\gam_1+\tau_2\gam_2,\:\theta_1,\:\theta_1+\gam_1,\:\theta_1+\gam_1+\gam_2)}\\
\\
&\bb &\bb
\qquad\qquad\times\:\:\vph(\theta_1+\tau_1\gam_1,\:\theta_1+\gam_1+\tau_2\gam_2,\:0)\:d\tau_{2}\:d\tau_{1}\\
\\
&\bb = &\bb
\ds{\int_{0}^{1}\int_{0}^{1}
\bigg(\frac{(1-\tau_1)\gam_1+\tau_2\gam_2}{\gam_1+\gam_2}\bigg)
S(\tau_1,\tau_2,\theta_1,\gam_1,\gam_2)
\vph(\theta_1+\tau_1\gam_1,\:\theta_1+\gam_1+\tau_2\gam_2,\:0)\:d\tau_{2}\:d\tau_{1}}
\end{array}
$$
where
$$
S(\tau_1,\tau_2,\theta_1,\gam_1,\gam_2)=\frac{(\theta_1+\tau_1\gam_1)(\theta_1+\gam_1+\tau_2\gam_2)(1+(\tau_1-1)\gam_1-\tau_2\gam_2)}{1-\gam_1}.
$$
This formula is valid throughout $T_3$ as long as the coordinates $\theta_j$ for $j=1,2,3$ are not all equal
and $\theta\ne(-1,0,0)$, equivalently as long as $\gam_1+\gam_2>0$ and $\gam_1\ne 1$.
Upon letting $\theta=(\theta_1,\theta_2,\theta_3)$ approach a given point $(\theta_*,\theta_*,\theta_*)$ in $T_3$ (say, along a sequence),
one sees that $\gam_1$ and $\gam_2$ approach $0$,
hence $S(\tau_1,\tau_2,\theta_1,\gam_1,\gam_2)$ approaches $\theta_*^2$ and
$\vph(\theta_1+\tau_1\gam_1,\theta_1+\gam_1+\tau_2\gam_2,0)$ approaches $\vph(\theta_*,\theta_*,0)$, uniformly in the range of
integration. The fact that
$$
\int_0^1\int_0^1 \bigg(\frac{(1-\tau_1)\gam_1+\tau_2\gam_2}{\gam_1+\gam_2}\bigg)\:d\tau_2\:d\tau_1=\frac{1}{2},
$$
with an integrand which is bounded uniformly for nonnegative $\gam_1$ and $\gam_2$, implies that
$$
(\BB_0\vph)(\theta)\to\frac{\theta_*^2\vph(\theta_*,\theta_*,0)}{2}=(\BB_0\vph)(\theta_*,\theta_*,\theta_*),
$$
where the above equality may be taken as the defnition of $(\BB_0\vph)(\theta_*,\theta_*,\theta_*)$.
This gives the first equation in~\eqref{btrip}.
We omit the proof of the second equation in~\eqref{btrip}, which is similar.
With this, the result is proved.~\fp
\vv

\noindent {\bf Remark.} Although $\gam_1\to 0$ and $\gam_2\to 0$ in the above proof,
there is no assumption about the relative rates at which these quantities converge.
Consequently, the ratio $((1-\tau_1)\gam_1+\tau_2\gam_2)/(\gam_1+\gam_2)$ in the integrand
above need not have a pointwise limit in $(\tau_1,\tau_2)$ as $\gam_1,\gam_2\to 0$.
\vv

The next two lemmas give partial information on continuity properties of $\BB_0\vph$ near $(-1,0,0)$.
\vv

\noindent {\bf Lemma~8.5.} \bxx
Let $\vph\in C(T_3)$ and suppose that $\vph(\theta_1,0,0)\equiv 0$ identically for $\theta_1\in (-1,0]$.
Then $(\BB_0\vph)(\theta)$ is continuous at each $\theta\in T_3$, that is, $\BB_0\vph\in C(T_3)$.
Moreover, $(\BB_0\vph)(-1,0,0)=0$, and thus $\BB_0\vph\in C_0(T_3)$.
\exx
\vv

\noindent {\bf Proof.}
From Lemma~8.4, the only point at which $\BB_0\vph$ can fail to be continuous
is $\theta=(-1,0,0)$. Now let $\theta\in T_3\setminus\{(-1,0,0)\}$ be such that $\theta_1\ne\theta_3$. Then
from the formula~\eqref{bm3} and the bounds~\eqref{phibnd} we have that
$$
|\BB_0\vph(\theta)|\le
\sup_{-1\le t_1\le\theta_2\le t_2\le 0}|\vph(t_1,t_2,0)|=\del(\theta_2),
$$
where the above equality serves as the definition of $\del(\theta_2)$.
This bound also holds when $\theta_1=\theta_3$, and thus throughout $T_3\setminus\{(-1,0,0)\}$,
as $\del(\theta_2)$ depends continuously on $\theta_2$. Now
$\del(0)=0$ from the assumptions on $\vph$, and so the result follows directly.~\fp
\vv

\noindent {\bf Lemma~8.6.} \bxx
Let $\vph\in C(T_3)$ and suppose that $\vph(\theta_1,\theta_2,0)\equiv\vph(\theta_1,0,0)$
identically for every $(\theta_1,\theta_2)\in T_2$. Then
\be{nnu1}
(\BB_0\vph)(\theta)=Q\nu_0(\theta)+\psi(\theta),\qquad
Q=\frac{1}{2}\:\int_{-1}^{0}t^2(1+t)\vph(t,0,0)\:dt,
\ee{nnu1}
for some $\psi\in C_0(T_3)$.
\exx
\vv

\noindent {\bf Proof.}
We first claim that
\be{pwr}
\tint_{\theta_{2}}^{\theta_{3}}
(t_2-t_1)t_1t_2(1+t_1-t_2)\:dt_2=
\frac{1}{2}\bigg(-t_1^2(1+t_1)(\theta_2+\theta_3)+R(t_1,\theta_2,\theta_3)\bigg),
\ee{pwr}
where the polynomial $R$ satisfies
\be{rest}
R(t_1,\theta_2,\theta_3)=O(t_1(\theta_2^2+\theta_3^2))
\ee{rest}
near the origin
in $\IR^3$. Indeed, this can be readily verified by direct calculation, by expanding
the integrand in~\eqref{pwr} in powers of $t_2$
about the origin. We omit the details.
We next observe that
$$
(\BB_0\vph)(\theta)=
\frac{1}{2(\theta_3-\theta_1)(1+\theta_1-\theta_2)}\:\tint_{\theta_{1}}^{\theta_{2}}
\bigg(-t_1^2(1+t_1)(\theta_2+\theta_3)+R(t_1,\theta_2,\theta_3)\bigg)
\vph(t_{1},0,0)\:dt_{1},
$$
which follows immediately from~\eqref{bm3}, \eqref{phiker}, and~\eqref{pwr},
using the assumption on $\vph$. Denoting
$$
\begin{array}{lcl}
\wt\psi(\theta) &\bb = &\bb \ds{\frac{1}{2(\theta_3-\theta_1)}\:\tint_{\theta_{1}}^{\theta_{2}}
t_1^2(1+t_1)\vph(t_{1},0,0)\:dt_{1},}\\
\\
\ovrl\psi(\theta) &\bb = &\bb
\ds{\frac{1}{2(\theta_3-\theta_1)(1+\theta_1-\theta_2)}\:\tint_{\theta_{1}}^{\theta_{2}}
R(t_1,\theta_2,\theta_3)\vph(t_{1},0,0)\:dt_{1},}
\end{array}
$$
we have that $(\BB_0\vph)(\theta)=\nu_0(\theta)\wt\psi(\theta)+\ovrl\psi(\theta)$.

Certainly $\wt\psi$ is continuous in a neighborhood of the point $(-1,0,0)$ in $T_3$.
Also, $\ovrl\psi$ is continuous in some neighborhood of $(-1,0,0)$, in fact with $\ovrl\psi(-1,0,0)=0$,
in light of the estimate~\eqref{rest}.
Letting $\psi(\theta)=\nu_0(\theta)[\wt\psi(\theta)-Q]+\ovrl\psi(\theta)$ where $Q=\wt\psi(-1,0,0)$,
we have that~\eqref{nnu1} holds. Also,
$\psi$ is continuous in some neighborhood
of $(-1,0,0)$ in $T_3$, and $\psi(-1,0,0)=0$, where the boundedness of $\nu_0$ near that point is used.
Further, $\BB_0\vph$ and $\nu_0$ are continuous at every point of
$T_3\setminus\{(-1,0,0)\}$, using Lemma~8.4, so $\psi$ is also continuous there, by~\eqref{nnu1}.
Thus $\psi\in C_0(T_3)$, as claimed.~\fp
\vv

\noindent {\bf Lemma~8.7.} \bxx
We have that $\BB_0\nu_0\in C_0(T_3)$, while
$$
\BB_0\nu_1=\frac{1}{24}\nu_0+\mu_0
$$
with $\mu_0\in C_0(T_3)$.
\exx
\vv

\noindent {\bf Proof.} This is a straightforward calculation.~\fp
\vv

\noindent {\bf Remark.}
In the spirit of the above lemma, one can also check that $\BB_1\nu_1\in C_0(T_3)$, while
$$
\BB_1\nu_0=\frac{1}{12}\nu_1+\mu_1
$$
with $\mu_1\in C_0(T_3)$.
\vv

\noindent {\bf Proposition~8.8.} \bxx
Let $\vph\in C_V$. Then $\BB_0\vph$ has the form~\eqref{nnu1} where $\psi\in C_0(T_3)$,
with $Q$ as in that formula.
\exx
\vv

\noindent {\bf Proof.}
If either $\vph=\nu_0$ or $\vph=\nu_1$, the result follows from Lemma~8.7 along with a calculation
of the constant $Q$. Thus it is enough to consider $\vph\in C(T_3)$. For such $\vph$
write $\vph(\theta)=\wt\vph(\theta)+\ovrl\vph(\theta)$ where
$$
\wt\vph(\theta_1,\theta_2,\theta_3)=\vph(\theta_1,\theta_2,\theta_3)-\vph(\theta_1,\theta_3,\theta_3),\qquad
\ovrl\vph(\theta_1,\theta_2,\theta_3)=\vph(\theta_1,\theta_3,\theta_3).
$$
Then $\wt\vph$ and $\ovrl\vph$ satisfy the conditions of Lemmas~8.5 and~8.6, respectively.
Further, $\ovrl\vph(\theta_1,0,0)\equiv\vph(\theta_1,0,0)$ identically on $(-1,0]$.
The result follows immediately.~\fp
\vv

\noindent {\bf Proof of Theorem~8.2.}
Let $\vph\in C_V$.
We have that $(\BB_0\vph)(\theta)=Q_0\nu_0(\theta)+\psi_0(\theta)$ with $Q_0$ as in~\eqref{49}
and where $\psi_0\in C_0(T_3)$, by Proposition~8.8. Also, by Lemma~8.4 we have that
$\wt\BB_1\vph\in C(T_3)$, and from~\eqref{bm3} and~\eqref{phiker}
we have that $(\wt\BB_1\vph)(-1,0,0)=Q_1$
for $Q_1$ as in~\eqref{49}. Thus upon letting $\psi_1(\theta)=\nu_1(\theta)[(\wt\BB_1\vph)(\theta)-Q_1]$,
we have from~\eqref{bm3} that $(\BB_1\vph)(\theta)=Q_1\nu_1(\theta)+\psi_1(\theta)$, so upon letting
$\psi(\theta)=\psi_0(\theta)+\psi_1(\theta)$ we have~\eqref{49}. The fact that
$\psi\in C_0(T_3)$, or equivalently, that $\psi_1\in C_0(T_3)$,
follows directly from the definition of $\psi_1$ using
the continuity of $\wt\BB_1\vph$ and the choice of $Q_1$.

To prove the final claim in the statement of the theorem, it is enough to show that
every pair of numbers $Q_0,Q_1\in\IR$ as in~\eqref{49} can be achieved for some $\vph\in C(T_3)$.
However, this follows easily from the explicit formulas~\eqref{49} for $Q_0$ and $Q_1$.~\fp
\vv

\noindent {\bf Lemma~8.9.} \bxx
Let $\vph\in C_{0,V}$ have the form
$$
\vph(\theta)=Q_0\nu_0(\theta)+Q_1\nu_1(\theta)+\psi(\theta)
$$
for $\theta\in T_3\setminus\{(-1,0,0)\}$ where $\psi\in C_0(T_3)$.
Assume that $\vph(\theta)>0$ for every $\theta\in T_3\setminus\{(-1,0,0)\}$
and also that $Q_0>0$ and $Q_1>0$. Then there exists $B>0$ such that
\be{79m}
\vph(\theta)\ge B
\ee{79m}
for every $\theta\in T_3\setminus\{(-1,0,0)\}$.
\exx
\vv

\noindent {\bf Proof.}
Denote $\theta_0=(-1,0,0)$. Then for any $r$ satisfying $0<r\le 1$, let
$$
\del(r)=\inf_{|\theta-\theta_0|\ge r}\vph(\theta),\qquad
\eps(r)=\sup_{0\le |\theta-\theta_0|\le r} |\psi(\theta)|,
$$
where $|\cdot|$ denotes the euclidean distance in $\IR^3$ and $\theta\in T_3$.
Since $\vph$ is continuous and positive throughout $T_m\setminus\{\theta_0\}$,
it follows that $\del(r)$
is positive and depends continuously on $r$. Also, since $\psi$ is continuous in $T_3$ and $\psi(\theta_0)=0$, we
have that $\eps(r)$ also depends continuously on $r$, and that $\ds{\lim_{r\to 0}}\eps(r)=0$.
(We note that it need not be the case that the function $\psi$ is nonnegative everywhere.)
Therefore,
$$
\begin{array}{lcl}
\ds{\inf_{0<|\theta-\theta_0|\le r}\vph(\theta)} &\bb \ge &\bb
\ds{\inf_{0<|\theta-\theta_0|\le r}\bigg(Q_0\nu_0(\theta)+Q_1\nu_1(\theta)\bigg)-\eps(r)}\\
\\
&\bb = &\bb
\ds{\inf_{0<|\theta-\theta_0|\le r}\bigg(\frac{-Q_0(\theta_2+\theta_3)+Q_1(1+\theta_1)}{1+\theta_1-\theta_2}\bigg)-\eps(r)
\ge\min\{Q_0,Q_1\}-\eps(r).}
\end{array}
$$
By choosing $r$ sufficiently small that $\eps(r)<\min\{Q_0,Q_1\}$, one sees immediately that
the desired inequality~\eqref{79m} holds throughout $T_m\setminus\{\theta_0\}$ with $B=\min\{\del(r),\:\min\{Q_0,Q_1\}-\eps(r)\}$.~\fp
\vv

\noindent {\bf Proof of Proposition~8.3.}
Just as in the proof of Proposition~8.1, it is sufficient here only to prove
the existence of $B_5$. Toward this end, let us first define the sets
$$
\begin{array}{lcl}
L_0 &\bb = &\bb T_3\setminus\{(-1,0,0)\},\\
\\
L_1 &\bb = &\bb \{\theta\in T_3\setminus\{(-1,0,0)\}\:|\:\theta_3=0\},\\
\\
L_2 &\bb = &\bb \{\theta\in T_3\setminus\{(-1,0,0)\}\:|\:\theta_2=\theta_3=0\},\\
\\
L_3 &\bb = &\bb \{(0,0,0)\},
\end{array}
$$
observing that $L_3\subseteq L_2\subseteq L_1\subseteq L_0$. Also define
$$
L_\del=\{\theta\in T_3\setminus\{(-1,0,0)\}\:|\:\theta_2\in[-\del,0)\mbox{ and }\theta_3=0\}
$$
for any $\del>0$. We claim the following facts hold for every $\vph\in C_V^+$.
\begin{itemize}
\item[{(1)}] If $\vph(\theta)>0$ for some $\theta\in L_0$, then $(\BB\vph)(\wt \theta)>0$ for some $\wt\theta\in L_1$;
\item[{(2)}] if $\vph(\theta)>0$ for some $\theta\in L_1$, then $(\BB\vph)(\wt \theta)>0$ for some $\wt\theta\in L_2$;
\item[{(3)}] if $\vph(\theta)>0$ for some $\theta\in L_2$,
then there exists $\del>0$ such that $(\BB\vph)(\wt \theta)>0$ for every $\wt\theta\in L_\del$;
\item[{(4)}] if there exists $\del>0$ such that $\vph(\theta)>0$ for every $\theta\in L_\del$,
then $(\BB\vph)(\wt \theta)>0$ for every $\wt\theta\in L_1\setminus L_2$; and
\item[{(5)}] if $\vph(\theta)>0$ for every $\theta\in L_1\setminus L_2$,
then $(\BB\vph)(\wt \theta)>0$ for every $\wt\theta\in L_0\setminus L_2$.
\end{itemize}
Additionally, we claim the following facts hold for every $\vph\in C_V^+$.
\begin{itemize}
\item[{(3$'$)}] If $\vph(\theta)>0$ for some $\theta\in L_2$, then $(\BB\vph)(\wt \theta)>0$ for $\wt\theta=(0,0,0)$,
that is, for $\wt\theta\in L_3$; and
\item[{(5$'$)}] if $\vph(\theta)>0$ for every $\theta\in L_1\setminus L_2$,
then $(\BB\vph)(\wt \theta)>0$ for every $\wt\theta\in L_2\setminus L_3$.
\end{itemize}
If one accepts the above facts then it is immediate from~(1)--(5) that if $\vph\in C_V^+\setminus\{0\}$
then $(\BB^5\vph)(\theta)>0$ for every $\theta\in L_0\setminus L_2$. (Note that if $\vph\in C_V^+\setminus\{0\}$
then necessarily $\vph(\theta)>0$ for some $\theta\in L_0$.)
It is also immediate from~(1)--(4) and~(5$'$) that if $\vph\in C_V^+\setminus\{0\}$
then $(\BB^5)(\theta)>0$ for every $\theta\in L_2\setminus L_3$ and thus for every
$\theta\in (L_0\setminus L_2)\cup(L_2\setminus L_3)=L_0\setminus L_3=L_0\setminus\{(0,0,0)\}$.
Finally, one sees that if $\vph\in C_V^+\setminus\{0\}$ then also $\BB^2\vph\in C_V^+\setminus\{0\}$,
and using~(1), (2), and~(3$'$) one concludes that $(\BB^5\vph)(\theta)>0$ at $\theta=(0,0,0)$.
One therefore concludes that if $\vph\in C_V^+\setminus\{0\}$, then $(\BB^5)(\theta)>0$ for every $\theta\in L_0=T_3\setminus\{(-1,0,0)\}$.

Now fixing any $\vph\in C_V^+\setminus\{0\}$, write
$$
(\BB^5\vph)(\theta)=Q_0\nu_0(\theta)+Q_1\nu_1(\theta)+\psi(\theta)
$$
with $\psi\in C_0(T_3)$ as per Theorem~8.2. Then
$$
Q_0=\frac{1}{2}\int_{-1}^0t^2(1+t)[(\BB^4\vph)(t,0,0)]\:dt,\qquad
Q_1=\int_{-1}^0t^2(1+t)[(\BB^4\vph)(-1,t,0)]\:dt.
$$
As $\BB^2\vph\in C_V^+\setminus\{0\}$, one has from~(1) and~(2) that $(\BB^4\vph)(\theta)>0$
for some $\theta\in L_2$, that is, for some $\theta=(t,0,0)$ with $t\in (-1,0]$. Thus $Q_0>0$.
Similarly, as $\BB\vph\in C_V^+\setminus\{0\}$, one has from~(1)--(3) that there exists $\del>0$
such that $(\BB^4\vph)(\theta)>0$ for every $\theta\in L_\del$, and in particular for every
$\theta=(-1,t,0)$ with $t\in[-\del,0)$. Thus $Q_1>0$. With this,
the existence of a uniform lower bound $B_5$ for $\BB^5\vph$ follows directly from Lemma~8.9.

There remains to establish the properties~(1)--(5) and~(3$'$) and~(5$'$). For the most part, these
follow rather straightforwardly from the formulas~\eqref{bm3}, \eqref{phiker} for $\BB_0\vph$
and $\BB_1\vph$. Let us rewrite~\eqref{bm3} as
\be{bm33}
\begin{array}{lcl}
(\BB_0\vph)(\wt\theta) &\bb = &\bb
\ds{\tint_{\wt\theta_{1}}^{\wt\theta_{2}}\tint_{\wt\theta_{2}}^{\wt\theta_{3}}
\Phi_0(t_1,t_2,\wt\theta_1,\wt\theta_2,\wt\theta_3)
\vph(t_{1},t_{2},0)\:dt_{2}\:dt_{1},}\\
\\
(\BB_1\vph)(\wt\theta) &\bb = &\bb
\ds{\nu_1(\wt\theta)\:\:\tint_{-1}^{\wt\theta_{1}}\tint_{\wt\theta_{1}}^{\wt\theta_{2}}\tint_{\wt\theta_{2}}^{\wt\theta_{3}}
\Phi_1(t_0,t_1,t_2,\wt\theta_1,\wt\theta_3)
\vph(t_{0},t_{1},t_{2})\:dt_{2}\:dt_1\:dt_{0},}
\end{array}
\ee{bm33}
using the variable $\wt\theta=(\wt\theta_1,\wt\theta_2,\wt\theta_3)$. We recall that the arguments of these
integrals are nonnegative throughout the range of integration, and so it is enough to prove for each
of the indicated $\wt\theta$ in the above claimed properties, that either $(\BB_0\vph)(\wt\theta)>0$
or $(\BB_1\vph)(\wt\theta)>0$. Generally, this will be done by exhibiting a point $(t_1,t_2)$ or $(t_0,t_1,t_2)$
in the range of integration at which the integrand is strictly positive. As we are taking average integrals, it
will not matter if the upper and lower limits of an integral are equal.

To prove~(1), we assume that $\vph(\theta)>0$ for some $\theta=(\theta_1,\theta_2,\theta_3)\in L_0$,
and without loss we may assume that
\be{spaced}
-1<\theta_1<\theta_2<\theta_3\le 0
\ee{spaced}
as $\vph$ is continuous on $L_0$.
Now letting $\wt\theta=(\wt\theta_1,\wt\theta_2,\wt\theta_3)=(\theta_1,\theta_2,0)\in L_1$,
one sees directly that $(\BB_1\vph)(\wt\theta)>0$. In particular, the relevant integrand in~\eqref{bm33} is strictly
positive at the point $(t_0,t_1,t_2)=(\theta_1,\theta_2,\theta_3)$
which lies within the range of integration, as we have from~\eqref{phiker} that
$$
\Phi_1(\theta_1,\theta_2,\theta_3,\theta_1,0)\vph(\theta_1,\theta_2,\theta_3)
=(\theta_2-\theta_1)\bigg(\frac{\theta_3-\theta_2}{-\theta_1}\bigg)(\theta_3-\theta_1)(1+\theta_1-\theta_2)\vph(\theta_1,\theta_2,\theta_3)>0.
$$
Observe that the assumptions~\eqref{spaced} are used in drawing this conclusion.
Additionally, $\nu_1(\wt\theta)>0$ as $\theta_1>-1$. With this~(1) is established.

The proof of~(2) is similar. We assume that $\vph(\theta)>0$ for some $\theta=(\theta_1,\theta_2,0)\in L_1$,
and without loss $-1<\theta_1<\theta_2<\theta_3=0$. Letting $\wt\theta=(\theta_1,0,0)\in L_2$,
one sees that the integrand of the second integral in~\eqref{bm33} is positive at $(t_0,t_1,t_2)=(\theta_1,\theta_2,0)$
and also $\nu_1(\wt\theta)>0$ as before. Thus again $(\BB_1\vph)(\wt\theta)>0$, and~(2) is proved.

The proof of~(3) is slightly different from the proofs of~(1) and~(2). First, assuming that $\vph(\theta)>0$
for some $\theta\in L_2$, we may assume that $\theta=(\theta_1,0,0)$ where $-1<\theta_1<0$. Further, by
continuity, there exists $\del>0$ such that $\vph(\theta_1,\gam,0)>0$ for every $\gam\in[-\del,0]$,
and where also $\theta_1<-\del$. Now let any point $\wt\theta\in L_\del$ be given, that is,
$\wt\theta=(\wt\theta_1,\wt\theta_2,0)$ where $-1\le\wt\theta_1\le\wt\theta_2<0$ and also $-\del\le\wt\theta_2<0$.
Two cases now arise. First, suppose that $\wt\theta_1\ge\theta_1$. Then as in the proofs of~(1) and~(2), one
shows that $(\BB_1\vph)(\wt\theta)>0$ by noting that the point $(t_0,t_1,t_2)=(\theta_1,\wt\theta_2,0)$
lies in the domain of integration and the relevant integrand is positive there, and again that $\nu_1(\wt\theta)>0$.
The fact that
\be{thinq}
\theta_1<-\del\le\wt\theta_2<0,
\ee{thinq}
in particular, is used here.
For the second case we assume that $\wt\theta_1\le\theta_1$, and here we show that $(\BB_0\vph)(\wt\theta)>0$.
Indeed, the relevant integrand is strictly positive at $(t_1,t_2)=(\theta_1,\wt\theta_2)$, again because of~\eqref{thinq}.
This establishes~(3).

To prove~(3$'$), under the same assumptions as for~(3), one uses the second equation in~\eqref{btrip} with $\theta_*=0$
to show that $(\wt\BB_1\vph)(0,0,0)>0$ and so $(\BB_1\vph)(0,0,0)>0$.

To prove~(4) we assume that $\vph(\theta)>0$ for every $\theta\in L_\delta$, for some $\delta>0$ and let
$\wt\theta=(\wt\theta_1,\wt\theta_2,0)\in L_1\setminus L_2$, and so $-1\le\wt\theta_1\le\wt\theta_2<0$. We claim
that $(\BB_0\vph)(\wt\theta)>0$. To see this, let $(t_1,t_2)$ be such that
$\wt\theta_1\le t_1\le\wt\theta_2<t_2<0$ and also $t_2\in[-\delta,0)$. Then the integrand in the first
formula of~\eqref{bm33} is positive at $(t_1,t_2)$, establishing the claim.

To prove~(5), assume that $\vph(\theta)>0$ for every $\theta\in L_1\setminus L_2$ and let
$\wt\theta=(\wt\theta_1,\wt\theta_2,\wt\theta_3)\in L_0\setminus L_2$. Thus
$-1\le\wt\theta_1\le\wt\theta_2\le\wt\theta_3\le 0$ with $\wt\theta_2<0$. Two cases must be considered.
First, if $\wt\theta_1=\wt\theta_2=\wt\theta_3=\theta_*$ for some $\theta_*$, then $(\BB_0\vph)(\wt\theta)>0$
by the first equation in~\eqref{btrip}, since $(\theta_*,\theta_*,0)\in L_1\setminus L_2$. For the
second case, where either $\wt\theta_1<\wt\theta_2$ or $\wt\theta_2<\wt\theta_3$, let $(t_1,t_2)$ be such that
$\wt\theta_1\le t_1\le\wt\theta_2\le t_2\le\wt\theta_3$ and also $t_2<t_1<0$.  Then again the integrand
in the first formula of~\eqref{bm33} is positive at $(t_1,t_2)$, so again $(\BB_0\vph)(\wt\theta)>0$.

Finally, to prove~(5$'$), under the same assumptions as for~(5), let $\wt\theta=(\wt\theta_1,0,0)\in L_2\setminus L_3$,
and so $-1<\wt\theta_1<0$. Then the integrand in the second formula of~\eqref{bm33} is positive at the point
$(t_0,t_1,t_2)=(-1,\wt\theta_1,0)$ and also $\nu_1(\wt\theta)>0$, and so $(\BB_1\vph)(\wt\theta)>0$.
With this the proof of the proposition is complete.~\fp
\vv

\noindent {\bf Proof of Theorem~7.4.}
Given any $\vph\in C(T_m)^+\setminus\{0\}$, then by~\eqref{b3} of Proposition~7.6
we have the upper bound $\AA^k\vph\le 2^k\|\vph\|u_m$
(the order here being with respect to the cone $C(T_m)^+$) for every $k\ge m-1$, and
in particular for every $k\ge 3$
if $m=2$ and for every $k\ge 5$ if $m=3$.

To obtain a lower bound for $\AA^k\vph$, fix $\psi\in C(T_m)^+\setminus\{0\}$ satisfying
$\psi\le\vph$ and such that the support of the function $\psi$ is contained in
the set $O_m$ in~\eqref{om}. Then upon defining
$$
\zeta(\theta)=\frac{\psi(\theta)}{u_m(\theta)}\qquad\hbox{for }\theta\in T_m,
$$
we have that $\zeta\in C(T_m)^+\setminus\{0\}$. Moreover, we have directly from
the formulas~\eqref{aa} and~\eqref{bformx} defining $\AA$ and $\BB$ that
$$
(\BB^k\zeta)(\theta)=\frac{(\AA^k\psi)(\theta)}{u_m(\theta)}\le \frac{(\AA^k\vph)(\theta)}{u_m(\theta)}
$$
for every $k\ge 1$ and for $\theta\in O_m$.
From this it follows, by Proposition~8.1 in the case $m=2$, and by Proposition~8.3
in the case $m=3$, that we have the lower bounds
\be{ab}
B_k u_m(\theta)\le(\AA^k\vph)(\theta)
\ee{ab}
for every $k\ge 3$ if $m=2$ and for every $k\ge 5$ if $m=3$.
The bounds~\eqref{ab} are valid for $\theta\in T_2$ if $m=2$, or for
$\theta\in T_3\setminus\{(-1,0,0)\}$ if $m=3$ and thus for every $\theta\in T_3$
as $\AA^k\vph$ is continuous in $T_3$.

We conclude that $\AA^k\vph\sim u_m$ for every $k\ge 3$ if $m=2$
and for every $k\ge 5$ if $m=2$, as desired.~\fp
\vv

\noindent {\bf Proof of Theorem~7.1.}
This follows directly
from Proposition~7.3 and Theorem~7.4.~\fp

\bigskip

\noindent {\small
\begin{tabular}{l}
John Mallet-Paret\\
Division of Applied Mathematics\\
Brown University\\
Providence, Rhode Island 02912\\
{\tt John\_Mallet-Paret@brown.edu}\\
\\
Roger D. Nussbaum\\
Department of Mathematics\\
Rutgers University\\
Piscataway, New Jersey 08854\\
{\tt nussbaum@math.rutgers.edu}
\end{tabular}}

\end{document}